\documentclass{amsart}

\usepackage{amssymb}

\allowdisplaybreaks

\newtheorem{theo-intro}{Theorem}
\newtheorem{theorem}{Theorem}[section]
\newtheorem{theoA}{Theorem}
\newtheorem{theoB}{Theorem}
\newtheorem{theoC}{Theorem}
\newtheorem{theoD}{Theorem}
\newtheorem{theoE}{Theorem}

\newtheorem{theoF}{Theorem}
\newtheorem{lemma}[theorem]{Lemma}
\newtheorem{proposition}[theorem]{Proposition}
\newtheorem{corollary}[theorem]{Corollary}
\newtheorem{remark}[theorem]{Remark}

\theoremstyle{definition}

\newcommand{\R}{\mathbb{R}}
\newcommand{\C}{\mathbb{C}}

\newcommand{\RR}{\mathcal{R}}
\newcommand{\LL}{\mathcal{L}}
\newcommand{\QQ}{\mathcal{Q}}
\newcommand{\Al}{\mathcal{A}}
\newcommand{\Be}{\mathcal{B}}

\newcommand{\Le}{\mathsf{L}}

\newcommand{\Es}{\mathsf{E}}

\newcommand{\ten}{\otimes}
\newcommand{\al}{\alpha}

\newcommand{\dem}{\noindent {\bf Proof. }}
\newcommand{\demA}{\noindent {\bf Proof of Theorem A. }}
\newcommand{\demB}{\noindent {\bf Proof of Theorem B. }}
\newcommand{\demC}{\noindent {\bf Proof of Theorem C. }}
\newcommand{\demD}{\noindent {\bf Proof of Theorem D. }}
\newcommand{\demE}{\noindent {\bf Proof of Theorem E. }}
\newcommand{\demF}{\noindent {\bf Sketch of the proof of Theorem F. }}

\newcommand{\fin}{\hspace*{\fill} $\square$ \vskip0.2cm}

\newcommand{\summ}{\sum\nolimits}

\def\bubl{{\displaystyle \mathop{\mathsf{A}}^{\circ}}}
\def\bubla{{\displaystyle \mathop{a}^{\circ}}}

\begin{document}

\title[Rosenthal type inequalities for free chaos]
{Rosenthal type inequalities for free chaos}

\author[Junge, Parcet and Xu]
{Marius Junge, Javier Parcet and Quanhua Xu}

\address{University of Illinois at Urbana-Champaign}

\email{junge@math.uiuc.edu}

\address{Centre de Recerca Matem{\`a}tica, Barcelona}

\email{jparcet@crm.es}

\address{Universit{\'e} Franche-Comt{\'e}, Besan\c{c}on}

\email{qx@math.univ-fcomte.fr}


\begin{abstract}
Let $\mathcal{A}$ denote the reduced amalgamated free product of a
family $\mathsf{A}_1, \mathsf{A}_2, \ldots, \mathsf{A}_n$ of von
Neumann algebras over a von Neumann subalgebra $\Be$ with respect
to normal faithful conditional expectations $\Es_k: \mathsf{A}_k
\to \Be$. We investigate the norm in $L_p(\Al)$ of homogeneous
polynomials of a given degree $d$. We first generalize
Voiculescu's inequality to arbitrary degree $d \ge 1$ and indices
$1 \le p \le \infty$. This can be regarded as a free analogue of
the classical Rosenthal inequality. Our second result is a
length-reduction formula from which we generalize recent results
of Pisier, Ricard and the authors. All constants in our estimates
are independent of $n$ so that we may consider infinitely many
free factors. As applications, we study square functions of free
martingales. More precisely we show that, in contrast with the
Khintchine and Rosenthal inequalities, the free analogue of the
Burkholder-Gundy inequalities does not hold on $L_\infty(\Al)$. At
the end of the paper we also consider Khintchine type inequalities
for Shlyakhtenko's generalized circular systems.
\end{abstract}

\maketitle

\

\section*{Introduction and main results}

A strong interplay between harmonic analysis, probability theory
and Banach space geometry can be found in the works of Burkholder,
Gundy, Kwapie\'n, Maurey, Pisier, Rosenthal and many others
carried out mostly in the 70's. Norm estimates for sums of
independent random variables as well as martingale inequalities
play a prominent role. Let us mention, for instance, the classical
Khintchine and Rosenthal inequalities, Fefferman's duality theorem
and the inequalities of Burkholder and Burkholder-Gundy for
martingales. On the other hand, in the last two decades the
noncommutative analogues of these aspects have been considerably
developed. Important tools in this process come from free
probability, operator space theory and theory of noncommutative
martingales.

\vskip5pt

\renewcommand{\theequation}{$\mathsf{R}_p$}
\addtocounter{equation}{-1}

In this paper, we continue this line of research by studying
$L_p$-estimates for homogeneous polynomials of free random
variables. Our results are motivated by the classical Rosenthal
inequality \cite{Ro}. That is, given a family $f_1, f_2, f_3,
\ldots$ of independent, mean-zero random variables over a
probability space $\Omega$, we have
\begin{equation} \label{Rosenthal}
\Big\| \sum_{k=1}^n f_k \Big\|_{L_p(\Omega)} \sim_{c_p} \Big(
\sum_{k=1}^n \|f_k\|_2^2 \Big)^{\frac12} + \Big( \sum_{k=1}^n
\|f_k\|_p^p \Big)^{\frac1p},
\end{equation}
for $2 \le p < \infty$ and where $\mathsf{A} \sim_c \mathsf{B}$
means that $c^{-1} \mathsf{A} \le \mathsf{B} \le c \mathsf{A}$.
The growth rate for the constant $c_p$ as $p \to \infty$ is
$p/\log p$ (see \cite{JSZ}) and so the Rosenthal inequality fails
on $L_\infty(\Omega)$. In sharp contrast are Voiculescu's
inequality \cite{V2} and its operator-valued analogue \cite{J2}
which are valid in $L_\infty$. Let $\Al = \mathsf{A}_1 *
\mathsf{A}_2 * \cdots * \mathsf{A}_n$ denote the reduced free
product of a family $\mathsf{A}_1, \mathsf{A}_2, \ldots,
\mathsf{A}_n$ of von Neumann algebras equipped with normal
faithful (\emph{n.f.} in short) states $\phi_1, \phi_2, \ldots,
\phi_n$, respectively. Then, given $a_1 \in \mathsf{A}_1, a_2 \in
\mathsf{A}_2, \ldots, a_n \in \mathsf{A}_n$ mean-zero random
variables (i.e. freely independent) in $\Al$ and a collection
$b_1, b_2, \ldots, b_n \in \Be(\mathcal{H})$ of bounded linear
operators on some Hilbert space $\mathcal{H}$, Voiculescu's
inequality claims that
\renewcommand{\theequation}{$\mathsf{V}_{\infty}$}
\addtocounter{equation}{-1}
\begin{eqnarray} \label{Voiculescu}
\Big\| \sum_{k=1}^n a_k \ten b_k \Big\|_{\Al \bar\ten
\Be(\mathcal{H})} & \sim_c & \sup_{1\le k \le n} \big\| a_k \ten
b_k \big\|_{\mathsf{A}_k \bar\ten \Be(\mathcal{H})} \\ \nonumber &
+ & \Big\| \Big( \sum_{k=1}^n  \phi_k(a_k^*a_k)b_k^* b_k
\Big)^{\frac12} \Big\|_{\Be(\mathcal{H})}
\\ \nonumber & + & \Big\| \Big( \sum_{k=1}^n  \phi_k(a_k a_k^*)b_k b_k^*
\Big)^{\frac12} \Big\|_{\Be(\mathcal{H})}
\end{eqnarray}
for some universal positive constant $c$. The equivalence
\eqref{Voiculescu} was proved by Voiculescu \cite{V2} in the
tracial scalar-valued case. The general case as stated above (or
more generally using amalgamated free product) can be found in
\cite{J2}. This result can be regarded as the operator-valued free
analogue of the Rosenthal inequality for homogeneous free
polynomials of degree $1$ and $p=\infty$. Quite surprisingly, the
$L_\infty$-estimates (which do not hold in the classical case) are
easier to obtain in the free case by virtue of the Fock space
representation. In contrast with the classical situation, the
passage from $L_\infty$ to $L_p$ in the free setting is much more
delicate. This is mainly because of the fact that a concrete Fock
space representation does not seem available for $L_p({\mathcal
A})$.

\vskip5pt

Our first contribution in this paper consists in generalizing
Voiculescu's inequality to homogeneous free polynomials of
arbitrary degree $d$ and to any index $1 \le p \le  \infty$. Let
us be more precise and fix some notations. Assume that $\Be$ is a
common von Neumann subalgebra of $\mathsf{A}_1, \mathsf{A}_2,
\ldots, \mathsf{A}_n$ such that there is a normal faithful
conditional expectation $\Es_k: \mathsf{A}_k \to \Be$ for each
$k$. Let $\Al$ be the reduced amalgamated free product $*_\Be
\mathsf{A}_k$ of $\mathsf{A}_1, \mathsf{A}_2, \ldots,
\mathsf{A}_n$ over $\mathcal{B}$ with respect to the $\Es_k$.
$\Es: \Al \to \Be$ will denote the corresponding conditional
expectation and $\mathbf{P}_\Al(p,d)$ the subspace of $L_p(\Al)$
of homogeneous free polynomials of degree $d$. Then, given $1 \le
k \le n$, we consider the map $\QQ_k$ on $\mathbf{P}_\Al(p,d)$
which collects all reduced words starting and ending with a letter
in $\mathsf{A}_k$. Then we have the following result.

\begin{theoA} \label{Theorem-Voiculescu}
If $2 \le p \le \infty$ and $a_1, a_2, \ldots, a_n \in
\mathbf{P}_\Al(p,d)$, we have
\begin{eqnarray*}
\Big\| \sum_{k=1}^n \QQ_k(a_k) \Big\|_p & \sim_{c d^7} & \Big(
\sum_{k=1}^n \big\| \QQ_k(a_k) \big\|_p^p \Big)^{\frac1p} \\ & + &
\Big\| \Big( \sum_{k=1}^n \mathsf{E} \big( \QQ_k(a_k)^* \QQ_k(a_k)
\big) \Big)^{\frac12}\Big\|_p \\ & + & \Big\| \Big( \sum_{k=1}^n
\Es \big( \QQ_k(a_k) \QQ_k(a_k)^* \big) \Big)^{\frac12} \Big\|_p.
\end{eqnarray*}
\end{theoA}

We note that the operator-valued case is also contemplated in
Theorem A since we are allowing amalgamation, see Remark
\ref{Remark-Operator-Valued} below for more details. On the other
hand, we also point out that, since freeness implies
noncommutative independence, the case of degree $1$ polynomials
for $2 \le p < \infty$ follows from the noncommutative analogue of
Rosenthal's inequality \cite{JX,JX2}. However, the constants
obtained in this way are not uniformly bounded as $p \to \infty$,
see Remark \ref{Remark-Free-Independence} for a more detailed
discussion. Finally, we should also emphasize that Theorem A can
be easily generalized to the case $1 \le p \le 2$ by duality, see
Remark \ref{Remark-Duality-AC} for the details.

\vskip5pt

Our second major result is a length-reduction formula for
homogeneous free polynomials in $L_p(\Al)$. Again, we need to fix
some notations. In what follows, $\Lambda$ will denote a finite
index set and we shall keep the terminology for $\Al$, $\Be$ and
$\Es: \Al \to \Be$. Then, we use the following notation suggested
by quantum mechanics
\begin{eqnarray*}
\Big\| \sum_{\alpha \in \Lambda} b(\alpha) \big\langle a(\alpha) |
\Big\|_p & = & \Big\| \Big( \sum_{\alpha,\beta \in \Lambda}
b(\alpha) \Es \big( a(\alpha) a(\beta)^* \big) b(\beta)^*
\Big)^{\frac12} \Big\|_p,
\\ \Big\| \sum_{\alpha \in \Lambda} | a(\alpha) \big\rangle b(\alpha)
\Big\|_p & = & \Big\| \Big( \sum_{\alpha,\beta \in \Lambda}
b(\alpha)^* \Es \big( a(\alpha)^* a(\beta) \big) b(\beta)
\Big)^{\frac12} \Big\|_p.
\end{eqnarray*}
Finally, given $1 \le k \le n$ we consider the map $\LL_k$ (resp.
$\RR_k$) on $\mathbf{P}_\Al(p,d)$ which collects the reduced words
starting (resp. ending) with a letter in $\mathsf{A}_k$. Thus we
have $$\QQ_k = \LL_k \RR_k = \RR_k \LL_k.$$ We shall write
$\mathbf{P}_\Al(d)$ for $\mathbf{P}_\Al(p,d)$ with $p=\infty$. Our
second result is the following.

\begin{theoB} \label{Theorem-Reduction}
Let $2 \le p \le \infty$ and let $x_k(\alpha) \in
L_p(\mathsf{A}_k)$ with $\Es(x_k(\alpha)) = 0$ for each $1 \le k
\le n$ and $\alpha$ running over a finite set $\Lambda$. Let
$w_k(\alpha) \in \mathbf{P}_\Al(d)$ for some $d \ge 0$ and
satisfying $\RR_k(w_k(\alpha)) = 0$ for all $1 \le k \le n$ and
every $\alpha \in \Lambda$. Then, we have the equivalence
 $$\Big\|
 \sum_{k, \alpha} w_k(\alpha) x_k(\alpha) \Big\|_{L_p(\Al)} \sim_{c
 d^{2}} \Big\| \sum_{k, \alpha} w_k(\alpha) \big\langle x_k(\alpha)
 | \Big\|_p + \Big\| \sum_{k, \alpha} | w_k(\alpha) \big\rangle
 x_k(\alpha) \Big\|_p.$$
Similarly, if $\LL_k(w_k(\alpha)) = 0$ we have
 $$\Big\| \sum_{k,
 \alpha} x_k(\alpha) w_k(\alpha) \Big\|_{L_p(\Al)} \sim_{c d^{2}}
 \Big\| \sum_{k, \alpha} | x_k(\alpha) \big\rangle w_k(\alpha)
 \Big\|_p + \Big\| \sum_{k, \alpha} x_k(\alpha) \big\langle
 w_k(\alpha) | \Big\|_p.$$
\end{theoB}

\vskip5pt

A large part of this paper will be devoted to the proofs of
Theorems A and B. One of the key points in both proofs is the main
complementation result in \cite{RX} (\emph{c.f.} Theorem
\ref{Theorem-Complementation} below) since it allows us to use
interpolation starting from the case $p=\infty$, for which  both
results hold with constants independent of $d$. Our main
application of Theorem B is a Khintchine type inequality. In the
classical case, Khintchine's inequality is a particular case of
Rosenthal's inequality with relevant constant $c_p \sim \sqrt{p}$
as $p \to \infty$. However, as in the Rosenthal/Voiculescu case,
the free analogue of Khintchine's inequality holds in $L_\infty$.
Indeed, the first example of this phenomenon was found by Leinert
\cite{Le}, who replaced the Bernoulli random variables by the
operators $\lambda(g_1), \lambda(g_2), \ldots, \lambda(g_n)$
arising from the generators $g_1, g_2, \ldots, g_n$ of a free
group $\mathbb{F}_n$ via the left regular representation
$\lambda$. More generally, if $\mathcal{W}_d$ denotes the subset
of reduced words in $\mathbb{F}_n$ of length $d$ and
$\mathrm{C}^*_\lambda(\mathbb{F}_n)$ stands for the reduced
$\mathrm{C}^*$-algebra on $\mathbb{F}_n$, Haagerup \cite{H0}
proved that
\renewcommand{\theequation}{$\mathsf{H}_p$}
\addtocounter{equation}{-1}
\begin{equation} \label{Haagerup}
\Big\| \sum_{w \in \mathcal{W}_d} \alpha_w \lambda(w)
\Big\|_{\mathrm{C}^*_\lambda(\mathbb{F}_n)} \sim_{1+d} \Big(
\sum_{w \in \mathcal{W}_d} |\alpha_w|^2 \Big)^{\frac12}.
\end{equation}

There are two ways to extend these inequalities. The \emph{first
step} consists in considering operator-valued coefficients. In the
classical case, the operator-valued analogue is the so-called
noncommutative Khintchine inequality, by Lust-Piquard and Pisier
\cite{Lu,LuP}. Leinert's result was extended to the
operator-valued case by Haagerup and Pisier in \cite{HP} while
Haagerup's inequality (\ref{Haagerup}) was generalized by Buchholz
\cite{Bu1}. Finally, the result in \cite{Bu1} has been recently
extended to arbitrary indices $1 \le p \le \infty$ by Pisier and
the second-named author in \cite{PP}.

\vskip5pt

The \emph{second step} consists in replacing the free generators
by arbitrary free random variables and
$\mathrm{C}^*_\lambda(\mathbb{F}_n)$ by a reduced amalgamated free
product von Neumann algebra $\Al$. In this case we find the recent
paper \cite{RX} by Ricard and the third-named author, where
Buchholz's result was extended to arbitrary reduced amalgamated
free products, see also Buchholz \cite{Bu2} and Nou \cite{N} for
the case of $q$-gaussians.

\vskip5pt

In this paper, we shall apply Theorem B to generalize the main
results of \cite{PP,RX}. More precisely, we have the following
result.

\begin{theoC}
Let $x$ be a $d$-homogeneous free polynomial
 $$x = \sum_{\alpha \in \Lambda}^{\null} \sum_{j_1 \neq
 j_2 \neq \cdots \neq j_d} x_{j_1}(\alpha) \cdots x_{j_d}(\alpha)
 \in L_p(\Al)$$
for some $2 \le p \le \infty$. Then we have
 $$\|x\|_p \sim_{c^d d!^2} \Sigma_1 + \Sigma_2$$
where $\Sigma_1$ is given by
 $$\sum_{s=0}^d \Big\| \sum_{\alpha
 \in \Lambda} \sum_{j_1 \neq \cdots \neq j_d} | x_{j_1}(\alpha)
 \cdots x_{j_s}(\alpha) \big\rangle \big\langle x_{j_{s+1}}(\alpha)
 \cdots x_{j_d}(\alpha) | \Big\|_p,$$
and $\Sigma_2$ has the form
 $$\sum_{s=1}^d \Big( \sum_{j_s=1}^n \Big\| \sum_{\alpha \in
 \Lambda} \sum_{\begin{subarray}{c} 1 \le j_1 \neq \cdots \neq
 j_{s-1} \le n
 \\ 1 \le j_{s+1} \neq \cdots \neq j_{d} \le n \\ j_{s-1} \neq j_s
 \neq j_{s+1} \end{subarray}} | x_{j_1}(\alpha) \cdots
 x_{j_{s-1}}(\alpha) \big\rangle x_{j_s}(\alpha) \big\langle
 x_{j_{s+1}}(\alpha) \cdots x_{j_d}(\alpha) | \Big\|_p^p
 \Big)^{\frac1p}_.$$
\end{theoC}

Our proof of Theorem C is an inductive application of Theorem B
and provides a natural explanation of the norms $\Sigma_1$ and
$\Sigma_2$. This leads naturally to $3$ terms if $d=1$, $5$ terms
if $d=2$, etc... We refer to Section \ref{Section3} below for a
more detailed explanation of the norms $\Sigma_1$ and $\Sigma_2$.
In the case of $p=\infty$, this result was obtained in \cite{RX}
in a slightly different form but with an equivalence constant
depending linearly on the degree $d$, which is essential for the
applications there. The inductive nature of our arguments leads to
worse constants, see Remark \ref{Remark-Better-Constant} for
details.

\renewcommand{\theequation}{\arabic{equation}}

\vskip5pt

In the last part of the paper, we shall apply our techniques to
studying square functions of free martingales and Khintchine type
inequalities for generalized circular systems. More precisely, we
first study the free analogue of the Burkholder-Gundy inequalities
\cite{BG}. The noncommutative version of these inequalities was
obtained by Pisier and the third-named author in \cite{PX1}. Thus,
since any free martingale is a noncommutative martingale, the only
interesting case seems to be $p=\infty$. In contrast with the free
Khintchine and Rosenthal inequalities, we shall prove that the
free analogue of the Burkholder-Gundy inequalities does not hold
in $L_\infty(\Al)$. To be more precise, let us consider an
infinite family $\mathsf{A}_1, \mathsf{A}_2, \mathsf{A}_3, \ldots$
of von Neumann algebras equipped with distinguished normal
faithful states and the associated reduced free product $\Al = *_k
\mathsf{A}_k$. We consider the natural filtration
 $$\Al_n = \mathsf{A}_1 * \mathsf{A}_2 * \cdots *
 \mathsf{A}_n \qquad \mbox{with consitional expectation} \;
 \Es_n: \Al \to \Al_n.$$
Any martingale adapted to this filtration is called a \emph{free
martingale}. Now, let $\mathcal{K}_n$ be the best constant for
which the lower estimate below holds for all free martingales
$x_1, x_2, \ldots$ in $L_\infty(\Al)$
 $$\max \left\{ \Big\| \Big( \sum_{k=1}^{2n} d x_k d x_k^*
 \Big)^{\frac12} \Big\|_\infty,\; \Big\| \Big( \sum_{k=1}^{2n} d
x_k^* d x_k \Big)^{\frac12} \Big\|_\infty \right\} \le
 \mathcal{K}_n \Big\| \sum_{k=1}^{2n} d x_k \Big\|_\infty.$$
Then we have the following result.

\begin{theoD}
$\mathcal{K}_n$ satisfies $\mathcal{K}_n \ge c \log n$ for some
absolute positive constant $c$.
\end{theoD}

The last section is devoted to Khintchine type inequalities for
Shlyakhtenko's generalized circular variables \cite{S} and Hiai's
generalized $q$-gaussians \cite{Hi}. In these particular cases,
the resulting inequalities are much nicer than those of Theorem C.
The Khintchine inequalities for $1$-homogeneous polynomials of
generalized gaussians were already proved in \cite{X2}, see
Theorem \ref{Th-OldE}  for an explicit formulation. We obtain here
its natural extension for Hiai's generalized $q$-gaussians.
Namely, let us consider a system of $q$-generalized circular
variables $gq_k = \lambda_k \ell_q(e_k) + \mu_k \ell_q^*(e_{-k})$
(see Section \ref{Section-Generalized-Circular} for precise
definitions) and let $\Gamma_q$ denote the von Neumann algebra
generated by these variables in the GNS-construction with respect
to the vacuum state $\phi_q(\cdot) = \langle \Omega,\; \cdot\,
\Omega \rangle_q$. Then, if $d_{\phi_q}$ denotes the density
associated to the state $\phi_q$, we have the following
inequalities for the $L_p$-variables
 $$gq_{k,p} = d_{\phi_q}^{\frac{1}{2p}}\, gq_k\,
 d_{\phi_q}^{\frac{1}{2p}}.$$

\begin{theoE} Let $\mathcal{N}$ be a von Neumann algebra and $1 \le p \le
\infty$. Let us consider a finite sequence $x_1, x_2, \ldots, x_n$
in $L_p(\mathcal{N})$. Then, the following equivalences hold up to
a constant $c_q$ depending only on $q$.
\begin{itemize}
\item[i)] If $1 \le p \le 2$, then
\begin{eqnarray*}
\lefteqn{\Big\| \sum_{k=1}^n x_k \ten gq_{k,p} \Big\|_p} \\ &
\sim_{c_q} & \inf_{x_k = a_k + b_k} \Big\| \Big( \sum_{k=1}^n
\lambda_k^{\frac{2}{p}} \mu_k^{\frac{2}{p'}} a_ka_k^*
\Big)^{\frac12} \Big\|_p + \Big\| \Big( \sum_{k=1}^n
\lambda_k^{\frac{2}{p'}} \mu_k^{\frac{2}{p}} b_k^*b_k
\Big)^{\frac12} \Big\|_p\;.
\end{eqnarray*}

\item[ii)] If $2\le p\le \infty$, then
\begin{eqnarray*}
\lefteqn{\Big\| \sum_{k=1}^n x_k \ten gq_{k,p} \Big\|_p} \\ &
\sim_{c_q} & \max \left\{ \Big\| \Big( \sum_{k=1}^n
\lambda_k^{\frac{2}{p}} \mu_k^{\frac{2}{p'}} x_k x_k^*
\Big)^{\frac12} \Big\|_p \, , \, \Big\| \Big( \sum_{k=1}^n
\lambda_k^{\frac{2}{p'}} \mu_k^{\frac{2}{p}} x_k^* x_k
\Big)^{\frac12} \Big\|_p \right\}.
\end{eqnarray*}
\end{itemize}
Moreover, if $\mathcal{G}q_{p}$ denotes the closed subspace of
$L_p(\Gamma_q)$ generated by the system of the generalized
$q$-gaussians $(gq_{k,p})_{k \ge 1}$, there exists a completely
bounded projection $\gamma q_p: L_p(\Gamma_q) \to \mathcal{G}q_p$
satisfying
 $$\|\gamma q_p\|_{cb} \le \Big( \frac{2}{\sqrt{1-|q|}}
 \Big)^{|1-\frac{2}{p}|}.$$
\end{theoE}

\vskip5pt

In our last result we calculate the Khintchine inequalities for
$2$-homogeneous polynomials of generalized free gaussians $g_k =
\lambda_k \ell(e_k) + \mu_k \ell^*(e_{-k})$ (corresponding to the
case of $q=0$). As we shall see, our method is also valid for
$d$-homogeneous polynomials and the resulting inequalities can be
regarded as asymmetric versions of the main inequalities in
\cite{PP}. Let $\Gamma$ denote the von Neumann algebra generated
by the system of $g_k$'s in the GNS-construction with respect to
the vacuum state $\phi(\cdot) = \langle \Omega,\; \cdot\, \Omega
\rangle$. Our result reads as follows.

\begin{theoF}
Let $\mathcal{N}$ be a von Neumann algebra and $1 \le p \le
\infty$. Let us consider a finite double indexed family $x =
(x_{ij})_{i,j \ge 1}$ in $L_p(\mathcal{N})$ and define the
following norms associated to $x$
\begin{eqnarray*}
\mathcal{M}_{p}(x) & = & \Big\| \sum_{i \neq j} (\mu_i
\lambda_j)^{\frac{1}{p}} (\lambda_i \mu_j)^{\frac{1}{p'}} x_{ij}
\ten e_{ij} \Big\|_{S_p(L_p(\mathcal{N}))},
\\ \mathcal{R}_{p}(x) & = & \Big\| \Big( \sum_{i \neq j}
(\mu_i \mu_j)^{\frac{2}{p'}} (\lambda_i \lambda_j)^{\frac{2}{p}}
x_{ij} x_{ij}^* \Big)^{\frac12} \Big\|_{L_p(\mathcal{N})}, \\
\mathcal{C}_{p}(x) & = & \Big\| \Big( \sum_{i \neq j} (\mu_i
\mu_j)^{\frac{2}{p}} (\lambda_i \lambda_j)^{\frac{2}{p'}} x_{ij}^*
x_{ij} \Big)^{\frac12} \Big\|_{L_p(\mathcal{N})}.
\end{eqnarray*}
Then, the following equivalences hold up to an absolute positive
constant $c$.
\begin{itemize}
\item[i)] If $1 \le p \le 2$, then $$\Big\| \sum_{i \neq j} x_{ij}
\ten d_\phi^{\frac{1}{2p}}\, g_i\, g_j\, d_\phi^{\frac{1}{2p}}
\Big\|_p \sim_c \inf_{x = a + b + c} \mathcal{R}_{p}(a) +
\mathcal{M}_{p}(b) + \mathcal{C}_{p}(c).$$

\item[ii)] If $2\le p\le \infty$, then
 $$\Big\| \sum_{i \neq j}
x_{ij} \ten d_\phi^{\frac{1}{2p}}\, g_i\, g_j\,
d_\phi^{\frac{1}{2p}} \Big\|_p \sim_c \ \max \ \Big\{
\mathcal{R}_{p}(x) \, , \, \mathcal{M}_{p}(x) \, , \,
\mathcal{C}_{p}(x) \Big\}.$$
\end{itemize}
Moreover, if $\mathcal{G}_{p,2}$ denotes the subspace of
$L_p(\Gamma)$ generated by the system $$\Big\{
d_\phi^{\frac{1}{2p}}\, g_i\, g_j\, d_\phi^{\frac{1}{2p}} \,
\big| \ 1  \le i \neq j < \infty \Big\},$$ there exists a
projection $\gamma_{p,2}: L_p(\Gamma) \to \mathcal{G}_{p,2}$ with
cb-norm uniformly bounded on $p$.
\end{theoF}

We conclude the Introduction with some general remarks. We shall
assume some familiarity with Voiculescu's free probability
\cite{V,V2,VDN} and Pisier's vector-valued noncommutative
integration \cite{P2}. In fact, we will be concerned only with the
vector-valued Schatten classes and their column/row subspaces. On
the other hand, since we are working over (amalgamated) free
product von Neumann algebras, we shall need to use Haagerup
noncommutative $L_p$-spaces \cite{H,T1}. As is well known,
Haagerup $L_p$-spaces have trivial intersection and thereby do not
form an interpolation scale. However, the complex interpolation
method will be a basic tool in this paper. This problem is solved
by means of Kosaki's definition of $L_p$-spaces, see \cite{Ko,T2}.
We also refer the reader to Chapter 1 in \cite{JP} or to the
survey \cite{PX2} for a quick review of Haagerup's and Kosaki's
definitions of noncommutative $L_p$-spaces and the compatibility
between them. In particular, using such compatibility, we shall
use in what follows the complex interpolation method without
further details. At some specific points, we shall also need some
basic notions from operator space theory \cite{ER,P3}, Hilbert
modules \cite{L} and Tomita's modular theory \cite{KR,PT}. Along
the paper, $c$ will denote an absolute positive constant that may
change from one instance to another.

\vskip8pt

\noindent \textbf{Acknowledgment.} We thank Ken Dykema, Gilles
Pisier, Ana Maria Popa and Eric Ricard for discussions related to
the content of this paper. The first-named author was partially
supported by the NSF DMS-0301116. The second-named author was
partially supported by the Project MTM2004-00678, Spain.

\section{Amalgamated free products} \label{Section1}

We begin by recalling the construction of the reduced amalgamated
free product of a family of von Neumann algebras. Amalgamated free
products of $\mathrm{C}^*$-algebras, which we also outline below,
were introduced by Voiculescu \cite{V}. Let $\mathsf{A}_1,
\mathsf{A}_2, \ldots, \mathsf{A}_n$ be a family of von Neumann
algebras and let $\mathcal{B}$ be a common von Neumann subalgebra
of all of them. We assume that there are normal faithful
conditional expectations $\mathsf{E}_k: \mathsf{A}_k \to
\mathcal{B}.$ In addition, we also assume the existence of a von
Neumann algebra $\mathcal{A}$ containing $\mathcal{B}$ with a
normal faithful conditional expectation $\mathsf{E}: \mathcal{A}
\to \mathcal{B}$ and the existence of $*$-homomorphisms $\pi_k:
\mathsf{A}_k \to \mathcal{A}$ such that
 $$\mathsf{E} \circ \pi_k =
 \mathsf{E}_k \quad \mbox{and} \quad {\pi_k}_{\mid_{\mathcal{B}}}
 =id_{\mathcal{B}}.$$
The family $\mathsf{A}_1, \mathsf{A}_2, \ldots, \mathsf{A}_n$ is
called \emph{freely independent} over $\mathsf{E}$ if
 $$\mathsf{E}
 \big( \pi_{j_1}(a_1) \pi_{j_2}(a_2) \cdots \pi_{j_m}(a_m) \big)
 =0$$
whenever $a_k\in \mathsf{A}_{j_k}$ are such that
$\mathsf{E}(\pi_{j_k}(a_k)) = 0$ for all $1 \le k \le m$ and $j_1
\neq j_2 \neq \cdots \neq j_m$. In what follows we may identify
$\mathsf{A}_k$ with the von Neumann subalgebra
$\pi_k(\mathsf{A}_k)$ of $\mathcal{A}$ with no risk of confusion.
In particular, we may use $\mathsf{E}$ or $\mathsf{E}_k$
indistinctively over $\mathsf{A}_k$. Moreover, for notational
convenience we shall only use $\mathsf{E}$ almost all the time. In
the scalar case, $\mathcal{B}$ is the complex field and the
conditional expectations $\mathsf{E}$ and $\mathsf{E}_1,
\mathsf{E}_2, \ldots, \mathsf{E}_n$ are replaced by normal
faithful states.

\vskip5pt

As in the scalar-valued case, operator-valued freeness admits a
natural Fock space representation. We first assume that
$\mathsf{A}_1, \mathsf{A}_2, \ldots, \mathsf{A}_n$ are
$\mathrm{C}^*$-algebras having $\mathcal{B}$ as a common
$\mathrm{C}^*$-subalgebra. Let us consider the mean-zero subspaces
 $$\bubl_k = \Big\{ a_k \in \mathsf{A}_k \, \big| \ \mathsf{E}(a_k)
 = 0 \Big\}.$$
We define the Hilbert $\mathcal{B}$-module
 $$\bubl_{j_1} \otimes_{\mathcal{B}} \bubl_{j_2}
 \otimes_{\mathcal{B}} \cdots \ten_{\mathcal{B}} \bubl_{j_m}$$
equipped with the $\mathcal{B}$-valued inner product
 $$\big\langle
 a_1 \otimes \cdots \otimes a_m,\; a'_1 \otimes \cdots \otimes a'_m
 \big\rangle = \mathsf{E}_{j_m} \big( a_m^* \cdots
 \mathsf{E}_{j_2}(a_2^* \, \mathsf{E}_{j_1}(a_1^* a'_1) \, a'_2)
 \cdots a'_m \big).$$
Then, the usual Fock space is replaced by the Hilbert
$\mathcal{B}$-module
 $$\mathcal{H}_{\mathcal{B}}
 = \mathcal{B} \oplus \bigoplus_{m \ge1}
 \bigoplus_{j_1 \neq j_2 \neq \cdots \neq j_m} \bubl_{j_1}
 \otimes_{\mathcal{B}} \bubl_{j_2} \otimes_{\mathcal{B}} \cdots
 \otimes_{\mathcal{B}} \bubl_{j_m}.$$
The direct sums above are assumed to be $\mathcal{B}$-orthogonal.
Let $\mathcal{L(H_B)}$ stand for the algebra of adjointable maps
on $\mathcal{H_B}$. Recall that a linear right
$\mathcal{B}$-module map $\mathrm{T}: \mathcal{H_B} \rightarrow
\mathcal{H_B}$ is called \emph{adjointable} if there exists
$\mathrm{S}: \mathcal{H_B} \rightarrow \mathcal{H_B}$ such that
 $$\langle x, \mathrm{T}y
 \rangle = \langle \mathrm{S} x,y \rangle \qquad \mbox{for all}
 \qquad x,y \in \mathcal{H_B}.$$
Let us also recall how elements of $\mathsf{A}_k$ act on
$\mathcal{H_B}$. We decompose any $a_k \in \mathsf{A}_k$ as
 $$a_k= \bubla_k + \mathsf{E}_k(a_k).$$
An element in $\mathcal{B}$ acts on $\mathcal{H_B}$ by left
multiplication. Therefore, it suffices to define the action of
mean-zero elements. The $*$-homomorphism $\pi_k: \mathsf{A}_k
\rightarrow \mathcal{L(H_B)}$ has the following form
 $$\pi_k(\bubla_k) (x_{j_1} \otimes \cdots \otimes x_{j_m}) =
 \left\{
 \begin{array}{ll} \, \bubla_k \otimes x_{j_1} \otimes x_{j_2} \ten
 \cdots \otimes x_{j_m}, & \mbox{if} \ k \neq j_1 \\ \\ \,
 \mathsf{E}(\bubla_k x_{j_1}) \, x_{j_2} \otimes \cdots \otimes
 x_{j_m} \ \oplus & \\ \big( \bubla_k x_{j_1} - \mathsf{E}(\bubla_k
 x_{j_1}) \big) \otimes x_{j_2} \otimes \cdots \otimes x_{j_m}, &
 \mbox{if} \ k = j_1.
 \end{array} \right.$$
This definition also applies for the empty word. Then, since the
algebra $\mathcal{L(H_B)}$ is a $\mathrm{C}^*$-algebra \cite{L},
we can define the \emph{reduced $\mathcal{B}$-amalgamated free
product} $\mathrm{C}^*(*_{\mathcal{B}} \mathsf{A}_k)$ as the
$\mathrm{C}^*$-closure of linear combinations of operators of the
form
 $$\pi_{j_1}(a_1) \pi_{j_2}(a_2) \cdots \pi_{j_m}(a_m).$$
The $\mathrm{C}^*$-algebra $\mathrm{C}^*(*_{\mathcal{B}}
\mathsf{A}_k)$ is usually denoted in the literature by
 $$*_k(\mathsf{A}_k,\mathsf{E}_k).$$
However, we shall use a more relaxed notation, see Remark
\ref{Remark-Notation-Free} below.

\vskip3pt

Now we assume that $\mathsf{A}_1, \mathsf{A}_2, \ldots,
\mathsf{A}_n$ and $\Be$ are von Neumann algebras and that
$\mathcal{B}$ comes equipped with a normal faithful state
$\varphi: \mathcal{B} \rightarrow \C$. This provides us with the
induced \emph{n.f.} states $\phi: \Al \rightarrow \C$ and $\phi_k:
\mathsf{A}_k \to \C$ given by
 $$\phi = \varphi \circ \mathsf{E} \quad \mbox{and} \quad
 \phi_k =\varphi \circ \mathsf{E}_k.$$
The Hilbert space
 $$L_2 \Big( \bubl_{j_1} \otimes_{\mathcal{B}}
 \bubl_{j_2} \otimes_{\mathcal{B}} \cdots \otimes_{\mathcal{B}}
 \bubl_{j_m}, \varphi \Big)$$
is obtained from $\bubl_{j_1} \otimes_{\mathcal{B}} \bubl_{j_2}
\otimes_{\mathcal{B}} \cdots \otimes_{\mathcal{B}} \bubl_{j_m}$ by
considering the inner product
 $$\big\langle a_1 \otimes \cdots
 \otimes a_m,\; a'_1 \otimes \cdots \otimes a'_m
 \big\rangle_{\varphi} = \varphi \Big( \big\langle a_1 \otimes
 \cdots \otimes a_m,\; a'_1 \otimes \cdots \otimes a'_m \big\rangle\Big).$$
Then we define the orthogonal direct sum
 $$\mathcal{H}_{\varphi} = L_2(\mathcal{B}) \oplus \bigoplus_{m \ge1}
 \bigoplus_{j_1 \neq j_2 \neq \cdots \neq j_m} L_2 \Big(
 \bubl_{j_1} \otimes_{\mathcal{B}} \bubl_{j_2}
 \otimes_{\mathcal{B}} \cdots \otimes_{\mathcal{B}} \bubl_{j_m},
 \varphi \Big).$$
Let us consider the $*$-representation $\lambda: \mathcal{L(H_B)}
\rightarrow \mathcal{B}(\mathcal{H}_{\varphi})$ defined by $\big(
\lambda(\mathrm{T}) x \big) = \mathrm{T} x$. The faithfulness of
$\lambda$ is implied by the fact that $\varphi$ is also faithful.
Indeed, assume that $\lambda(\mathrm{T}^* \mathrm{T}) = 0$, then
we have
 $$\langle \mathrm{T}^* \mathrm{T} x, x \rangle_{\varphi} =
 \varphi \big( \langle \mathrm{T} x, \mathrm{T} x \rangle \big) = 0
 \qquad \mbox{for all} \qquad x \in \mathcal{H_B}.$$
Since $\varphi$ is faithful, $\mathrm{T} x = 0$ (as an element in
$\mathcal{H_B}$) for all $x \in \mathcal{H_B}$, and so $\mathrm{T}
= 0$. Then, the \emph{$\mathcal{B}$-amalgamated reduced free
product} $*_{\mathcal{B}} \mathsf{A}_k$ is defined as the weak$^*$
closure of $\mathrm{C}^* (*_{\mathcal{B}} \mathsf{A}_k)$ in
$\mathcal{L(H_B)}$. After decomposing
 $$a_k = \bubla_k +\mathsf{E}(a_k)$$
and identifying $\bubl_k$ with $\lambda \big( \pi_k(\bubl_k)
\big)$, we can think of $*_{\mathcal{B}} \mathsf{A}_k$ as
 $$*_{\mathcal{B}} \mathsf{A}_k = \Big( \mathcal{B} \oplus
 \bigoplus_{m \ge 1} \bigoplus_{j_1 \neq j_2 \neq \cdots \neq j_m}
 \bubl_{j_1} \bubl_{j_2} \cdots \bubl_{j_m} \Big)''.$$
Again, the usual notation for $*_{\mathcal{B}} \mathsf{A}_k$ is a
bit more explicit one
 $\bar*_k (\mathsf{A}_k, \mathsf{E}_k).$

\vskip5pt

Let us consider the orthogonal projections
\begin{eqnarray*}
\mathcal{Q}_{\emptyset}: \mathcal{H}_{\varphi} & \rightarrow &
L_2(\mathcal{B}), \\ \mathcal{Q}_{j_1 \cdots j_m}:
\mathcal{H}_{\varphi} & \rightarrow & L_2 \Big( \bubl_{j_1}
\otimes_{\mathcal{B}} \bubl_{j_2} \otimes_{\mathcal{B}} \cdots
\otimes_{\mathcal{B}} \bubl_{j_m}, \varphi \Big).
\end{eqnarray*}
 Then $\Es: *_\Be \mathsf{A}_k \to \Be$ is given by $\Es(a) =
\QQ_\emptyset a \QQ_\emptyset$ and the mappings
 $$\mathcal{E}_{\mathsf{A}_k}:  *_{\mathcal{B}} \mathsf{A}_k\ni a
 \mapsto \mathcal{Q}_{\mathsf{A}_k} a \mathcal{Q}_{\mathsf{A}_k}
 \in \mathsf{A}_k \quad (\mathcal{Q}_{\mathsf{A}_k} =
 \mathcal{Q}_{\emptyset} + \mathcal{Q}_k),$$
are \emph{n.f.} conditional expectations. In particular, it turns
out that $\mathsf{A}_1, \mathsf{A}_2, \ldots, \mathsf{A}_n$ are
von Neumann subalgebras of $*_{\mathcal{B}} \mathsf{A}_k$ freely
independent over $\mathsf{E}$. Reciprocally, if $\mathsf{A}_1,
\mathsf{A}_2, \ldots, \mathsf{A}_n$ is a collection of von Neumann
subalgebras of $\mathcal{A}$ freely independent over $\mathsf{E}:
\mathcal{A} \rightarrow \mathcal{B}$ and generating $\mathcal{A}$,
then $\mathcal{A}$ is isomorphic to $*_{\mathcal{B}}
\mathsf{A}_k$.

\begin{remark} \label{Remark-Operator-Valued}
\emph{Let $\mathcal{A} = \mathsf{A}_1 * \mathsf{A}_2 * \cdots *
\mathsf{A}_n$ be a reduced free product von Neumann algebra (i.e.
$\mathcal{A}$ is amalgamated over the complex field) equipped with
its natural \emph{n.f.} state $\phi$. Let $\mathcal{B}$ be another
von Neumann algebra, non necessarily included in $\mathcal{A}$. A
relevant example of the construction outlined above is the
following. Let us consider the conditional expectation
$\mathsf{E}: \mathcal{A} \bar\ten \mathcal{B} \rightarrow
\mathcal{B}$ defined by $\mathsf{E}(a \otimes b) = \phi(a)
1_\mathcal{A} \otimes b$. Then, it is well-known that $
\mathsf{A}_1 \bar\ten \mathcal{B}, \mathsf{A}_2 \bar\ten
\mathcal{B}, \ldots, \mathsf{A}_n \bar\ten \mathcal{B}$ are freely
independent subalgebras of $\mathcal{A} \bar\ten \mathcal{B}$ over
$\mathsf{E}$. In particular, we obtain
 $$\mathcal{M} = \mathcal{A}\bar\ten \mathcal{B}
 = *_{\mathcal{B}} \big( \mathsf{A}_k \bar\ten
 \mathcal{B} \big).$$
Therefore, taking $\mathcal{B}$ to be $\mathcal{B}(\ell_2)$, it
turns out that the complete boundedness of a map $u:
L_p(\mathcal{A}) \to L_p(\mathcal{A})$ is equivalent to the
boundedness (with the same norm) of the map $u \ten
id_{\mathcal{B}}: L_p(\mathcal{M}) \to L_p(\mathcal{M})$. In other
words, since our results are presented for general amalgamated
free products, complete boundedness follows automatically and is
instrumental in some of our arguments. This will be used below
without any further reference. }
\end{remark}

\begin{remark} \label{Remark-Notation-Free}
\emph{Let $\mathcal{A}$ be a von Neumann algebra equipped with a
\emph{n.f.} state $\phi$ and $\mathcal{B}$ a von Neumann
subalgebra of $\mathcal{A}$. According to Takesaki \cite{Ta2}, the
existence and uniqueness of a \emph{n.f.} conditional expectation
$\mathsf{E}: \mathcal{A} \rightarrow \mathcal{B}$ is equivalent to
the invariance of $\mathcal{B}$ under the action of the modular
automorphism group $\sigma_t^{\phi}$ associated to
$(\mathcal{A},\phi)$. Moreover, in that case we have $\phi \circ
\mathsf{E} = \phi$ and following Connes \cite{Co}
 $$\mathsf{E}
 \circ \sigma_t^\phi = \sigma_t^\phi \circ \mathsf{E}.$$
In what follows we shall assume this invariance in all the von
Neumann subalgebras considered. Hence, we may think of a
\emph{natural} conditional expectation $\mathsf{E}: \mathcal{A}
\to \mathcal{B}$. This somehow justifies our relaxed notation for
reduced amalgamated free products, where we do not make explicit
the associated conditional expectations. This should not cause any
confusion since only reduced free product is considered in this
paper.}
\end{remark}

\section{Rosenthal/Voiculescu type inequalities}
\label{Section2}

In this section we present a free analogue of Rosenthal's
inequality (\ref{Rosenthal}). Let $\mathcal{A}$ be the amalgamated
reduced free product $*_{\mathcal{B}} \mathsf{A}_k$ with $1 \le k
\le n$ and $\mathcal{B}$ a common von Neumann subalgebra of the
$\mathsf{A}_k$'s, equipped with a \emph{n.f.} state $\varphi$. As
we have already seen, the state $\varphi$ induces a \emph{n.f.}
state $\phi$ on $\mathcal{A}$ by setting $\phi = \varphi \circ
\mathsf{E}$. Given a non-negative integer $d$ we shall write
$\mathbf{P}_{\mathcal{A}}(d)$ for the closure of elements of the
form
\begin{equation} \label{Eq-d-Polynomial}
a = \sum_{\alpha \in \Lambda}^{\null} \sum_{j_1 \neq j_2 \neq
\cdots \neq j_d} a_{j_1}(\alpha) a_{j_2}(\alpha) \cdots
a_{j_d}(\alpha),
\end{equation}
with $a_{j_k}(\alpha) \in \bubl_{j_k}$ and $\alpha$ running over a
finite set $\Lambda$. In other words, $\mathbf{P}_\mathcal{A}(d)$
is the subspace of $\mathcal{A}$ of homogeneous free polynomials
of degree $d$. When $d$ is $0$, the expression
\eqref{Eq-d-Polynomial} does not make sense.
$\mathbf{P}_\mathcal{A}(0)$ is meant to be $\mathcal{B}$. Then we
define the space $\mathbf{P}_{\mathcal{A}}(p,d)$ as the closure in
$L_p(\Al)$ of
 $$\mathbf{P}_{\mathcal{A}}(d) \, d_{\phi}^{\frac1p}\,,$$
where $d_{\phi}$ denotes the density of the state $\phi$. Note
that, by using approximation with analytic elements, we might have
well located the density $d_\phi$ on the left of
$\mathbf{P}_\Al(d)$ with no consequence in the definition of
$\mathbf{P}_\Al(p,d)$.

\vskip5pt

Similarly, $\mathbf{Q}_\Al(d)$ denotes the subspace of polynomials
of degree less than or equal to $d$ in $\Al$ and
 $$\mathbf{Q}_\Al(p,d) = \bigoplus_{k=0}^d \mathbf{P}_\Al(p,k)
 \quad \mbox{with} \quad \mathbf{P}_{\mathcal{A}}(p,0)
 = L_p(\mathcal{B}).$$

The complementation result below from \cite{RX} is crucial for our
further purposes. Indeed, it was proved there that
$\mathbf{P}_{\mathcal{A}}(d)$ and $\mathbf{Q}_{\mathcal{A}}(d)$
are complemented in $\Al$ with projection constants controlled by
$4d$ and $2d+1$ respectively. Thus, transposition and complex
interpolation yield the following result for $1 \le p \le \infty$.

\begin{theorem} \label{Theorem-Complementation}
The following results hold:
\begin{itemize}
\item[(a)] \hskip1pt $\mathbf{P}_{\mathcal{A}}(p,d)$ is
complemented in $L_p(\mathcal{A})$ with projection constant $\le
4d$.

\item[(b)] $\mathbf{Q}_{\mathcal{A}}(p,d)$ is complemented in
$L_p(\mathcal{A})$ with projection constant $\le 2d + 1$.
\end{itemize}
\end{theorem}

\begin{remark} \label{Remark-Proyecciones}
\emph{In what follows we shall write
$$\quad \Pi_\mathcal{A}(p,d):  L_p(\mathcal{A}) \to
\mathbf{P}_\mathcal{A}(p,d)\quad\mbox{and}\quad
\Gamma_\mathcal{A}(p,d):  L_p(\mathcal{A}) \to
\mathbf{Q}_\mathcal{A}(p,d)$$ for the natural projections
determined by Theorem \ref{Theorem-Complementation}. It is worthy
of mention that both projections above are completely determined
by the natural projections $\Pi_\mathcal{A}(\infty,d)$ and
$\Gamma_\mathcal{A}(\infty,d)$ from \cite{RX}. More precisely,
given $x \in \mathcal{A}$ we have
\begin{equation} \label{Eq-Density-Proyecciones}
\Pi_\mathcal{A}(p,d) \big( x d_\phi^{\frac1p} \big) =
\Pi_\mathcal{A}(\infty,d)(x) d_\phi^{\frac1p} \quad \mbox{and}
\quad \Gamma_\mathcal{A}(p,d) \big( x d_\phi^{\frac1p} \big) =
\Gamma_\mathcal{A}(\infty,d)(x) d_\phi^{\frac1p}.
\end{equation}
In particular, by the density of the subspace $\mathcal{A}\,
d_\phi^{1/p}$ in $L_p(\mathcal{A})$, the relations above
completely determine the projections $\Pi_\mathcal{A}(p,d)$ and
$\Gamma_\mathcal{A}(p,d)$. This will be essential in what follows
for the interpolation of  the spaces
$\mathbf{P}_{\mathcal{A}}(p,d)$ and
$\mathbf{Q}_{\mathcal{A}}(p,d)$ by the complex method. Another
relevant fact implicitly used in the sequel is that both
$\Pi_\mathcal{A}(\infty,d)$ and $\Gamma_\mathcal{A}(\infty,d)$
commute with the modular automorphism group of $\phi$.}
\end{remark}

\subsection{The mappings $\mathcal{L}_k$ and $\mathcal{R}_k$}

Elements in $\bubl_k$ will be called mean-zero letters of
$\mathsf{A}_k$. Given $1 \le k \le n$, we consider the map
$\mathcal{L}_k$ on $\mathbf{P}_\mathcal{A}(p,d)$ which collects
all the reduced words starting with a mean-zero letter in
$\mathsf{A}_k$. Similarly, the map $\mathcal{R}_k$ collects all
the reduced words ending with a mean-zero letter in
$\mathsf{A}_k$. That is, if $a$ is given by the expression
(\ref{Eq-d-Polynomial}) we have
\begin{eqnarray*}
\mathcal{L}_k(a) & = & \sum_{\alpha \in \Lambda} \sum_{j_1 = k
\neq j_2 \neq \cdots \neq j_d} a_{j_1}(\alpha) a_{j_2}(\alpha)
\cdots a_{j_d}(\alpha), \\ \mathcal{R}_k(a) & = & \sum_{\alpha \in
\Lambda} \sum_{j_1 \neq j_2 \neq \cdots \neq k = j_d}
a_{j_1}(\alpha) a_{j_2}(\alpha) \cdots a_{j_d}(\alpha).
\end{eqnarray*}
Of course, both $\mathcal{L}_k$ and $\mathcal{R}_k$ vanish on
$\mathbf{P}_\mathcal{A}(p,0)$. The mappings $\LL_k$ and $\RR_k$
were introduced by Voiculescu \cite{V2}. They are clearly
$\Be$-bimodule maps which commute with the modular automorphism
group and with densities as in \eqref{Eq-Density-Proyecciones}.
Note also that $\mathcal{L}_k$ and $\mathcal{R}_k$ can also be
regarded as orthogonal projections on $L_2(\mathcal{A})$. Thus,
when $p=2$ we need no restriction to the subspaces of homogeneous
polynomials. In this particular case, we shall denote
$\mathcal{L}_k$ and $\mathcal{R}_k$ respectively by $\mathsf{L}_k$
and $\mathsf{R}_k$. Now we prove some fundamental freeness
relations that will be used throughout the whole paper with no
further reference.

\begin{lemma} \label{KeyLemma}
If $1 \le i,j \le n$ and $a_i, a_j \in \mathbf{P}_\mathcal{A}(d)$,
we have
\begin{equation} \label{Eq-Buble1}
\mathsf{L}_i (1 - \mathcal{R}_i)(a_i)^* (1 - \mathcal{R}_j)(a_j)
\mathsf{L}_j = \delta_{ij}\, \mathsf{E} \Big( (1 -
\mathcal{R}_i)(a_i)^* (1 - \mathcal{R}_j)(a_j) \Big) \mathsf{L}_j,
\end{equation}
\begin{equation}  \label{Eq-Buble2}
(1-\mathsf{L}_i) \mathcal{R}_i(a_i)^* \mathcal{R}_j(a_j)
(1-\mathsf{L}_j) = \delta_{ij}\, \mathsf{E} \Big(
\mathcal{R}_i(a_i)^* \mathcal{R}_j(a_j) \Big) (1-\mathsf{L}_j).
\end{equation}
\end{lemma}

\dem By the GNS construction on $(\mathcal{A},\phi)$ we know that
$\mathcal{A}$ acts on $L_2(\mathcal{A})$ by left multiplication.
Thus, we may regard the left hand sides of \eqref{Eq-Buble1} and
\eqref{Eq-Buble2} as mappings on $L_2(\mathcal{A})$. To prove
\eqref{Eq-Buble1}, we first note that
 $$(1 - \mathcal{R}_i)(a_i)^*(1 - \mathcal{R}_j)(a_j)$$
is a linear combination of words of the following form
 $$w_{xy} = x_{i_d}^* x_{i_{d-1}}^* \cdots x_{i_1}^* y_{j_1} \cdots
 y_{j_{d-1}} y_{j_d},$$
where $x_{i_s} \in \bubl_{i_s}$, $y_{j_s} \in \bubl_{j_s}$ and
 $$i_1 \neq i_2 \neq \cdots \neq i_d \neq i,$$
 $$j \neq j_d \neq \cdots \neq j_2 \neq j_1.$$
When $i_1 \neq j_1$, it turns out that $w_{xy}$ is a reduced word
and, since $j_d \neq j$, the map $w_{xy} \mathsf{L}_j$ can only
act as a tensor so that the range of $w_{xy} \mathsf{L}_j$ lies in
the ortho-complement of $\mathsf{L}_i(L_2(\mathcal{A}))$, since
$i_d \neq i$. In other words, in that case we have
 $$\mathsf{L}_iw_{xy} \mathsf{L}_j = 0
 = \mathsf{L}_i \Es(w_{xy}) \mathsf{L}_j.$$
When $i_1 = j_1$ we may write $w_{xy} = w'_{xy} + w''_{xy}$ with
 $$w'_{xy} = x_{i_d}^* x_{i_{d-1}}^* \cdots x_{i_2}^*\,
 \mathsf{E}(x_{i_1}^*y_{j_1})\,y_{j_{2}}
 \cdots y_{j_{d-1}} y_{j_d}.$$
If $d\ge2$, the argument above implies again that $\mathsf{L}_i
w''_{xy} \mathsf{L}_j = 0$ since $w''_{xy}$ is a reduced word not
starting with mean-zero letters in $\mathsf{A}_i$ nor ending with
mean-zero letters in $\mathsf{A}_j$. Then it is clear that we can
iterate the same argument and obtain
 $$\mathsf{L}_i w_{xy} \mathsf{L}_j = \mathsf{L}_i
 \mathsf{E}(x_{i_d}^* \cdots \mathsf{E}(x_{i_1}^* y_{j_1}) \cdots
 y_{j_d}) \mathsf{L}_j = \mathsf{L}_i \mathsf{E}(w_{xy})
 \mathsf{L}_j = \delta_{ij}\, \mathsf{E}(w_{xy}) \mathsf{L}_j.$$
The second identity follows easily by freeness. Summing up we
obtain (\ref{Eq-Buble1}).

\vskip5pt

The proof of \eqref{Eq-Buble2} is quite similar. Indeed, now we
may write $\mathcal{R}_i(a_i)^* \mathcal{R}_j(a_j)$ as a linear
combination of words $w_{xy}$ with the form given above and
satisfying
 $$i_1 \neq i_2 \neq \cdots \neq i_d = i,$$
 $$j = j_d \neq \cdots\neq j_2 \neq j_1.$$
Then the arguments above lead to the following identity
\begin{eqnarray*}
(1 - \mathsf{L}_i) w_{xy} (1 - \mathsf{L}_j) & = & (1 -
\mathsf{L}_i) \mathsf{E}(x_{i_d}^* \cdots \mathsf{E}(x_{i_1}^*
y_{j_1}) \cdots y_{j_d}) (1 - \mathsf{L}_j) \\ & = & \delta_{ij}\,
\mathsf{E}(x_{i_d}^* \cdots \mathsf{E}(x_{i_1}^* y_{j_1}) \cdots
y_{j_d}) (1 - \mathsf{L}_j) \\ & = & \delta_{ij}\,
\mathsf{E}(w_{xy}) (1 - \mathsf{L}_j),
\end{eqnarray*}
where the second identity holds because the only way not to have
 $$\mathsf{E}(x_{i_d}^* \cdots \mathsf{E}(x_{i_1}^* y_{j_1})
 \cdots y_{j_d}) = 0$$
is the case where the indices $i_s$ and $j_s$ fit, i.e. $i_s =
j_s$ for $1 \le s \le d$. Therefore, since $i = i_d$ and $j_d =
j$, the symbol $\delta_{ij}$ appears. Summing up one more time, we
obtain the identity \eqref{Eq-Buble2}. This completes the proof.
\fin

\begin{remark} \label{Remark-Dykema}
\emph{The assumption that $a_i$ and $a_j$ are homogeneous and of
the same degree is essential in Lemma \ref{KeyLemma}. Indeed, the
following counterexample was brought to our attention by Ken
Dykema. Let $\mathbb{F}_2$ denote a free group on two generators
$g_1,g_2$ and let $\lambda: \mathbb{F}_2 \to
\mathcal{B}(\ell_2(\mathbb{F}_2))$ stand for the left regular
representation. Let $\mathsf{A}_k$ be the von Neumann algebra
generated by $\lambda(g_k)$ for $k=1,2$. In this case,
$\mathcal{A} = \mathsf{A}_1 * \mathsf{A}_2$ is the von Neumann
algebra generated by $\lambda$ and the conditional expectation
$\mathsf{E}$ is just $1_\mathcal{A}\, \tau$, where $\tau$ is the
natural trace on $\mathcal{A}$. Then we consider the
(non-homogeneous) polynomial $a = \lambda(g_2) +
\lambda(g_2g_1g_2)$. Clearly, we have $\mathcal{R}_1(a) = 0$ and
 $$a^*a = \big( \lambda(g_2)^* + \lambda(g_2g_1g_2)^* \big)
 \big( \lambda(g_2) + \lambda(g_2g_1g_2) \big) = \mathsf{E}(a^*a) +
 \lambda(g_1g_2) + \lambda(g_1g_2)^*.$$
Taking for instance $h = \lambda(g_1g_2)$, we see that
 $$\mathsf{L}_1 a^* a \mathsf{L}_1 (h) = \mathsf{L}_1
 \mathsf{E}(a^* a) \mathsf{L}_1 (h) + \lambda(g_1g_2g_1g_2) \neq
 \mathsf{L}_1 \mathsf{E}(a^* a) \mathsf{L}_1 (h).$$
Thus identity (\ref{Eq-Buble1}) does not hold for $a$. A similar
counterexample can be constructed for (\ref{Eq-Buble2}). In
particular, since identities (\ref{Eq-Buble1}) and
(\ref{Eq-Buble2}) are essential in most of our results below, this
explains why this paper is written in terms of homogeneous
polynomials.}
\end{remark}

\begin{lemma} \label{Lemma-Complementation-Boundedness}
If $1 \le p \le \infty$ and $a_1,a_2, \ldots, a_n \in
\mathbf{P}_\mathcal{A}(p,d)$,
\begin{eqnarray*}
\Big\| \Big( \sum_{k=1}^n \RR_k(a_k)^* \RR_k(a_k) \Big)^{\frac12}
\Big\|_p & \le & c d^2 \, \Big\| \Big( \sum_{k=1}^n a_k^* a_k
\Big)^{\frac12} \Big\|_p, \\ \Big\| \Big( \sum_{k=1}^n \RR_k(a_k)
\RR_k(a_k)^* \Big)^{\frac12} \Big\|_p & \le & c d^2 \, \Big\|
\Big( \sum_{k=1}^n a_k a_k^* \Big)^{\frac12} \Big\|_p.
\end{eqnarray*}
Moreover, the same inequalities hold with the operator $\LL_k$
instead of $\RR_k$.
\end{lemma}

\dem It is clear that any $a \in \mathbf{P}_\mathcal{A}(p,d)$
satisfies $\mathcal{L}_k(a^*) = \mathcal{R}_k(a)^*$. Consequently,
it suffices to prove the inequalities for the $\RR_k$'s. On the
other hand, in the row/column terminology (i.e. taking $R_p^n$ and
$C_p^n$ to be the first row and column of the Schatten class
$S_p^n$), the two terms on the right hand side are the norms of
$(a_1,a_2,\ldots,a_n)$ in $R_p^n(L_p(\Al))$ and $C_p^n(L_p(\Al))$,
respectively. According to \cite{P2}, both spaces embed
isometrically into
 $$S_p^n(L_p(\Al)) = L_p(\mathrm{M}_n \ten \Al)
 = L_p \Big(*_{\mathrm{M}_n \ten \Be}
 (\mathrm{M}_n \ten \mathsf{A}_k) \Big) = L_p(\mathcal{A}_n).$$
Therefore, by means of Theorem \ref{Theorem-Complementation}
(applied to the amplified algebra $\mathcal{A}_n$), we know that
$\mathbf{P}_{\mathcal{A}_n}(p,d)$ is complemented in
$L_p(\mathcal{A}_n)$ with projection constant $4d$. Using the same
projection restricted to $R_p^n(L_p(\Al))$ and $C_p^n(L_p(\Al))$,
we conclude that the respective subspaces of homogeneous
polynomials $R_p^n(\mathbf{P}_\Al(p,d))$ and
$C_p^n(\mathbf{P}_\Al(p,d))$ form interpolation scales with
equivalent norms up to a constant $4d$. By this observation, it
suffices to show that the assertion holds for $p=1$ and $p=\infty$
with constant in both cases controlled by $c d$. In fact, in the
latter case we shall even prove that the constant does not depend
on $d$. This will be used sometimes in the paper without further
reference. We prove the desired estimates in several steps.

\vskip5pt

\noindent \textbf{Step 1.} Let us prove the first inequality of
$\RR_k$'s for $p=\infty$. The GNS construction on
$(\mathcal{A},\phi)$ implies that $\mathcal{A}$ acts on
$L_2(\mathcal{A})$ by left multiplication. Thus, we may regard
$a_k \Le_k$, $\RR_k(a_k)(1-\Le_k)$ and $(id_\mathcal{A} -
\mathcal{R}_k)(a_k) \mathsf{L}_k$ as mappings on
$L_2(\mathcal{A})$. In particular, since we have $\RR_k(a_k) = a_k
\Le_k + \RR_k(a_k)(1-\Le_k) - (id_\mathcal{A} -
\mathcal{R}_k)(a_k) \mathsf{L}_k$, we obtain by triangle
inequality (with $e_{ij}$ denoting the usual matrix units in
$\mathcal{B}(\ell_2)$)
\begin{eqnarray*}
\Big\| \Big( \sum_{k=1}^n \RR_k(a_k)^* \RR_k(a_k) \Big)^{\frac12}
\Big\|_\infty & \le & \Big\| \sum_{k=1}^n e_{k1} \ten a_k \Le_k
\Big\|_\infty \\ & + & \Big\| \sum_{k=1}^n e_{k1} \ten \RR_k(a_k)
(1-\Le_k) \Big\|_\infty \\ & + & \Big\| \sum_{k=1}^n e_{k1} \ten
(id_\mathcal{A} - \mathcal{R}_k)(a_k) \mathsf{L}_k \Big\|_\infty.
\end{eqnarray*}
If $\mathrm{A}_1, \mathrm{A}_2, \mathrm{A}_3$ denote respectively
the terms on the right, we have
 $$\mathrm{A}_1 = \Big\| \Big( \sum_{k=1}^n
 e_{kk} \ten a_k \Big) \Big( \sum_{k=1}^n e_{k1} \ten \Le_k \Big)
 \Big\|_\infty \le \max_{1 \le k \le n} \|a_k\|_\infty \Big\|
 \sum_{k=1}^n \Le_k \Big\|_\infty^{\frac12}.$$
Thus, since $\sum_k \Le_k = 1 - \Es$, we find
 $$\mathrm{A}_1 \le
 \Big\| \Big( \sum_{k=1}^n a_k^* a_k \Big)^{\frac12}
 \Big\|_\infty.$$
On the other hand, by (\ref{Eq-Buble2}) we have
\begin{eqnarray*}
\mathrm{A}_2 & = & \Big\| \sum_{k=1}^n (1-\Le_k) \RR_k(a_k)^*
\RR_k(a_k) (1-\Le_k) \Big\|_\infty^{\frac12} \\ & = & \Big\|
\sum_{k=1}^n \Es \big( \RR_k(a_k)^* \RR_k(a_k) \big) (1-\Le_k)
\Big\|_\infty^{\frac12}.
\end{eqnarray*}
Now, since $\Le_k$ commutes with $\Be$, the last term is
 $$\summ_k \mathsf{E} (\RR_k(a_k)^* \RR_k(a_k))^{\frac12} (1-\Le_k)
 \mathsf{E} (\RR_k(a_k)^* \RR_k(a_k))^{\frac12} \le \summ_k
 \mathsf{E} (\RR_k(a_k)^* \RR_k(a_k)).$$
Next, we observe that
 \begin{equation} \label{Eq-Expect-Order}
\Es \big(\RR_k(a_k)^* \RR_k(a_k) \big) \le \Es(a_k^*a_k).
\end{equation}
Indeed, since $a_k$ is mean-zero
 $$a_k = \summ_j\mathcal{R}_j(a_k);$$
so, by freeness
 $$\mathsf{E}(a_k^*a_k) = \sum_{i,j} \mathsf{E} \big(
 \mathcal{R}_i(a_k)^* \mathcal{R}_j(a_k) \big) = \sum_j \mathsf{E}
 \big( \mathcal{R}_j(a_k)^* \mathcal{R}_j(a_k) \big) \ge \mathsf{E}
 \big( \mathcal{R}_k(a_k)^* \mathcal{R}_k(a_k) \big).$$
This proves (\ref{Eq-Expect-Order}). Combining the estimates
above, we find
 $$\mathrm{A}_2 \le \Big\| \sum_{k=1}^n \Es(a_k^*a_k)
 \Big\|_\infty^{\frac12} \le \Big\| \Big( \sum_{k=1}^n a_k^*a_k
 \Big)^{\frac12} \Big\|_\infty.$$
The estimate for $\mathrm{A}_3$ is similar to the one for
$\mathrm{A}_2$ and we leave it to the reader.

\vskip5pt

\noindent \textbf{Step 2.} Now we prove the second inequality for
$\RR_k$'s. As above, we have
\begin{eqnarray*}
\Big\| \Big( \sum_{k=1}^n \RR_k(a_k) \RR_k(a_k)^* \Big)^{\frac12}
\Big\|_\infty & \le & \Big\| \sum_{k=1}^n e_{1k} \ten a_k \Le_k
\Big\|_\infty \\ & + & \Big\| \sum_{k=1}^n e_{1k} \ten \RR_k(a_k)
(1-\Le_k) \Big\|_\infty \\ & + & \Big\| \sum_{k=1}^n e_{1k} \ten
(id_\mathcal{A} - \mathcal{R}_k)(a_k) \mathsf{L}_k \Big\|_\infty.
\end{eqnarray*}
We write $\mathrm{B}_1, \mathrm{B}_2, \mathrm{B}_3$ for the terms
on the right. The estimate of $\mathrm{B}_1$ is trivial
 $$\mathrm{B}_1 = \Big\| \Big( \sum_{k=1}^n e_{1k} \ten a_k \Big)
 \Big( \sum_{k=1}^n e_{kk} \ten \Le_k \Big) \Big\|_\infty \le
 \Big\| \Big( \sum_{k=1}^n a_k a_k^* \Big)^{\frac12}\Big\|_\infty.$$
On the other hand, by (\ref{Eq-Buble2}) and
\eqref{Eq-Expect-Order} we may write
\begin{eqnarray*}
\mathrm{B}_2 & = & \Big\| \sum_{i,j=1}^n e_{ij} \ten (1-\Le_i)
\RR_i(a_i)^* \RR_j(a_j) (1-\Le_j) \Big\|_\infty^{\frac12} \\ & = &
\Big\| \sum_{k=1}^n e_{kk} \ten \Es \big( \RR_k(a_k)^* \RR_k(a_k)
\big) (1-\Le_k) \Big\|_\infty^{\frac12} \\ & \le & \max_{1 \le k
\le n} \big\| \Es \big( \RR_k(a_k)^* \RR_k(a_k) \big)
\big\|_\infty^{\frac12} \le \max_{1 \le k \le n} \big\| \Es \big(
a_k^* a_k \big) \big\|_\infty^{\frac12} \le \max_{1 \le k \le n}
\|a_k\|_\infty.
\end{eqnarray*}
Finally, to estimate $\mathrm{B}_3$ we use \eqref{Eq-Buble1} and
the proof of \eqref{Eq-Expect-Order}
\begin{eqnarray*} \mathrm{B}_3 & = & \Big\|
\sum_{i,j=1}^n e_{ij} \ten \Le_i (a_i-\RR_i(a_i))^*
(a_j-\RR_j(a_j)) \Le_j \Big\|_\infty^{\frac12}
\\ & = & \Big\| \sum_{k=1}^n e_{kk} \ten \Es \big(
(a_k-\RR_k(a_k))^* (a_k-\RR_k(a_k)) \big) \Le_k
\Big\|_\infty^{\frac12}
\\ & \le & \max_{1 \le k \le n} \big\| \Es \big(
(a_k-\RR_k(a_k))^* (a_k-\RR_k(a_k)) \big) \big\|_\infty^{\frac12} \\
& \le & \max_{1 \le k \le n} \big\| \Es \big( a_k^* a_k \big)
\big\|_\infty^{\frac12} \le \max_{1 \le k \le n} \|a_k\|_\infty.
\end{eqnarray*}

\vskip5pt

\noindent \textbf{Step 3.} Now we use a duality argument to prove
the same estimates in $L_1(\Al)$. Recall $d_\phi$ denotes the
density associated to the state $\phi$ of $\Al$. Let $a \in
\mathbf{P}_\mathcal{A}(d)$ and $x \in \mathcal{A}$ be a finite sum
of reduced words. Then we have
\begin{eqnarray}
\label{Eq-Dualizacion} \big\langle x, \mathcal{R}_k(d_\phi a)
\big\rangle & = & \mbox{tr}_\mathcal{A} \big( x^*
\mathcal{R}_k(d_\phi a) \big) = \mbox{tr}_\mathcal{A} \big( d_\phi
\mathcal{R}_k(a) x^* \big) \\ \nonumber & = & \varphi \big(
\mathsf{E} (\mathcal{R}_k(a) x^*) \big) \ \, =  \varphi \big(
\mathsf{E} (a \mathcal{L}_k(x^*)) \big) \\ \nonumber & = &
\mbox{tr}_\mathcal{A} \big( d_\phi a \mathcal{R}_k(x)^* \big) =
\big\langle \mathcal{R}_k(x), d_\phi a \big\rangle.
\end{eqnarray}
Moreover, arguing as we did in the proof of Lemma \ref{KeyLemma},
it can be checked that any word in $x$ of length different from
$d$ does not contribute to the quantity considered in
\eqref{Eq-Dualizacion}. Let us define
\begin{eqnarray*}
R_\Al & = & \Big\{ (x_1,x_2,\ldots,x_n) \, \big| \ x_k \in \Al, \
\Big\| \summ_k x_k x_k^* \Big\|_\infty \le 1 \Big\}, \\
C_\Al & = & \Big\{ (x_1,x_2,\ldots,x_n) \, \big| \ x_k \in \Al, \
\Big\| \summ_k x_k^* x_k \Big\|_\infty \le 1 \Big\},
\end{eqnarray*}
and let us also consider the sets
\begin{eqnarray*}
R_\Al(d) & = & \Big\{ (x_1,x_2,\ldots,x_n) \, \big| \ x_k \in
\mathbf{P}_\Al(d), \ \Big\| \summ_k x_k x_k^* \Big\|_\infty \le 1
\Big\}, \\ C_\Al(d) & = & \Big\{ (x_1,x_2,\ldots,x_n) \, \big| \
x_k \in \mathbf{P}_\Al(d), \ \Big\| \summ_k x_k^* x_k
\Big\|_\infty \le 1 \Big\}.
\end{eqnarray*}
Arguing as at the beginning of this proof, we know from Remarks
\ref{Remark-Notation-Free} and \ref{Remark-Proyecciones} that the
projection $\Pi_\mathcal{A}(\infty,d)$ is in fact completely
bounded. This means that we may work on the amplified algebra
$\mathcal{A}_n = \mathrm{M}_n \ten \mathcal{A}$ and obtain
projections
 $$\begin{array}{rrcl} \Pi_d =
 \Pi_{\mathcal{A}_n}(\infty,d): &
 R_n(\mathcal{A}) & \to & R_n(\mathbf{P}_\mathcal{A}(d)), \\
 \Pi_d = \Pi_{\mathcal{A}_n}(\infty,d): & C_n(\mathcal{A}) & \to &
 C_n(\mathbf{P}_\mathcal{A}(d)), \end{array}$$
bounded by $4d$. In particular, we obtain the inclusions $\Pi_d
(R_\mathcal{A}) \subset 4d  \, R_\mathcal{A}(d)$ and $\Pi_d
(C_\mathcal{A}) \subset 4d \, C_\mathcal{A}(d)$. Therefore, since
the words of length different from $d$ do not contribute in
\eqref{Eq-Dualizacion}
\begin{eqnarray*}
\lefteqn{\Big\| \summ_k e_{1k} \ten \mathcal{R}_k(a_k) \Big\|_1}
\\ & = & \sup_{x \in R_\Al} \summ_k \big\langle x_k,
\mathcal{R}_k(a_k) \big\rangle \\ & = & \sup_{x \in
\Pi_d(R_\Al)} \, \summ_k \big\langle \RR_k(x_k), a_k \big\rangle \\
& \le & \sup_{x \in
4d R_\Al(d)} \summ_k \big\langle \RR_k(x_k), a_k \big\rangle \\
& \le & \sup_{x \in 4d R_\Al(d)} \Big\| \Big( \summ_k \RR_k(x_k)
\RR_k(x_k)^* \Big)^{\frac12} \Big\|_\infty \Big\| \Big( \summ_k
a_k a_k^* \Big)^{\frac12} \Big\|_1.
\end{eqnarray*}
By Step 2,
 $$\Big\| \Big( \summ_k \RR_k(a_k) \RR_k(a_k)^*
 \Big)^{\frac12} \Big\|_{1} \le c d \Big\| \Big( \summ_k a_k a_k^*
 \Big)^{\frac12} \Big\|_{1}$$
Similarly, using Step 1 and the space $C_\Al(d)$, we obtain the
remaining estimate. \fin

\begin{remark} \label{Remark-Diffp>2}
\emph{A detailed reading of the proof of Lemma
\ref{Lemma-Complementation-Boundedness} shows that the constant is
controlled by $cd^2$ for $1 \le p \le 2$ and by $c d$ for $2 \le p
\le \infty$. Moreover, the same arguments are valid to show that
Lemma \ref{Lemma-Complementation-Boundedness} also holds replacing
$\LL_k$ or $\RR_k$ by $\QQ_k = \LL_k \RR_k = \RR_k \LL_k$ (to be
used below). These generalizations of Lemma
\ref{Lemma-Complementation-Boundedness} will be used several times
in the sequel.}
\end{remark}

\begin{remark} \emph{In Remark \ref{Remark-Dykema} we
have partially justified why this paper is written in terms of
homogeneous polynomials. On the other hand, Lemma
\ref{Lemma-Complementation-Boundedness} for $n = 1$ shows that
$\LL_k$ and $\RR_k$ are bounded operators when acting on
$\mathbf{P}_\Al(p,d)$ for any $1 \le p \le \infty$ and $d \ge 1$.
Another relevant fact which justifies the use of homogeneous
polynomials is that $\LL_k$ and $\RR_k$ are not bounded on
$L_\infty(\Al)$. The following simple counterexample was brought
to our attention by Ana Maria Popa. Consider again the free group
$\mathbb{F}_2$ with two generators  $g_1,g_2$ and  keep the
terminology employed in Remark \ref{Remark-Dykema}. Let
$\mathsf{H}$ be the subgroup of $\mathbb{F}_2$ generated by
$w=g_1g_2$. Of course, it is clear that $\mathsf{H}$ is isomorphic
to $\mathbb{Z}$ and that $\lambda(\mathsf{H})'' \simeq
L_\infty(\mathbb{T})$. Moreover, we obviously have
$\LL_1(\lambda(w^k)) = \delta_{k>0} \, \lambda(w^k)$. In
particular, if $\mathcal{A} = \lambda(\mathbb{F}_2)''$ denotes the
reduced group von Neumann algebra, it turns out that the
restriction of $\LL_1: \Al \to \Al$ to $\lambda(\mathsf{H})$
behaves as $\frac12(id_{L_\infty(\mathbb{T})} + \mathbf{H})$,
where $\mathbf{H}$ denotes the Hilbert transform on the circle.
The claim follows since the Hilbert transform is known to be
unbounded on $L_\infty(\mathbb{T})$. Moreover, the same example
also shows that the map $\mathcal{Q}_k = \mathcal{R}_k
\mathcal{L}_k$ (to be used below) is unbounded on
$L_\infty(\mathcal{A})$. Indeed, $\mathcal{Q}_1$ is not bounded on
the subspace $\lambda(\mathsf{H} g_1)''$ since
 $$\mathcal{Q}_1 \big( \lambda(w^k g_1) \big) =
 \delta_{k>0} \lambda(w^k g_1).$$}
\end{remark}

\vskip5pt

\begin{proposition} \label{Proposition-Sum-L&R}
If $a_1, a_2, \ldots, a_n \in \mathbf{P}_\mathcal{A}(d)$, we have
\begin{eqnarray*}
\Big\| \sum_{k=1}^n \mathcal{L}_k(a_k) \Big\|_{\infty} \!\!\!\! &
\sim_c & \!\!\! \max \left\{ \Big\| \sum_{k=1}^n
\mathcal{L}_k(a_k)^* \mathcal{L}_k(a_k) \Big\|_{\infty}^{\frac12},
\ \Big\| \sum_{k=1}^n \mathsf{E}(\mathcal{L}_k(a_k)
\mathcal{L}_k(a_k)^*) \Big\|_{\infty}^{\frac12} \hskip2.5pt
\right\}, \\ \Big\| \sum_{k=1}^n \mathcal{R}_k(a_k)
\Big\|_{\infty} \!\!\!\! & \sim_c & \!\!\! \max \left\{ \Big\|
\sum_{k=1}^n \mathcal{R}_k(a_k) \mathcal{R}_k(a_k)^*
\Big\|_{\infty}^{\frac12}, \Big\| \sum_{k=1}^n
\mathsf{E}(\mathcal{R}_k(a_k)^* \mathcal{R}_k(a_k))
\Big\|_{\infty}^{\frac12} \right\}.
\end{eqnarray*}
\end{proposition}

\dem Once more, we only prove the assertion for $\mathcal{R}_k$.
We have
\begin{eqnarray*}
\lefteqn{\Big\| \sum_{k=1}^n \mathcal{R}_k(a_k) \Big\|_\infty}
\\ & \le & \Big\| \sum_{k=1}^n \RR_k(a_k) \Le_k \Big\|_\infty +
\Big\| \sum_{k=1}^n \RR_k(a_k) (1-\Le_k) \Big\|_\infty \\ & = &
\Big\| \sum_{k=1}^n \RR_k(a_k) \Le_k \RR_k(a_k)^*
\Big\|_\infty^{\frac12} + \Big\| \sum_{i,j=1}^n (1-\Le_i)
\RR_i(a_i)^* \RR_j(a_j) (1-\Le_j) \Big\|_\infty^{\frac12} \\ & = &
\Big\| \sum_{k=1}^n \RR_k(a_k) \Le_k \RR_k(a_k)^*
\Big\|_\infty^{\frac12} + \Big\| \sum_{k=1}^n \mathsf{E} \big(
\RR_k(a_k)^* \RR_k(a_k) \big) (1-\Le_k) \Big\|_\infty^{\frac12}.
\end{eqnarray*}
The first term is clearly bounded by $\sum_k \RR_k(a_k)
\RR_k(a_k)^*$. For the second term we argue as in the proof of
Lemma \ref{Lemma-Complementation-Boundedness}. That is, using that
$\Le_k$ commutes with $\Be$, we can write $\sum_k \mathsf{E}
(\RR_k(a_k)^* \RR_k(a_k)) (1-\Le_k)$ as
 $$\summ_k \mathsf{E} \big(\RR_k(a_k)^* \RR_k(a_k) \big)^{\frac12}
 (1-\Le_k) \mathsf{E} \big(\RR_k(a_k)^* \RR_k(a_k) \big)^{\frac12}.$$
Thus, we obtain the upper estimate
 $$\Big\| \sum_{k=1}^n\mathcal{R}_k(a_k) \Big\|_\infty
 \le \Big\| \sum_{k=1}^n\RR_k(a_k) \RR_k(a_k)^*
 \Big\|_\infty^{\frac12} + \Big\|\sum_{k=1}^n \mathsf{E}
 \big( \RR_k(a_k)^* \RR_k(a_k) \big)\Big\|_\infty^{\frac12}.$$

\vskip5pt

For the lower estimate, using freeness, we clearly have
\begin{equation} \label{Eq-Est2}
\Big\| \sum_{k=1}^n \Es \big( \RR_k(a_k)^* \RR_k(a_k) \big)
\Big\|_\infty \! =  \Big\| \sum_{i,j=1}^n \Es \big( \RR_i(a_i)^*
\RR_j(a_j) \big) \Big\|_\infty \le \Big\| \sum_{k=1}^n \RR_k(a_k)
\Big\|_\infty^2.
\end{equation}
Thus, it remains to show that
 $$\Big\| \sum_{k=1}^n \RR_k(a_k)
 \RR_k(a_k)^* \Big\|_\infty^{\frac12} \le c \Big\| \sum_{k=1}^n
 \RR_k(a_k) \Big\|_\infty.$$
To that aim we observe from (\ref{Eq-Est2}) and the calculation
above that
 $$\Big\| \sum_{k=1}^n \RR_k(a_k) \Le_k \Big\|_\infty
 \le \Big\| \sum_{k=1}^n \RR_k(a_k) \Big\|_\infty + \Big\|
 \sum_{k=1}^n \RR_k(a_k) (1-\Le_k) \Big\|_\infty \le 2 \Big\|
 \sum_{k=1}^n \RR_k(a_k) \Big\|_\infty.$$
Hence, since the term
 $$\mathcal{S}(a,\varepsilon) = \Big\| \sum_{k=1}^n \varepsilon_k
 \RR_k(a_k) \Le_k \Big\|_\infty + \Big\| \sum_{k=1}^n \Es \big(
 (\varepsilon_k \RR_k(a_k))^* (\varepsilon_k \RR_k(a_k)) \big)
 \Big\|_\infty^{\frac12}$$
is independent of any choice of signs $\varepsilon=(\varepsilon_1,
\ldots, \varepsilon_n) \in \Omega = \{ \pm 1\}^n$, we find
\begin{equation} \label{Eq-Est3}
\Big\| \sum_{k=1}^n \varepsilon_k \RR_k(a_k) \Big\|_\infty \le
\mathcal{S}(a,\varepsilon) \le 3 \Big\| \sum_{k=1}^n \RR_k(a_k)
\Big\|_\infty.
\end{equation}
Therefore, we obtain
\begin{eqnarray*}
\lefteqn{\Big\| \sum_{k=1}^n \RR_k(a_k) \RR_k(a_k)^*
\Big\|_\infty} \\ & = & \Big\| \int_{\Omega} \sum_{i,j=1}^n
\varepsilon_i \RR_i(a_i) \varepsilon_j \RR_j(a_j)^* d \varepsilon
\Big\|_\infty \\ & \le & \int_{\Omega} \Big\| \sum_{i,j=1}^n
\varepsilon_i \RR_i(a_i) \varepsilon_j \RR_j(a_j)^* \Big\|_\infty
d \varepsilon \\ & \le & \int_{\Omega} \Big\| \sum_{i=1}^n
\varepsilon_i \RR_i(a_i) \Big\|_\infty \Big\| \sum_{j=1}^n
\varepsilon_j \RR_j(a_j)^* \Big\|_\infty d \varepsilon \le 9
\Big\| \sum_{k=1}^n \RR_k(a_k) \Big\|_\infty^2.
\end{eqnarray*}
This is the remaining inequality to complete the proof of the
lower estimate. \fin

\begin{corollary} \label{Corollary-Sum-L&R}
If $2 \le p \le \infty$ and $a_1,a_2, \ldots, a_n \in
\mathbf{P}_\Al(p,d)$, we have
\begin{eqnarray*}
\Big\| \sum_{k=1}^n \LL_k(a_k) \Big\|_p \!\!\! & \le & \!\!\! c
d^2 \max \left\{ \Big\| \sum_{k=1}^n \LL_k(a_k) \LL_k(a_k)^*
\Big\|_{p/2}^{\frac12}, \, \Big\| \sum_{k=1}^n
\LL_k(a_k)^* \LL_k(a_k) \Big\|_{p/2}^{\frac12} \, \right\}, \\
\Big\| \sum_{k=1}^n \RR_k(a_k) \Big\|_p \!\!\! & \le & \!\!\! c
d^2 \max \left\{ \Big\| \sum_{k=1}^n \RR_k(a_k) \RR_k(a_k)^*
\Big\|_{p/2}^{\frac12}, \! \Big\| \sum_{k=1}^n \RR_k(a_k)^*
\RR_k(a_k) \Big\|_{p/2}^{\frac12} \! \right\}.
\end{eqnarray*}
\end{corollary}

\dem We only prove the second inequality. According to Proposition
\ref{Proposition-Sum-L&R}, the case $p=\infty$ follows with
constant $3$ while the case $p=2$ holds with constant $1$ by
orthogonality. Therefore, it suffices to show that we can
interpolate. To that aim we observe that the term of the right
hand side can be rewritten as
 $$\max \left\{ \Big\| \sum_{k=1}^n e_{1k} \ten \RR_k(a_k)
 \Big\|_{S_p^n(L_p(\Al))}, \Big\| \sum_{k=1}^n e_{k1} \ten
 \RR_k(a_k) \Big\|_{S_p^n(L_p(\Al))} \right\}.$$
In other words, this is the norm of $\big( \RR_k(a_k) \big)$ in
 $$RC_p^n(L_p(\Al)) = R_p^n(L_p(\Al)) \cap C_p^n(L_p(\Al)).$$
On the other hand, by Theorem \ref{Theorem-Complementation} we
know that $R_p^n(\mathbf{P}_\Al(p,d))$ and
$C_p^n(\mathbf{P}_\Al(p,d))$ are complemented respectively in
$R_p^n(L_p(\Al))$ and $C_p^n(L_p(\Al))$ with projection cons\-tant
less than or equal to $4d$. Thus, taking the same projection on
$RC_p^n(L_p(\Al))$ and using that $RC_p^n(L_p(\Al))$ is an
interpolation scale (see \cite{LuP,P2}), we conclude that
$RC_p^n(\mathbf{P}_\Al(p,d))$ is an interpolation scale with
equivalent norms up to a constant controlled by $cd$. Then we need
to consider the subspace of $RC_p^n(\mathbf{P}_\Al(p,d))$ made up
of elements for which the $k$-th component is in
$\RR_k(L_p(\Al))$. The associated projection is
 $$\Pi_\RR =\summ_k \delta_k \ten \RR_k,$$
where $(\delta_k)$ denotes the common basis of $R$ and $C$ when
$(R,\, C)$ is viewed as a compatible couple. According to Lemma
\ref{Lemma-Complementation-Boundedness} and Remark
\ref{Remark-Diffp>2}, the projection $\Pi_\RR$ is bounded and of
norm $\le c d$. Therefore, the family of spaces
$\Pi_\RR(RC_p^n(\mathbf{P}_\Al(p,d)))$, $2\le p\le\infty$, forms
an interpolation scale with equivalent norms up to a constant
controlled by $c d^2$. This completes the proof. \fin

\subsection{Proof of Theorem A and applications}

We now study generalizations of Voiculescu's inequality \cite{V2},
originally formulated for $1$-homogeneous polynomials in a free
product von Neumann algebra. Our main result is Theorem A (stated
in the Introduction), which extends Voiculescu's inequality in
three aspects: we allow amalgamation, homogeneous free polynomials
of arbitrary degree and our inequalities hold in
$L_p(\mathcal{A})$ for $2 \le p \le \infty$. In particular,
Theorem A can be regarded as a generalization of Rosenthal's
inequality (\ref{Rosenthal}) in the free setting.

\vskip5pt

The notation
  $$\QQ_k = \RR_k \LL_k = \LL_k \RR_k$$
for the projection onto words starting and ending in $\bubl_k$ is
crucial for our analysis.

\begin{lemma} \label{Lemma-Est1-Q}
If $a \in \mathbf{P}_\Al(d)$, we have
 $$\max_{1\le k \le n} \big\|\QQ_k(a) \big\|_\infty
 + \Big\| \sum_{k=1}^n \Es \big( \QQ_k(a)^*\QQ_k(a) \big)
 \Big\|_\infty^{\frac12} + \Big\| \sum_{k=1}^n \Es\big( \QQ_k(a)
 \QQ_k(a)^* \big) \Big\|_\infty^{\frac12} \le c\|a\|_\infty.$$
Moreover, if $a_1,a_2,\ldots,a_n \in \mathbf{P}_\Al(d)$, we have
\begin{eqnarray*}
\Big\| \sum_{k=1}^n \QQ_k(a_k) \Big\|_\infty & \sim_c & \max_{1\le
k \le n} \big\| \QQ_k(a_k) \big\|_\infty \\ & + & \Big\| \Big(
\sum_{k=1}^n \Es \big( \QQ_k(a_k)^* \QQ_k(a_k) \big)
\Big)^{\frac12} \Big\|_\infty \\ & + & \Big\| \Big( \sum_{k=1}^n
\Es \big( \QQ_k(a_k) \QQ_k(a_k)^* \big) \Big)^{\frac12}
\Big\|_\infty.
\end{eqnarray*}

\end{lemma}

\dem According to the proof of Lemma
\ref{Lemma-Complementation-Boundedness}, we know that $\LL_k$ and
$\RR_k$ are bounded maps on $\mathbf{P}_\Al(d)$ with constant $3$.
In particular, we find $\|\QQ_k(a)\|_\infty \le 9 \|a\|_\infty$.
On the other hand, using the identities
 $$\Es(a^*a) = \summ_k \Es \big( \LL_k(a)^* \LL_k(a) \big) =
 \summ_k \Es \big( \RR_k(a)^* \RR_k(a) \big),$$
for homogeneous polynomials (\emph{c.f.} the proof of
(\ref{Eq-Expect-Order})) we easily obtain the estimate
\begin{eqnarray*}
\Big\| \sum_{k=1}^n \Es \big( \QQ_k(a)^* \QQ_k(a) \big)
\Big\|_\infty & = & \Big\| \sum_{k=1}^n \Es \big(
\RR_k(\LL_k(a))^* \RR_k(\LL_k(a)) \big) \Big\|_\infty \\ & \le &
\Big\| \sum_{k=1}^n \Es \big( \LL_k(a)^* \LL_k(a) \big)
\Big\|_\infty = \big\| \Es(a^*a) \big\|_\infty \le \|a\|_\infty^2.
\end{eqnarray*}
Using this estimate for $a^*$, we deduce the first assertion.

\vskip5pt

To prove the second one we note that $\QQ_k(a_k) = \QQ_k(a)$ for
$a = \sum_k \QQ_k(a_k)$. In particular, the lower estimate follows
from the first assertion. For the upper estimate we use
 $$\sum_{k=1}^n \QQ_k(a_k) = \sum_{k=1}^n \Le_k \QQ_k(a_k) \Le_k +
 \sum_{k=1}^n \QQ_k(a_k) (1-\Le_k) + \sum_{k=1}^n
 (1-\Le_k)\QQ_k(a_k) \Le_k.$$
The first term gives the maximum. The remaining terms are
estimated by (\ref{Eq-Buble2}). \fin

\begin{lemma} \label{Lemma-Unconditionallity}
Let $a_k \in \mathbf{P}_\Al(p,d)$ and signs $\varepsilon_k =
\pm1.$
\begin{itemize}
\item[i)] If $1 \le p< 2$, we have
 $$\Big\| \sum_{k=1}^n\varepsilon_k \QQ_k(a_k)
 \Big\|_p \le c d^2 \Big\| \sum_{k=1}^n\QQ_k(a_k) \Big\|_p.$$

\item[ii)] If $2 \le p \le \infty$, we have
 $$\Big\| \sum_{k=1}^n\varepsilon_k \QQ_k(a_k)
 \Big\|_p \le c d \, \Big\| \sum_{k=1}^n\QQ_k(a_k) \Big\|_p.$$
\end{itemize}
\end{lemma}

\dem If $a \in \mathbf{P}_\mathcal{A}(p,d)$, we claim that
\begin{equation} \label{Eq-Est33}
\begin{array}{rclcl} \displaystyle \Big\| \sum_{k=1}^n
\varepsilon_k \RR_k(a) \Big\|_p & \le & \displaystyle c d^2 \Big\|
\sum_{k=1}^n \RR_k(a) \Big\|_p & \mbox{for} & 1 \le p < 2, \\
[10pt] \displaystyle \Big\| \sum_{k=1}^n \varepsilon_k \RR_k(a)
\Big\|_p & \le & \displaystyle c d \hskip4.3pt \Big\| \sum_{k=1}^n
\RR_k(a) \Big\|_p & \mbox{for} & 2 \le p \le \infty.
\end{array}
\end{equation}
The second inequality clearly holds with constant $1$ for $p=2$.
On the other hand, according to (\ref{Eq-Est3}), it also holds for
$p=\infty$ with constant $3$. Therefore, since any $a \in
\mathbf{P}_\Al(p,d)$ satisfies $a = \sum_{k} \RR_k(a)$, our claim
follows for $2 \le p \le \infty$ by complex interpolation from
Theorem \ref{Theorem-Complementation}.

\vskip5pt

Then a duality argument yields the first inequality in the claim.
Indeed, by Theo\-rem \ref{Theorem-Complementation} one more time,
we have $\mathbf{P}_\mathcal{A}(p,d)^* \simeq
\mathbf{P}_\mathcal{A}(p',d)$ with equivalence constant controlled
by $4d$. Therefore, given $1 \le p \le 2$, an element $a \in
\mathbf{P}_\Al(p,d)$ and signs $\varepsilon_1, \varepsilon_2,
\ldots, \varepsilon_n$, we choose $x \in \mathbf{P}_\Al(p',d)$ of
norm one such that
\begin{eqnarray*}
\Big\| \sum_{k=1}^n \varepsilon_k \RR_k(a) \Big\|_p & \le & 4d \,
\mbox{tr}_\Al \Big( x^* \sum_{k=1}^n \varepsilon_k \RR_k(a) \Big)
\\ & = & 4d \, \mbox{tr}_\Al \Big(
\sum_{k=1}^n \varepsilon_k \LL_k(x^*) \, a \Big) \\ & \le & 4d \,
\|a\|_p \Big\| \sum_{k=1}^n \varepsilon_k \RR_k(x) \Big\|_{p'} \le
c d^2 \|a\|_p.
\end{eqnarray*}

Taking $a = \sum_k \RR_k(a_k)$, we see that  \eqref{Eq-Est33}
implies
\begin{equation} \label{Eq-Est333}
\begin{array}{rclcl} \displaystyle \Big\| \sum_{k=1}^n
\varepsilon_k \RR_k(a_k) \Big\|_p & \le & \displaystyle c d^2
\Big\| \sum_{k=1}^n \RR_k(a_k) \Big\|_p & \mbox{for} & 1 \le p < 2, \\
[10pt] \displaystyle \Big\| \sum_{k=1}^n \varepsilon_k \RR_k(a_k)
\Big\|_p & \le & \displaystyle c d \hskip4.3pt \Big\| \sum_{k=1}^n
\RR_k(a_k) \Big\|_p & \mbox{for} & 2 \le p \le \infty.
\end{array}
\end{equation}
Therefore, the lemma immediately follows from \eqref{Eq-Est333}
since $\QQ_k = \RR_k \QQ_k$. \fin

\begin{lemma} \label{Lemma-Type}
If $1 \le p \le 2$ and $a_1, a_2, \ldots, a_n \in
\mathbf{Q}_\Al(p,d)$, we have
 $$\Big\| \sum_{k=1}^n \QQ_k(a_k)\Big\|_p
 \le c d^4 \Big( \sum_{k=1}^n \|a_k\|_p^p\Big)^{\frac1p}.$$
\end{lemma}

\dem Using the boundedness of the projection
$\Gamma_\mathcal{A}(p,d)$ from Remark \ref{Remark-Proyecciones}
and complex interpolation, it suffices to see that the
inequalities associated to the extremal indices hold with constant
controlled by $c d^3$. In the case $p=2$, this follows by
orthogonality with constant $1$. When $p=1$, we decompose the
$a_k$'s into their homogeneous parts and use the boundedness of
 $$\QQ_k \circ \Pi_\mathcal{A}(1,s): L_1(\mathcal{A}) \to
 \mathbf{P}_\mathcal{A}(1,s).$$
Indeed, by Step 3 in the proof of Lemma
\ref{Lemma-Complementation-Boundedness} and Remark
\ref{Remark-Diffp>2} we have
$$\big\| \QQ_k \circ \Pi_\mathcal{A}(1,s) \big\|_1 \le c(1+s)
\big\| \Pi_\mathcal{A}(1,s) \big\|_1.$$ Therefore, we find
\begin{eqnarray*}
\Big\| \sum_{k=1}^n \QQ_k(a_k) \Big\|_1 & \le & \sum_{k=1}^n
\big\| \QQ_k(a_k) \big\|_1  \le  \sum_{k=1}^n \sum_{s=0}^d \big\|
\QQ_k \big( \Pi_\mathcal{A}(1,s)(a_k) \big) \big\|_1 \\ & \le &
\sum_{k=1}^n \sum_{s=0}^d c(1+s) \big\| \Pi_\mathcal{A}(1,s)(a_k)
\big\|_1  \le   c \sum_{k=1}^n \sum_{s=0}^d (1+s)^2 \|a_k\|_1
\\ & = & c \, \Big( \sum_{s=0}^d (1+s)^2 \Big) \Big( \sum_{k=1}^n
\|a_k\|_1 \Big) \le  c d^3 \sum_{k=1}^n \|a_k\|_1.
\end{eqnarray*}
This proves the remaining estimate. The proof is complete. \fin

\demA Lemma \ref{Lemma-Est1-Q} implies the assertion for
$p=\infty$. Thus, we may assume in what follows that $2 \le p <
\infty$. Let us prove the lower estimate. First we observe that
$L_p(\Al)$ has Rademacher cotype $p$ for $2 \le p < \infty$. This,
combined with Lemma \ref{Lemma-Unconditionallity} yields
\begin{equation} \label{Eq-QQpp}
\Big( \sum_{k=1}^n \big\| \QQ_k(a_k) \big\|_p^p \Big)^{\frac1p}
\le \int_\Omega \Big\| \sum_{k=1}^n \varepsilon_k \QQ_k(a_k)
\Big\|_p d\varepsilon \le c d \, \Big\| \sum_{k=1}^n \QQ_k(a_k)
\Big\|_p.
\end{equation}
For the second term we use
 $$\sum_{k=1}^n \mathsf{E} \big(\QQ_k(a_k)^* \QQ_k(a_k) \big)
 = \sum_{i,j=1}^n \mathsf{E} \big(\QQ_i(a_i)^* \QQ_j(a_j) \big).$$
Hence, by the contractivity of $\Es$
 $$\Big\| \sum_{k=1}^n\mathsf{E} \big( \QQ_k(a_k)^* \QQ_k(a_k)\big)
  \Big\|_{p/2} \le\Big\| \sum_{k=1}^n \QQ_k(a_k) \Big\|_p^2.$$
The third term is estimated in the same way. Therefore, the lower
estimate holds with constant $c d$. Now we prove the upper
estimate. To that aim we proceed in two steps. First we prove the
case $2 \le p \le 4$ and after that we shall apply an induction
argument.

\vskip5pt

\noindent \textbf{Step 1.} Since $\RR_k(\QQ_k(a_k)) = \QQ_k(a_k)$,
we may apply Corollary \ref{Corollary-Sum-L&R} and obtain
\begin{equation} \label{Equation-Start}
\Big\| \sum_{k=1}^n \QQ_k(a_k) \Big\|_p \le c d^2 \left( \Big\|
\sum_{k=1}^n \QQ_k(a_k) \QQ_k(a_k)^*
\Big\|_{\frac{p}{2}}^{\frac12} + \Big\| \sum_{k=1}^n \QQ_k(a_k)^*
\QQ_k(a_k) \Big\|_{\frac{p}{2}}^{\frac12} \right).
\end{equation}
Then we observe that
\begin{eqnarray}
\label{left} \QQ_k(a_k) \QQ_k(a_k)^* & = & \Es \big( \QQ_k(a_k)
\QQ_k(a_k)^* \big) + \QQ_k \big( \QQ_k(a_k) \QQ_k(a_k)^* \big), \\
\label{right} \QQ_k(a_k)^* \QQ_k(a_k) & = & \Es \big( \QQ_k(a_k)^*
\QQ_k(a_k) \big) + \QQ_k \big( \QQ_k(a_k)^* \QQ_k(a_k) \big).
\end{eqnarray}
Let us first assume that $2 \le p \le 4$. Note that $\QQ_k(a_k)
\QQ_k(a_k)^*$ is not necessarily homogeneous. However, it is not
difficult to see that it is a polynomial in $L_{p/2}(\Al)$ of
degree $2d-1$. Therefore, it follows from Lemma \ref{Lemma-Type}
that
\begin{eqnarray*}
\Big\| \sum_{k=1}^n \QQ_k \big( \QQ_k(a_k) \QQ_k(a_k)^* \big)
\Big\|_{p/2} & \le & c d^4 \Big( \sum_{k=1}^n
\big\| \QQ_k(a_k) \QQ_k(a_k)^* \big\|_{p/2}^{p/2} \Big)^{\frac2p} \\
& = & c d^4 \Big( \sum_{k=1}^n \big\| \QQ_k(a_k) \big\|_p^p
\Big)^{\frac2p}.
\end{eqnarray*}
By (\ref{left}) and the triangle inequality, we deduce
\begin{eqnarray*}
\Big\| \sum_{k=1}^n \QQ_k(a_k) \QQ_k(a_k)^* \Big\|_{p/2}^{\frac12}
& \le & \Big\| \sum_{k=1}^n \Es \big( \QQ_k(a_k) \QQ_k(a_k)^*
\big) \Big\|_{p/2}^{\frac12} \\ & + & c d^2 \Big( \sum_{k=1}^n
\big\| \QQ_k(a_k) \big\|_p^p \Big)^{\frac1p}.
\end{eqnarray*}
Taking adjoints, we obtain a similar estimate for the last term of
(\ref{Equation-Start}). Hence, given any index $2 \le p \le 4$, we
have proved that the assertion holds with $\mathcal{C}_p(d) \le
c_0 d^4$ for some absolute constant $c_0$.

\vskip5pt

\noindent \textbf{Step 2.} Now we proceed by induction and assume
the assertion is proved in $L_{p/2}(\Al)$ with constant
$\mathcal{C}_{p/2}(d)$ for some $4 < p < \infty$. Of course, we
still have (\ref{Equation-Start}), (\ref{left}) and (\ref{right})
at our disposal. Thus, arguing as above it suffices to estimate
the term
 $$\Big\| \sum_{k=1}^n \QQ_k \big( \QQ_k(a_k) \QQ_k(a_k)^*
 \big) \Big\|_{p/2}^{\frac12}.$$
Let us write $x_k = \QQ_k(a_k) \QQ_k(a_k)^*$. As observed above,
we know that $x_k$ is a polynomial of degree $2d-1$. Hence, we may
use the projections $\Pi_\Al(p,s)$ from Remark
\ref{Remark-Proyecciones} and obtain the following inequality for
$x_{ks} = \Pi_\Al(p,s)(x_k)$
 $$\Big\| \sum_{k=1}^n \QQ_k(x_k) \Big\|_{p/2} \le
 \sum_{s=1}^{2d-1} \Big\| \sum_{k=1}^n \QQ_k(x_{ks})\Big\|_{p/2}.$$
By the induction hypothesis, we have
 $$\sum_{s=1}^{2d-1} \Big\| \sum_{k=1}^n \QQ_k(x_{ks}) \Big\|_{p/2}
 \le \sum_{s=1}^{2d-1} \mathcal{C}_{p/2}(s) \big( \mathrm{A}_s +
 \mathrm{B}_s + \mathrm{C}_s \big).$$
By Remark \ref{Remark-Diffp>2}, the first term on the right is
estimated by
\begin{eqnarray*}
\mathrm{A}_s & = & \Big( \sum_{k=1}^n \big\| \QQ_k(x_{ks})
\big\|_{p/2}^{p/2} \Big)^{\frac2p} \\ & \le & c s \Big(
\sum_{k=1}^n \big\| \Pi_\Al(p,s) (x_k) \big\|_{p/2}^{p/2}
\Big)^{\frac2p} \le c s^2 \Big( \sum_{k=1}^n \big\| \QQ_k(a_k)
\big\|_p^p \Big)^{\frac2p}.
\end{eqnarray*}
The second term is given by
 $$\mathrm{B}_s = \Big\| \sum_{k=1}^n\Es \big(
 \QQ_k(x_{ks})^* \QQ_k(x_{ks}) \big)\Big\|_{p/4}^{\frac12}.$$
Using $x_k = \sum_s x_{ks}$, freeness and (\ref{right}), we have
for all $1 \le s \le 2d-1$
\begin{eqnarray*}
\Es \big( \QQ_k(x_{ks})^* \QQ_k(x_{ks}) \big) & \le & \summ_r \Es
\big( \QQ_k(x_{kr})^* \QQ_k(x_{kr}) \big)
\\ & = & \summ_{q,r} \Es
\big( \QQ_k(x_{kq})^* \QQ_k(x_{kr}) \big) \\ & = & \Es \big(
\QQ_k(x_k)^* \QQ_k(x_k) \big) \\ & = & \Es \Big( \big( x_k -
\Es(x_k) \big)^* \big( x_k - \Es(x_k) \big) \Big) \\ & = & \Es
(x_k^* x_k) - \Es (x_k)^* \Es(x_k) \le \Es (x_k^* x_k).
\end{eqnarray*}
 Then we apply \cite[Lemma 5.2]{JX} and then obtain
\begin{eqnarray*}
\mathrm{B}_s & \le & \big\| \sum_{k=1}^n \Es (x_k^*x_k)
\big\|_{p/4}^{\frac12} = \big\| \sum_{k=1}^n \Es \big|
\QQ_k(a_k)^* \big|^4 \big\|_{p/4}^{\frac12} \\ & \le & \big\|
\sum_{k=1}^n \Es \big( \QQ_k(a_k) \QQ_k(a_k)^* \big)
\big\|_{p/2}^{\frac{p-4}{2p-4}} \big( \sum_{k=1}^n \big\|
\QQ_k(a_k) \big\|_p^p \big)^{\frac{2}{2p-4}}.
\end{eqnarray*}
The same estimate holds for $\mathrm{C}_s$. Now, by homogeneity we
may assume that
 $$\Big( \sum_{k=1}^n \big\| \QQ_k(a_k) \big\|_p^p\Big)^{\frac1p}
 + \Big\| \sum_{k=1}^n \mathsf{E} [\QQ_k(a_k)^* \QQ_k(a_k) ]
 \Big\|_{\frac p2}^{\frac12} + \Big\|\sum_{k=1}^n \Es
 [ \QQ_k(a_k) \QQ_k(a_k)^* ]\Big\|_{\frac p2}^{\frac12}=1.$$
 Then combining the inequalities so far obtained, we deduce
 $$\sum_{s=1}^{2d-1} \Big\| \sum_{k=1}^n \QQ_k(x_{ks}) \Big\|_{p/2}
 \le \sum_{s=1}^{2d-1} \mathcal{C}_{p/2}(s) \big( 2 + c s^2\big).$$
Chasing through the inequalities above, we obtain the estimate
 $$\mathcal{C}_p(d) \le \sqrt{c} \, d^{\frac72}
 \sqrt{\mathcal{C}_{\frac p2}(d)},$$
for some absolute constant $c$. Taking $c$ big enough so that $c_0
\le c$ and recalling that $\mathcal{C}_p(d) \le c_0 d^4 \le c d^7$
for $2 \le p \le 4$, it turns out that the growth of the constant
$\mathcal{C}_p(d)$ as $d \to \infty$ is controlled by $c d^7$.
This proves the assertion. \fin

\begin{remark} \label{Remark-Free-Independence}
\emph{A noncommutative analogue of Rosenthal's inequality for
general von Neumann algebras (non necessarily free products) was
obtained in \cite{JX,JX2}, see also \cite{X} for the proof and the
notion of noncommutative independence employed in it. As we have
pointed out in the Introduction, recalling that freeness implies
this notion of independence, Theorem A for $d=1$ and $2 \le p <
\infty$ follows from the noncommutative Rosenthal inequality.
However, the constants in \cite{JX,JX2} are not uniformly bounded
as $p \to \infty$, in sharp contrast with Theorem A. Similarly,
one could try to derive Theorem A for $d \ge 1$ and $2 \le p <
\infty$ by proving that $\QQ_1(a_1), \QQ_2(a_2), \ldots,
\QQ_n(a_n)$ are independent in the sense of \cite{JX2}.
Nevertheless, this alternative approach to Theorem A would provide
constants depending on $p$, rather than on $d$.}
\end{remark}

Since any $a \in \mathbf{P}_\Al(p,d)$ satisfies
 $$a = \sum_{k=1}^n\LL_k(a) = \sum_{k=1}^n \RR_k(a),$$
the following result characterizes the $L_p$ norm of all
homogeneous free polynomials.

\begin{corollary} \label{Corollary-Voiculescu1}
If $2 \le p \le \infty$ and $a_1, a_2, \ldots, a_n \in
\mathbf{P}_\Al(p,d)$, we have
\begin{eqnarray*}
\Big\| \sum_{k=1}^n \LL_k(a_k) \Big\|_p & \sim_{c d^{7}} & \Big\|
\sum_{k=1}^n \LL_k(a_k)^* \LL_k(a_k) \Big\|_{p/2}^{\frac12} \ +
\Big\| \sum_{k=1}^n \Es \big( \LL_k(a_k) \LL_k(a_k)^* \big)
\Big\|_{p/2}^{\frac12}, \\ \Big\| \sum_{k=1}^n \RR_k(a_k) \Big\|_p
& \sim_{c d^{7}} & \Big\| \sum_{k=1}^n \RR_k(a_k) \RR_k(a_k)^*
\Big\|_{p/2}^{\frac12} + \Big\| \sum_{k=1}^n \Es \big(
\RR_k(a_k)^* \RR_k(a_k) \big) \Big\|_{p/2}^{\frac12}.
\end{eqnarray*}
\end{corollary}

\dem By (\ref{Eq-Est33}) we have
\begin{eqnarray*}
\lefteqn{\Big\| \sum_{k=1}^n \RR_k(a_k) \RR_k(a_k)^*
\Big\|_{p/2}^{\frac12}} \\ & = & \Big\| \int_{\Omega}
\sum_{i,j=1}^n \varepsilon_i \RR_i(a_i) \varepsilon_j \RR_j(a_j)^*
d \varepsilon \Big\|_{p/2}^{\frac12} \\ & \le & \Big(
\int_{\Omega} \Big\| \sum_{i=1}^n \varepsilon_i \RR_i(a_i)
\Big\|_p \Big\| \sum_{j=1}^n \varepsilon_j \RR_j(a_j)^* \Big\|_p d
\varepsilon \Big)^{\frac12} \le c d \, \Big\| \sum_{k=1}^n
\RR_k(a_k) \Big\|_p.
\end{eqnarray*}
On the other hand, by freeness
 $$\sum_{k=1}^n \Es \big(\RR_k(a_k)^* \RR_k(a_k) \big)
 = \Es \Big( \big( \sum_{i=1}^n\RR_i(a_i) \big)^*
 \big( \sum_{j=1}^n \RR_j(a_j) \big) \Big).$$
Therefore, by the contractivity of $\Es$
 $$\Big\|\sum_{k=1}^n \Es \big( \RR_k(a_k)^* \RR_k(a_k) \big)
 \Big\|_{p/2}^{\frac12} \le \Big\| \sum_{k=1}^n
 \RR_k(a_k)\Big\|_p.$$
This gives the lower estimate.

\vskip5pt

For the upper estimate we assume that $2 \le p < \infty$, since
the case $p=\infty$ was already proved in Proposition
\ref{Proposition-Sum-L&R}. Now we use the second inequality stated
in Corollary \ref{Corollary-Sum-L&R}
 $$\Big\| \sum_{k=1}^n \RR_k(a_k) \Big\|_p \le c d^2
 \left( \Big\|\sum_{k=1}^n \RR_k(a_k) \RR_k(a_k)^*
 \Big\|_{p/2}^{\frac12} +\Big\| \sum_{k=1}^n \RR_k(a_k)^* \RR_k(a_k)
 \Big\|_{p/2}^{\frac12}\right).$$
On the other hand, it is clear that
\begin{equation} \label{Eq-E-Q}
\sum_{k=1}^n \RR_k(a_k)^* \RR_k(a_k) = \sum_{k=1}^n \Es \big(
\RR_k(a_k)^* \RR_k(a_k) \big) + \sum_{k=1}^n \QQ_k \big(
\RR_k(a_k)^* \RR_k(a_k) \big).
\end{equation}
Hence, it suffices to estimate the last term on the right. This
part of the proof is similar to the corresponding one of the proof
of Theorem A.  Again, we observe that $x_k = \RR_k(a_k)^*
\RR_k(a_k)$ is no longer homogeneous but a polynomial of degree $
\le 2d$. Our argument for this term depends on the value of $p$.

\vskip5pt

\noindent \textbf{Step 1.} If $2 \le p \le 4$, we apply Lemma
\ref{Lemma-Type} and obtain
\begin{eqnarray*}
\Big\| \sum_{k=1}^n \QQ_k(x_k) \Big\|_{p/2}^{\frac12} & \le & c
d^2 \Big( \sum_{k=1}^n \|x_k\|_{p/2}^{p/2} \Big)^{\frac1p}
\\ & = & c d^2 \Big( \sum_{k=1}^n \big\| \RR_k(a_k) \big\|_{p}^{p}
\Big)^{\frac1p} \le  c d^2 \Big\| \sum_{k=1}^n \RR_k(a_k)
\RR_k(a_k)^* \Big\|_{p/2}^{\frac12},
\end{eqnarray*}
where the last inequality holds for $2 \le p \le \infty$ and
follows by complex interpolation. Hence, in the case $2 \le p \le
4$, we have proved the upper estimate with constant $c d^4$.

\vskip5pt

\noindent \textbf{Step 2.} If $4 < p < \infty$, we take $x_{ks} =
\Pi_\Al(p,s)(x_k)$ and write
 $$\Big\| \sum_{k=1}^n \QQ_k(x_k) \Big\|_{p/2}^{\frac12}
 \le \Big(\sum_{s=1}^{2d} \Big\| \sum_{k=1}^n \QQ_k(x_{ks})
 \Big\|_{p/2}\Big)^{\frac12}
 \le \sqrt{2d} \max_{1 \le s \le 2d}
 \Big\|\sum_{k=1}^n \QQ_k(x_{ks}) \Big\|_{p/2}^{\frac12}.$$
By Theorem \ref{Theorem-Voiculescu}, we have
\begin{eqnarray*}
\Big\| \sum_{k=1}^n \QQ_k(x_{ks}) \Big\|_{p/2} & \sim_{c s^7} &
\Big( \sum_{k=1}^n \big\| \QQ_k(x_{ks}) \big\|_{p/2}^{p/2}
\Big)^{\frac2p} \\ & + & \Big\| \sum_{k=1}^n \mathsf{E} \big(
\QQ_k(x_{ks})^* \QQ_k(x_{ks}) \big) \Big\|_{p/4}^{\frac12} \\ & +
& \Big\| \sum_{k=1}^n \Es \big( \QQ_k(x_{ks}) \QQ_k(x_{ks})^*
\big) \Big\|_{p/4}^{\frac12} = \mathrm{A}_s + \mathrm{B}_s +
\mathrm{C}_s.
\end{eqnarray*}
These terms are estimated as in the proof of Theorem
\ref{Theorem-Voiculescu} (Step 2)
\begin{eqnarray*}
\mathrm{A}_s & \le & c s^2 \Big( \sum_{k=1}^n \big\| \RR_k(a_k)
\big\|_{p}^{p} \Big)^{\frac2p}.
\end{eqnarray*}
Similarly, we have
 $$\max (\mathrm{B}_s, \mathrm{C}_s)
 \le \Big\|\sum_{k=1}^n \Es \big( \RR_k(a_k)^* \RR_k(a_k)\big)
 \Big\|_{p/2}^{\frac{p-4}{2p-4}} \Big( \sum_{k=1}^n \big\|
 \RR_k(a_k) \big\|_p^p \Big)^{\frac{2}{2p-4}}.$$
On the other hand, by homogeneity we may assume that
 $$\Big\| \sum_{k=1}^n \RR_k(a_k) \RR_k(a_k)^*
 \Big\|_{p/2}^{\frac12} +\Big\| \sum_{k=1}^n \Es \big( \RR_k(a_k)^*
 \RR_k(a_k) \big)\Big\|_{p/2}^{\frac12} = 1.$$
Using the estimates above and
 $$\Big( \sum_{k=1}^n \big\| \RR_k(a_k) \big\|_p^p \Big)^{\frac1p}
 \le\Big\|\sum_{k=1}^n\RR_k(a_k)\RR_k(a_k)^*\Big\|_{p/2}^{\frac12},$$
we obtain
 $$\Big\|\sum_{k=1}^n \QQ_k(x_k) \Big\|_{p/2}^{\frac12}
 \le \sqrt{2d}\max_{1\le s \le 2d}\Big[ c s^7\big( 2 + c s^2\big)
 \Big]^{\frac12},$$
for $4 < p < \infty$. Therefore, by Corollary
\ref{Corollary-Sum-L&R} and (\ref{Eq-E-Q}) we find
 $$\Big\| \sum_{k=1}^n \RR_k(a_k) \Big\|_p \le c d^7.$$
This and Step 1 yield the assertion for $\RR_k$'s. For $\LL_k$'s
we take adjoints. \fin

\begin{corollary} \label{Corollary-Voiculescu2}
If $2 \le p \le \infty$ and $a \in \mathbf{P}_\Al(p,d)$, we have
 $$\|a\|_p \sim_{c d^{14}} \Big\| \sum_{i,j=1}^n
 e_{ij} \ten \LL_i \RR_j(a) \Big\|_{S_p^n(L_p(\Al))}
 + \big\|\Es(aa^*)^{\frac12} \big\|_{L_p(\Be)}
 + \big\| \Es(a^*a)^{\frac12}\big\|_{L_p(\Be)}.$$
\end{corollary}

\dem We use $a = \sum_{k=1}^n \RR_k(a)$ and Corollary
\ref{Corollary-Voiculescu1}
 $$\|a\|_p\sim_{c d^7}\Big\|\sum_{k=1}^n\RR_k(a)\RR_k(a)^*
 \Big\|_{p/2}^{\frac12} + \Big\|\sum_{k=1}^n \Es
 \big( \RR_k(a)^* \RR_k(a) \big)\Big\|_{p/2}^{\frac12}
 = \mathrm{A} + \mathrm{B}.$$
To estimate $\mathrm{A}$ we use Corollary
\ref{Corollary-Voiculescu1} for the $\LL_k$'s
\begin{eqnarray*}
\lefteqn{\Big\| \sum_{k=1}^n \RR_k(a) \RR_k(a)^*
\Big\|_{p/2}^{\frac12}} \\ & = & \Big\| \sum_{k=1}^n e_{1k} \ten
\RR_k(a) \Big\|_{S_p^n(L_p(\Al))}
\\ & \sim_{c d^7} & \Big\| \sum_{i=1}^n e_{i1} \ten \LL_i \Big(
\sum_{j=1}^n e_{1j} \ten \RR_j(a) \Big)
\Big\|_{S_p^{n^2}(L_p(\Al))} \\ & + & \Big\| \sum_{i=1}^n \Es
\Big( \big( \sum_{j=1}^n e_{1j} \ten \LL_i \RR_j(a) \big) \big(
\sum_{j=1}^n e_{1j} \ten \LL_i \RR_j(a) \big)^* \Big)
\Big\|_{S_{p/2}^n(L_{p/2}(\Be))}^{\frac12} \\ & = & \Big\|
\sum_{i,j=1}^n e_{ij} \ten \LL_i \RR_j(a)
\Big\|_{S_p^{n}(L_p(\Al))}
\\ & + & \Big\| e_{11} \ten \sum_{i,j=1}^n \Es \big( (\LL_i \RR_j(a))
(\LL_i \RR_j(a))^* \big)
\Big\|_{S_{p/2}^n(L_{p/2}(\Be))}^{\frac12}
\\ & = & \Big\| \sum_{i,j=1}^n e_{ij} \ten \LL_i \RR_j(a)
\Big\|_{S_p^{n}(L_p(\Al))} + \big\| \Es (aa^*)^{\frac12}
\big\|_{L_p(\Be)}.
\end{eqnarray*}
On the other hand, it is clear that $$\mathrm{B} = \Big\|
\sum_{k=1}^n \Es \big( \RR_k(a)^* \RR_k(a) \big)
\Big\|_{p/2}^{\frac12} = \big\| \Es \big( a^*a \big)^{\frac12}
\big\|_{L_p(\Be)}.$$ Thus, since we have used equivalences at each
step, the proof is completed. \fin

\begin{remark}
\emph{By decomposing a free polynomial of degree $d$ into its
homogeneous parts, we automatically obtain trivial generalizations
of Theorem A and Corollaries \ref{Corollary-Voiculescu1} and
\ref{Corollary-Voiculescu2} for non-homogeneous free polynomials
of a fixed degree $d$. Most of the forthcoming results in this
paper are susceptible of this kind of generalization.}
\end{remark}

\section{A length-reduction formula}
\label{Section3}

In this section we prove a length-reduction formula for
polynomials in the free product. One more time, our standard
assumptions are that $\Al = *_{\Be} \mathsf{A}_k$ where $1 \le k
\le n$, $\Be$ is equipped with a \emph{n.f.} state $\varphi$ which
induces a \emph{n.f.} state $\phi = \varphi \circ \Es$ on $\Al$
and $\Es: \Al \to \Be$ is a \emph{n.f.} conditional expectation.
As usual, $d_\phi$ denotes the density of the state $\phi$. We
will need some preliminary facts on certain module maps. First,
given $2 \le p \le \infty$, we define on $\Al \ten_\Be L_p(\Al)$
the $L_{p/2}(\Al)$-valued inner product
 $$\big\langle \big\langle x_1 \ten y_1,\; x_2 \ten y_2 \big\rangle
 \big\rangle = y_1^* \Es(x_1^*x_2) y_2.$$
This allows us to define $L_p^c(\Al \ten_\Be \Al, \Es)$ and
$L_p^r(\Al \ten_\Be \Al, \Es)$ as the completion of the space $\Al
\ten_\Be L_p(\Al)$ with respect to the norms
\begin{eqnarray*}
\|z\|_{L_p^c(\Al \ten_\Be \Al, \Es)} & = & \big\| \langle \langle
z,\; z \rangle \rangle \big\|_{L_{p/2}(\Al)}^{\frac12}, \\
\|z\|_{L_p^r(\Al \ten_\Be \Al, \Es)} & = & \big\| \langle \langle
z^*,\; z^* \rangle \rangle \big\|_{L_{p/2}(\Al)}^{\frac12}.
\end{eqnarray*}
Let $C_\infty(\Be)$ be the column subspace of the $\Be$-valued
Schatten class $S_\infty(\Be)$
 $$C_\infty(\Be) = \Big\{ \summ_ke_{k1} \ten b_k \in
 \mathcal{B}(\ell_2) \ten_{\mathrm{min}} \Be\Big\}.$$
By \cite{Pas}, there exists a normal right $\Be$-module map $u:
\Al \to C_\infty(\Be)$ satisfying
\begin{equation} \label{Eq-u-E}
\Es(x^*y) = \sum_{k=1}^\infty u_k(x)^* u_k(y) = u(x)^* u(y) \quad
\mbox{for all} \quad x,y \in \Al,
\end{equation}
where $u_k$ stands for the $k$-th coordinate of $u$. Note that,
according to \cite{J1,JS}, this map canonically extends to
$L_p(\Al)$. On the other hand, recalling that amalgamation gives
$C_\infty(\Be) \ten_\Be L_p(\Al) = C_p (L_p(\Al))$, we have an
isometry
\begin{equation} \label{Eq-Auxxx1}
\widehat{u} = u \ten id_{L_p(\Al)}: L_p^c(\Al \ten_\Be \Al, \Es)
\to C_p(L_p(\Al)).
\end{equation}
Indeed, note that
 $$\big( \widehat{u} (x_1 \ten y_1) \big)^* \big( \widehat{u}
 (x_2 \ten y_2) \big) = \sum_{k=1}^\infty y_1^* u_k(x_1)^* u_k(x_2)
 y_2 = y_1^* \Es \big( x_1^* x_2 \big) y_2.$$
Thus, linearity gives
 $$\big\| \widehat{u}(z) \big\|_{C_p(L_p(\Al))} =
 \big\| \langle \langle z,\; z \rangle \rangle\big\|_{L_{p/2}(\Al)}^{\frac12}
 = \|z\|_{L_p^c(\Al \ten_{\Be} \Al,\Es)}.$$
A similar argument holds in the row case and by \cite[Proposition
2.8]{J1} we deduce

\begin{lemma} \label{Lemma-Complementation2}
Let $\Al$ and $\Be$ be as above. Then,
 $$L_p^r(\Al \ten_\Be \Al, \Es) \quad \mbox{and}
 \quad L_p^c(\Al\ten_\Be \Al, \Es)$$
are contractively complemented in the space $S_p(L_p(\Al))$ for
any $2 \le p \le \infty$.
\end{lemma}

\noindent In what follows, $\Lambda$ will always denote a finite
index set.

\begin{lemma} \label{Lemma-Complementation3}
If $2 \le p \le \infty$, the space $$\mathcal{W}_p =
\overline{\Big\{ \sum_{\alpha \in \Lambda} \sum_{k=1}^n
x_k(\alpha) \ten w_k(\alpha) \, \big| \ x_k(\alpha) \in
\bubl_k\Big\}}$$ is contractively complemented in $L_p^r(\Al
\ten_\Be \Al, \Es)$ as well as in $L_p^c(\Al \ten_\Be \Al, \Es)$.
\end{lemma}

\dem By definition, $L_p^c(\Al \ten_\Be \Al, \Es)$ is the closure
of
 $$\sum_{\alpha \in \Lambda} x(\alpha) \ten w(\alpha),$$
where $x(\alpha) \in \Al$ and $w(\alpha) \in L_p(\Al)$. Let us
recall the notation $\Pi_\Al(p,d)$, introduced in Remark
\ref{Remark-Proyecciones} for the projection from $L_p(\Al)$ onto
the homogeneous polynomials of degree $d$. Then we clearly have
 $$x(\alpha) = \Es(x(\alpha)) + \Pi_\Al(p,1)(x(\alpha))
 +\sum_{d\ge2} \Pi_\Al(p,d)(x(\alpha)) = x(\alpha,0)
 + x(\alpha,1) +x(\alpha,2).$$
Now we define
\begin{eqnarray*}
\mathrm{A} & = & \sum_{\alpha \in \Lambda} x(\alpha,1) \ten
w(\alpha), \\ \mathrm{B} & = & \sum_{\alpha \in \Lambda} 1_\Al
\ten x(\alpha,0) w(\alpha) + \sum_{\alpha \in \Lambda} x(\alpha,2)
\ten w(\alpha).
\end{eqnarray*}
Note that $\sum_\alpha x(\alpha) \ten w(\alpha) = \mathrm{A} +
\mathrm{B}$ and $\mathrm{A} \in \mathcal{W}_p$. On the other hand,
by freeness
 $$\big\langle \big\langle \mathrm{A} + \mathrm{B},\;
 \mathrm{A} + \mathrm{B} \big\rangle \big\rangle
 = \big\langle\big\langle \mathrm{A},\; \mathrm{A}
 \big\rangle \big\rangle +\big\langle \big\langle \mathrm{B},\;
 \mathrm{B} \big\rangle\big\rangle.$$
Therefore, by positivity
\begin{eqnarray*}
\|\mathrm{A}\|_{L_p^c(\Al \ten_\Be \Al, \Es)}^2 & = & \big\|
\langle \langle \mathrm{A},\; \mathrm{A} \rangle \rangle
\big\|_{p/2} \\ & \le & \big\| \langle \langle \mathrm{A} +
\mathrm{B},\; \mathrm{A} + \mathrm{B} \rangle \rangle \big\|_{p/2}
\\ & = & \Big\| \sum_{\alpha \in \Lambda} x(\alpha) \ten w(\alpha)
\Big\|_{L_p^c(\Al \ten_\Be \Al, \Es)}^2.
\end{eqnarray*}
By continuity, we find a contractive projection $L_p^c(\Al
\ten_\Be \Al, \Es) \to \mathcal{W}_p$ for any given index $2 \le p
\le \infty$. Obviously, the argument above also works for
$L_p^r(\Al \ten_\Be \Al, \Es)$. \fin

\begin{lemma} \label{Lemma-Complementation4}
If $2 \le p \le \infty$, the space
 $$\mathcal{Z}_{p,d}^r =\overline{\Big\{ \sum_{\alpha \in \Lambda}
 \sum_{k=1}^nx_k(\alpha) \ten w_k(\alpha) \in \mathcal{W}_p \,
 \big| \ w_k(\alpha) \in \mathbf{P}_\Al(p,d), \ \RR_k(w_k(\alpha))
 = 0\Big\}}$$
is complemented in $L_p^r(\Al \ten_\Be \Al, \Es)$. Similarly, the
space
 $$\mathcal{Z}_{p,d}^c
 = \overline{\Big\{ \sum_{\alpha \in \Lambda}
 \sum_{k=1}^n x_k(\alpha) \ten w_k(\alpha) \in \mathcal{W}_p \,
 \big| \ w_k(\alpha) \in \mathbf{P}_\Al(p,d), \ \LL_k(w_k(\alpha))
 = 0 \Big\}}$$
is complemented in $L_p^c(\Al \ten_\Be \Al, \Es)$. In both cases,
the projection constant is $\le c d^2$.
\end{lemma}

\dem Both complementation results can be proved using the same
arguments. Thus, we only prove the second assertion. According to
Lemma \ref{Lemma-Complementation3}, it suffices to check that
$\mathcal{Z}_{p,d}^c$ is complemented (with projection constant
$\le c d^2$) in $\mathcal{W}_p$ equipped with the norm inherited
from
 $$L_p^c(\Al \ten_\Be \Al, \Es).$$
To that aim, we consider the intermediate space
 $$\mathcal{W}_{p,d} = \overline{\Big\{ \sum_{\alpha \in \Lambda}
 \sum_{k=1}^n x_k(\alpha) \ten w_k(\alpha) \in \mathcal{W}_p \,
 \big| \ w_k(\alpha) \in \mathbf{P}_\Al(p,d) \Big\}}.$$
$\mathcal{W}_{p,d}$ is complemented in $\mathcal{W}_p$ with
constant $4 d$. Indeed, using one more time the projection
$\Pi_\Al(p,d)$ onto the $d$-homogeneous polynomials, we write
$w_{kd}(\alpha)$ for $\Pi_\Al(p,d)(w_k(\alpha))$ and obtain from
Lemma \ref{Lemma-Complementation2} and the discussion preceding it
\begin{eqnarray*}
\lefteqn{\Big\| \sum_{k,\alpha} x_k(\alpha) \ten w_{kd}(\alpha)
\Big\|_{L_p^c(\Al \ten_\Be \Al, \Es)}} \\ & = & \Big\|
\sum_{j=1}^\infty e_{j1} \ten \sum_{k,\alpha} u_j(x_k(\alpha))
w_{kd}(\alpha) \Big\|_{C_p(L_p(\Al))} \\ & = & \Big\|
\sum_{j=1}^\infty e_{j1} \ten \Pi_\Al(p,d) \Big( \sum_{k,\alpha}
u_j(x_k(\alpha)) w_{k}(\alpha) \Big) \Big\|_{C_p(L_p(\Al))} \\ &
\le & \big\| id_{C_p} \ten \Pi_\Al(p,d)
\big\|_{\mathcal{B}(C_p(L_p(\Al)))} \Big\| \sum_{j=1}^\infty
e_{j1} \ten \sum_{k,\alpha} u_j(x_k(\alpha)) w_{k}(\alpha)
\Big\|_{C_p(L_p(\Al))}.
\end{eqnarray*}
On the other hand, combining Remarks \ref{Remark-Operator-Valued}
and \ref{Remark-Proyecciones} we deduce that $\Pi_\Al(p,d)$ is a
completely bounded map on $L_p(\Al)$ with cb-norm less than or
equal to $4d$. Therefore, we deduce our claim
 $$\Big\| \sum_{k,\alpha} x_k(\alpha) \ten w_{kd}(\alpha)
 \Big\|_{L_p^c(\Al \ten_\Be \Al, \Es)} \le 4d \, \Big\|
 \sum_{k,\alpha} x_k(\alpha) \ten w_{k}(\alpha) \Big\|_{L_p^c(\Al
 \ten_\Be \Al, \Es)}.$$
It remains to see that $\mathcal{Z}_{p,d}^c$ is complemented (with
projection constant less than or equal to $c d$) in
$\mathcal{W}_{p,d}$ with the norm inherited from $L_p^c(\Al
\ten_\Be \Al, \Es)$. In other words, we are interested in proving
the following inequality
 $$\Big\| \sum_{k,\alpha} x_k(\alpha) \ten
 (id_\Al - \LL_k)(w_{kd}(\alpha)) \Big\|_{L_p^c(\Al \ten_\Be \Al,
 \Es)} \! \le c d \, \Big\| \sum_{k,\alpha} x_k(\alpha) \ten
 w_{kd}(\alpha) \Big\|_{L_p^c(\Al \ten_\Be \Al, \Es)}.$$
However, this follows from Lemma
\ref{Lemma-Complementation-Boundedness}, Remark
\ref{Remark-Diffp>2} and triangle inequality. \fin

\begin{remark} \label{Remark-Replacing-Density}
\emph{In our definition of the spaces $L_p^r(\Al \ten_\Be \Al,
\Es)$ and $L_p^c(\Al \ten_\Be \Al, \Es)$ as well as in Lemmas
\ref{Lemma-Complementation3} and \ref{Lemma-Complementation4}, we
have used tensors $x \ten w$ with $x \in \Al$ and $w \in
L_p(\Al)$. Note that, according to the definition of the inner
product $\langle \langle \ , \, \rangle \rangle$, it is relevant
to distinguish between the first and second components of these
tensors. However, in some forthcoming results (see e.g. the proof
of Lemma \ref{Lemma-Prelim1} below) we shall need to work with
tensors $x \ten w$ where $x \in L_p(\Al)$ and $w \in \Al$. Thus,
we have to understand which element of $L_p^r(\Al \ten_\Be \Al,
\Es)$ or $L_p^c(\Al \ten_\Be \Al, \Es)$ do we mean when writing $x
\ten w$. Let us consider a sequence $(x_n)_{n \ge 1}$ in $\Al$
such that
 $$x_n d_\phi^{\frac1p} \to x \quad \mbox{as}
 \quad n \to \infty$$ in $L_p(\Al)$. Then we set
  $$x \ten w =\lim_{n \to \infty} x_n \ten d_\phi^{\frac1p} w.$$
To make sure our definition makes sense, we must see that the
sequence on the right converges in the norms of $L_p^r(\Al
\ten_\Be \Al, \Es)$ and $L_p^c(\Al \ten_\Be \Al, \Es)$. Let us see
this for the first space, the other follows in the same way. By
completeness, it suffices to show that we have a Cauchy sequence.
This easily follows since}
\begin{eqnarray*}
\Big\| (x_n - x_m) \ten d_\phi^{\frac1p} w \Big\|_{L_p^r(\Al
\ten_\Be \Al, \Es)} & = & \Big\| w^* d_{\phi}^{\frac1p} \mathsf{E}
\big( (x_n-x_m)^* (x_n - x_m) \big) d_\phi^{\frac1p} w
\Big\|_{L_{p/2}(\Al)}^{\frac12} \\ & \le & \|w\|_{\Al} \Big\|
d_{\phi}^{\frac1p} \mathsf{E} \big( (x_n-x_m)^* (x_n - x_m) \big)
d_\phi^{\frac1p} \Big\|_{L_{p/2}(\Be)}^{\frac12} \\ & \le &
\|w\|_{\Al} \big\| (x_n - x_m) d_\phi^{\frac1p} \big\|_{L_p(\Al)}
\end{eqnarray*}
\emph{and the right hand side converges to $0$ as $n,m \to
\infty$.}
\end{remark}

\subsection{Preliminary estimates}

This paragraph is devoted to some necessary estimates that will be
used below. In the following we shall use the notation already
defined in the Introduction
\begin{eqnarray*}
\Big\| \sum_{k, \alpha} b_k(\alpha) \big\langle a_k(\alpha) |
\Big\|_p & = & \Big\| \Big( \sum_{i,j,\alpha,\beta} b_i(\alpha)
\Es \big( a_i(\alpha) a_j(\beta)^* \big) b_j(\beta)^*
\Big)^{\frac12} \Big\|_p, \\ \Big\| \sum_{k, \alpha} | a_k(\alpha)
\big\rangle b_k(\alpha) \Big\|_p & = & \Big\| \Big(
\sum_{i,j,\alpha,\beta} b_i(\alpha)^* \Es \big( a_i(\alpha)^*
a_j(\beta) \big) b_j(\beta) \Big)^{\frac12} \Big\|_p.
\end{eqnarray*}
In other words,
\begin{eqnarray*} \Big\| \sum_{k, \alpha}
b_k(\alpha) \big\langle a_k(\alpha) | \Big\|_p & = & \Big\|
\sum_{k, \alpha} a_k(\alpha) \ten b_k(\alpha) \Big\|_{L_p^r(\Al
\ten_\Be \Al, \Es)}, \\ \Big\| \sum_{k, \alpha} | a_k(\alpha)
\big\rangle b_k(\alpha) \Big\|_p & = & \Big\| \sum_{k, \alpha}
a_k(\alpha) \ten b_k(\alpha) \Big\|_{L_p^c(\Al \ten_\Be \Al,
\Es)}.
\end{eqnarray*}

\begin{lemma} \label{Lemma-Prelim1}
Let $2 \le p,q \le \infty$ be two indices related by $1/2 = 1/p +
1/q$. Let $x_k(\alpha)$ be a mean-zero element in $\mathsf{A}_k$
for each $1 \le k \le n$ and $\alpha$ running over a finite set
$\Lambda$. Let $w_k(\alpha) \in \mathbf{P}_\Al(d)$ for some $d \ge
0$ and satisfying $\RR_k(w_k(\alpha)) = 0$ for all $1 \le k \le n$
and every $\alpha \in \Lambda$. Then
\begin{eqnarray*}
\Big\| \sum_{k, \alpha} w_k(\alpha) \mathsf{L}_k x_k(\alpha)
d_\phi^{\frac1p} \Big\|_{\mathcal{B}(L_q(\Al), L_2(\Al))} & \le &
\Big\| \sum_{k, \alpha} | w_k(\alpha) \big\rangle x_k(\alpha)
d_\phi^{\frac1p} \Big\|_p, \\ \Big\| \sum_{k, \alpha} w_k(\alpha)
(1 - \mathsf{L}_k) x_k(\alpha) d_\phi^{\frac1p}
\Big\|_{\mathcal{B}(L_q(\Al), L_2(\Al))} & \le & c d^{2} \Big\|
\sum_{k, \alpha} w_k(\alpha) \big\langle x_k(\alpha)
d_\phi^{\frac1p} | \Big\|_p.
\end{eqnarray*}
\end{lemma}

\dem In what follows we use $x_k'(\alpha) = x_k(\alpha)
d_\phi^{1/p}$. Given $z \in L_q(\Al)$, we have $$h_k(\alpha) =
x_k'(\alpha) z \in L_2(\Al)$$ and the vector $\mathsf{L}_k
h_k(\alpha)$ is a linear combination of reduced words in
$L_2(\Al)$ starting with a mean-zero letter in $\mathsf{A}_k$.
Therefore, since $\RR_k(w_k(\alpha)) = 0$, the operator
$w_k(\alpha)$ acts on $\mathsf{L}_k h_k(\alpha)$ by tensoring from
the left. In particular, the $(d+1)$-th letter in the words of
$w_k(\alpha) \mathsf{L}_k h_k(\alpha)$ is always in $\mathsf{A}_k$
and the inequality below follows by freeness, (\ref{Eq-u-E}) and
the fact that $\mathsf{L}_k$ commutes with $\Be$
\begin{eqnarray*}
\Big\| \Big( \sum_{k, \alpha} w_k(\alpha) \mathsf{L}_k
x_k'(\alpha) \Big) (z) \Big\|_2^2 & = & \! \sum_{i,j,\alpha,\beta}
\!\! \mbox{tr}_{\Al} \Big( h_i(\alpha)^* \mathsf{L}_i
w_i(\alpha)^* w_j(\beta) \mathsf{L}_j h_j(\beta) \Big) \\ & = &
\sum_{k,\alpha,\beta} \mbox{tr}_{\Al} \Big( h_k(\alpha)^*
\mathsf{L}_k \Es \big( w_k(\alpha)^* w_k(\beta) \big) \mathsf{L}_k
h_k(\beta) \Big) \\ & \le & \sum_{k,\alpha, \beta} \mbox{tr}_\Al
\Big( h_k(\alpha)^* \Es \big( w_k(\alpha)^* w_k(\beta) \big)
h_k(\beta) \Big) \\ & = & \mbox{tr}_\Al \Big( z^*
\sum_{i,j,\alpha, \beta} x_i'(\alpha)^* \Es \big( w_i(\alpha)^*
w_j(\beta) \big) x_j'(\beta) z \Big) \\ & \le & \|z\|_q^2 \Big\|
\sum_{i,j, \alpha, \beta} x_i'(\alpha)^* \Es \big( w_i(\alpha)^*
w_j(\beta) \big) x_j'(\beta) \Big\|_{p/2}.
\end{eqnarray*}
This proves the first inequality.

\vskip5pt

Let us prove the second one. According to Lemma
\ref{Lemma-Complementation2}, we know that the spaces $L_p^r(\Al
\ten_\Be \Al, \Es)$ form an interpolation scale for $2 \le p \le
\infty$. Moreover, it follows from Lemma
\ref{Lemma-Complementation4} that the spaces $\mathcal{Z}_{p,d}^r$
also form (up to a constant $c d^2$) an interpolation scale for $2
\le p \le \infty$. Therefore, since (see Remark
\ref{Remark-Replacing-Density})
 $$\Big\| \sum_{k, \alpha} w_k(\alpha) \big\langle x_k(\alpha)
 d_\phi^{\frac1p} | \Big\|_p = \Big\| \sum_{k, \alpha} x_k(\alpha)
 d_\phi^{\frac1p} \ten w_k(\alpha) \Big\|_{\mathcal{Z}_{p,d}^r},$$
it suffices to see (by complex interpolation) that the assertion
holds when $p=2$ and $p=\infty$ with some constant not depending
on $d$. Let us use the same terminology for $x_k'(\alpha)$ as
above. If $p=2$, we have $q = \infty$ and the triangle inequality
gives
\begin{eqnarray*}
\Big\| \sum_{k, \alpha} w_k(\alpha) (1 - \mathsf{L}_k)
x_k'(\alpha) \Big\|_{\mathcal{B}(L_\infty(\Al), L_2(\Al))} & \le &
\Big\| \sum_{k, \alpha} w_k(\alpha) x_k'(\alpha)
\Big\|_{\mathcal{B}(L_\infty(\Al), L_2(\Al))} \\ & + & \Big\|
\sum_{k, \alpha} w_k(\alpha) \mathsf{L}_k x_k'(\alpha)
\Big\|_{\mathcal{B}(L_\infty(\Al), L_2(\Al))}.
\end{eqnarray*}
The first term equals
\begin{eqnarray*}
\Big\| \sum_{k, \alpha} w_k(\alpha) x_k'(\alpha) \Big\|_2 & = &
\Big( \sum_{i,j,\alpha,\beta} \mbox{tr}_\Al \Big[ w_i(\alpha)
x_i'(\alpha) x_j'(\beta)^* w_j(\beta)^* \Big] \Big)^{\frac12} \\ &
= & \Big( \sum_{i,j,\alpha,\beta} \mbox{tr}_\Al \Big[ w_i(\alpha)
\Es \big( x_i'(\alpha) x_j'(\beta)^* \big) w_j(\beta)^* \Big]
\Big)^{\frac12} \\ & = & \Big\| \sum_{k, \alpha} w_k(\alpha)
\big\langle x_k'(\alpha) | \Big\|_2.
\end{eqnarray*}
To estimate the second term, we use the first inequality proved in
this lemma
\begin{eqnarray*}
\Big\| \sum_{k, \alpha} w_k(\alpha) \mathsf{L}_k x_k'(\alpha)
\Big\|_{\mathcal{B}(L_\infty(\Al), L_2(\Al))}^2 & \le & \Big\|
\sum_{k, \alpha} | w_k(\alpha) \big\rangle x_k'(\alpha) \Big\|_2^2
\\ & = & \sum_{i,j,\alpha,\beta} \mbox{tr}_\Al \Big( x_i'(\alpha)^*
\Es \big( w_i(\alpha)^* w_j(\beta) \big) x_j'(\beta) \Big) \\ & =
& \sum_{i,j,\alpha,\beta} \mbox{tr}_\Al \Big( x_i'(\alpha)^*
w_i(\alpha)^* w_j(\beta) x_j'(\beta) \Big) \\ & = &
\sum_{i,j,\alpha,\beta} \mbox{tr}_\Al \Big( w_j(\beta) x_j'(\beta)
x_i'(\alpha)^* w_i(\alpha)^* \Big) \\ & = &
\sum_{i,j,\alpha,\beta} \mbox{tr}_\Al \Big( w_j(\beta) \Es \big(
x_j'(\beta) x_i'(\alpha)^* \big) w_i(\alpha)^* \Big) \\ & = &
\Big\| \sum_{k, \alpha} w_k(\alpha) \big\langle x_k'(\alpha) |
\Big\|_2^2.
\end{eqnarray*}
Therefore, we have proved that
 $$\Big\| \sum_{k, \alpha}w_k(\alpha)
 (1 - \mathsf{L}_k) x_k'(\alpha)\Big\|_{\mathcal{B}(L_\infty(\Al),
 L_2(\Al))} \le 2 \Big\|\sum_{k, \alpha} w_k(\alpha)
 \big\langle x_k'(\alpha) |\Big\|_2.$$
To prove the assertion for $p=\infty$ and $q=2$, we first note
that $$(1-\mathsf{L}_k) x_k(\alpha) = (1- \mathsf{L}_k)
x_k(\alpha) \mathsf{L}_k.$$ This implies
\begin{eqnarray*}
\lefteqn{\Big\| \sum_{k, \alpha} w_k(\alpha) (1- \mathsf{L}_k)
x_k(\alpha) \Big\|_{\mathcal{B}(L_2(\Al), L_2(\Al))}^2} \\ & = &
\Big\| \sum_{k, \alpha} w_k(\alpha) (1- \mathsf{L}_k) x_k(\alpha)
\mathsf{L}_k \Big\|_\infty^2 \\ & = & \Big\| \sum_{k, \alpha,
\beta} w_k(\alpha) (1- \mathsf{L}_k) x_k(\alpha) x_k(\beta)^* (1-
\mathsf{L}_k) w_k(\beta)^* \Big\|_\infty \\ & = & \Big\| \sum_{k,
\alpha, \beta} w_k(\alpha) (1- \mathsf{L}_k) \Es \big( x_k(\alpha)
x_k(\beta)^* \big) (1- \mathsf{L}_k) w_k(\beta)^* \Big\|_\infty \\
& \le & \Big\| \sum_{k, \alpha, \beta} w_k(\alpha) \Es \big(
x_k(\alpha) x_k(\beta)^* \big) w_k(\beta)^* \Big\|_\infty.
\end{eqnarray*}
Hence, we have seen that $$\Big\| \sum_{k, \alpha} w_k(\alpha) (1-
\mathsf{L}_k) x_k(\alpha) \Big\|_{\mathcal{B}(L_2(\Al), L_2(\Al))}
\le \Big\| \sum_{k, \alpha} w_k(\alpha) \big\langle x_k(\alpha) |
\Big\|_\infty.$$ This proves the assertion for $p=\infty$. The
general case follows by interpolation. \fin

\subsection{Proof of Theorems B and C}

Now we prove the second major result of this paper, a
length-reduction formula for homogeneous polynomials on free
random variables. As consequence, we extend the main results in
\cite{PP,RX}.

\vskip5pt

\demB The second reduction formula clearly follows from the first
one by taking adjoints. Thus, it suffices to prove the first
reduction formula. We begin by proving the upper estimate. If $1/p
+ 1/q = 1/2$, we have
\begin{eqnarray*}
\Big\| \sum_{k, \alpha} w_k(\alpha) x_k(\alpha) \Big\|_{L_p(\Al)}
& = & \Big\| \sum_{k, \alpha} w_k(\alpha) x_k(\alpha)
\Big\|_{\mathcal{B}(L_q(\Al), L_2(\Al))} \\ & \le & \Big\|
\sum_{k, \alpha} w_k(\alpha) \mathsf{L}_k x_k(\alpha)
\Big\|_{\mathcal{B}(L_q(\Al), L_2(\Al))} \\ & + & \Big\| \sum_{k,
\alpha} w_k(\alpha) (1 - \mathsf{L}_k) x_k(\alpha)
\Big\|_{\mathcal{B}(L_q(\Al), L_2(\Al))}.
\end{eqnarray*}
If we approximate $x_k(\alpha)$ by elements of the form
 $$z_k(\alpha) d_\phi^{\frac1p} \quad \mbox{with} \quad
 z_k(\alpha) \in \bubl_k,$$
the upper estimate follows from the inequalities in Lemma
\ref{Lemma-Prelim1}
 $$\Big\| \sum_{k, \alpha} w_k(\alpha) x_k(\alpha)
 \Big\|_{L_p(\Al)} \le \Big\| \sum_{k, \alpha} | w_k(\alpha)
 \big\rangle x_k(\alpha) \Big\|_p + c d^2 \Big\| \sum_{k, \alpha}
 w_k(\alpha) \big\langle x_k(\alpha) | \Big\|_p.$$
To prove the lower estimate we use the projection
 $$\Gamma_\Al(p,d): L_p(\Al) \to \mathbf{Q}_\Al(p,d)$$
which, according to Remark \ref{Remark-Proyecciones}, is bounded
by $2d+1$. Then we observe
 $$\sum_{i,j,\alpha,\beta} x_i(\alpha)^* \Es \big( w_i(\alpha)^*
 w_j(\beta) \big) x_j(\beta) = \Gamma_\Al(p/2,2)(a^*a)$$
for $a = \sum_{k, \alpha} w_k(\alpha) x_k(\alpha) \in L_p(\Al)$.
In particular, we deduce
 $$\Big\| \sum_{k, \alpha} | w_k(\alpha) \big\rangle x_k(\alpha)
 \Big\|_p = \big\| \Gamma_\Al(p/2,2)(a^*a) \big\|_{p/2}^{\frac12}
 \le \sqrt{5} \, \Big\| \sum_{k,\alpha}   w_k(\alpha) x_k(\alpha)
 \Big\|_p.$$
Thus, it remains to prove the estimate
 $$\Big\| \sum_{k, \alpha} w_k(\alpha) \big\langle x_k(\alpha) |
 \Big\|_p \le \sqrt{4d+1} \, \Big\| \sum_{k,\alpha} w_k(\alpha)
 x_k(\alpha) \Big\|_p.$$
To that aim, we use again the projection $\Gamma_\Al(p/2,2d)$ and
Remark \ref{Remark-Proyecciones}
\begin{eqnarray*}
\lefteqn{\Big\| \sum_{k,\alpha} w_k(\alpha) \big\langle
x_k(\alpha) | \Big\|_p^2} \\ & = & \Big\| \sum_{i,j,\alpha,\beta}
w_i(\alpha) \Es \big( x_i(\alpha) x_j(\beta)^* \big) w_j(\beta)^*
\Big\|_{p/2} \\ & = & \Big\| \Gamma_\Al(p/2,2d) \Big[ \Big(
\sum_{k,\alpha} w_k(\alpha) x_k(\alpha) \Big) \Big(
\sum_{k,\alpha} w_k(\alpha) x_k(\alpha) \Big)^* \Big]
\Big\|_{p/2}.
\end{eqnarray*}
Therefore, the assertion follows from the estimate $\big\|
\Gamma_\Al(p/2,2d) \big\| \le 4d + 1$. \fin

Our aim now is to iterate Theorem B to obtain a Khintchine type
inequality, stated as Theorem C in the Introduction, which
generalizes the main results of \cite{Bu1,PP,RX}. Before that, we
analyze in more detail the meaning of the brackets $| \ \rangle$
and $\langle \ |$. That is, according to the mapping $u: \Al \to
C_\infty(\Be)$, we can always write
\begin{equation} \label{Eq-Auxxx2}
| a \big\rangle = u(a) \quad \mbox{and} \quad \big\langle a | =
u(a^*)^*.
\end{equation}
This remark allows us to combine and iterate the brackets $| \
\big\rangle$ and $\big\langle \ |$. In particular, our expressions
for the norms $\Sigma_1$ and $\Sigma_2$ in the statement of
Theorem C (\emph{c.f.} the Introduction) are explained by
\eqref{Eq-Auxxx1} and \eqref{Eq-Auxxx2}.

\begin{lemma} \label{Lemma-Side-Side}
Let $2 \le p \le \infty$ and let $x_k(\alpha), z_k(\alpha)$ and
$w_k(\alpha)$ be homogeneous free polynomials of degree $d_1, d_2$
and $d_3$ respectively for all $1 \le k \le n$ and $\alpha$
running over a finite set $\Lambda$. Assume that $\sum_{k,\alpha}
x_k(\alpha) z_k(\alpha) w_k(\alpha) \in L_p(\Al)$. Then, if
$\RR_k(x_k(\alpha)) = x_k(\alpha)$ and $\LL_k(z_k(\alpha)) = 0$
for all $(k,\alpha)$, we have
 $$\Big\|\sum_{k,\alpha} \big| \, | x_k(\alpha)
 \big\rangle z_k(\alpha)\Big\rangle w_k(\alpha)
 \Big\|_{C_p(C_p(L_p(\Al)))} = \Big\|\sum_{k,\alpha} |
 x_k(\alpha) \, z_k(\alpha) \big\rangle w_k(\alpha)
 \Big\|_{C_p(L_p(\Al))}.$$
Similarly, we have:
\begin{itemize}
\item if $\LL_k(x_k(\alpha)) = x_k(\alpha)$ and
$\RR_k(z_k(\alpha)) = 0$,
 $$\Big\| \sum_{k,\alpha}w_k(\alpha) \Big\langle z_k(\alpha)
 \big\langle x_k(\alpha) | \big| \, \Big\|_{p}
 = \Big\| \sum_{k,\alpha} w_k(\alpha)\big\langle z_k(\alpha) \,
 x_k(\alpha) | \Big\|_{p};$$

\item if $\RR_k(x_k(\alpha)) = 0$ and $\LL_k(z_k(\alpha)) =
z_k(\alpha)$,
 $$\Big\| \sum_{k,\alpha} \big| \, | x_k(\alpha)
 \big\rangle z_k(\alpha) \Big\rangle w_k(\alpha)\Big\|_{p}
 = \Big\| \sum_{k,\alpha} | x_k(\alpha) \, z_k(\alpha)
 \big\rangle w_k(\alpha) \Big\|_{p};$$

\item if $\LL_k(x_k(\alpha)) = 0$ and $\RR_k(z_k(\alpha)) =
z_k(\alpha)$,
 $$\Big\| \sum_{k,\alpha} w_k(\alpha)
 \Big\langle z_k(\alpha) \big\langle x_k(\alpha) | \big| \,
 \Big\|_{p} = \Big\| \sum_{k,\alpha} w_k(\alpha) \big\langle
 z_k(\alpha) \, x_k(\alpha) | \Big\|_{p}.$$
\end{itemize}
\end{lemma}

\dem By freeness we have
\begin{eqnarray*}
\lefteqn{\Big\| \sum_{k,\alpha} | x_k(\alpha) \, z_k(\alpha)
\big\rangle w_k(\alpha) \Big\|_p} \\ & = & \Big\| \Big(
\sum_{i,j,\alpha,\beta} w_i(\alpha)^* \Es \big( z_i(\alpha)^*
x_i(\alpha)^* x_j(\beta) z_j(\beta) \big) w_j(\beta)
\Big)^{\frac12} \Big\|_{p} \\ & = & \Big\| \Big(
\sum_{i,j,\alpha,\beta} w_i(\alpha)^* \Es \big( z_i(\alpha)^*
\Es(x_i(\alpha)^* x_j(\beta)) z_j(\beta) \big) w_j(\beta)
\Big)^{\frac12} \Big\|_{p}.
\end{eqnarray*}
Thus, using the defining property of $u: \Al \to C_\infty(\Be)$,
we obtain
\begin{eqnarray*}
\Big\| \sum_{k,\alpha} | x_k(\alpha) \, z_k(\alpha) \big\rangle
w_k(\alpha) \Big\|_p & = & \Big\| \sum_{k,\alpha} u \Big(
u(x_k(\alpha)) z_k(\alpha) \Big) w_k(\alpha) \Big\|_{p} \\ & = &
\Big\| \sum_{k,\alpha} \big| \, | x_k(\alpha) \big\rangle
z_k(\alpha) \Big\rangle w_k(\alpha) \Big\|_{p}.
\end{eqnarray*}
The three remaining identities follow similarly. This completes
the proof. \fin

In the proof of Theorem C below, we shall use a shorter notation
to write sums like those appearing in the term $\Sigma_2$ (see the
statement of Theorem C in the Introduction) as follows. For a
fixed value $k$ of $j_s$ in $\{1,2, \ldots,n\}$ we shall write
 $$\sum_{\begin{subarray}{c} 1 \le j_1 \neq \cdots \neq j_{s-1} \le
 n \\ 1 \le j_{s+1} \neq \cdots \neq j_{d} \le n \\ j_{s-1} \neq
 j_s = k \neq j_{s+1} \end{subarray}} \quad \mbox{as} \quad
 \sum_{\begin{subarray}{c} j_1 \neq \cdots \neq j_d \\ [j_s = k]
 \end{subarray}}.$$

\vskip5pt

\demC The case of degree $1$ follows automatically from Theorem A.
Now we proceed by induction on $d$. Assume the assertion is true
for degree $d-1$ with relevant constant $\mathcal{C}_p(d-1)$. Then
we apply Theorem B and obtain
\begin{eqnarray*}
\|x\|_p & \sim_{c d^2} & \Big\| \sum_{\alpha \in \Lambda}
\sum_{j_1 \neq \cdots \neq j_d} x_{j_1}(\alpha)
\big\langle x_{j_2}(\alpha) \cdots x_{j_d}(\alpha) | \Big\|_p \\
& + & \Big\| \sum_{\alpha \in \Lambda} \sum_{j_1 \neq \cdots \neq
j_d} | x_{j_1}(\alpha) \big\rangle x_{j_2}(\alpha) \cdots
x_{j_d}(\alpha) \Big\|_p = \mathrm{A} + \mathrm{B}.
\end{eqnarray*}
The resulting terms are homogeneous polynomials of degree $1$ and
$d-1$ respectively. The first one belongs to $R_p(L_p(\Al))$ while
the second one lives in $C_p(L_p(\Al))$. We estimate the first
term by applying Theorem A one more time on the amplified space
$S_p(L_p(\Al))$
\begin{eqnarray*}
\mathrm{A} & \sim_{c} & \Big\| \sum_{\alpha \in \Lambda} \sum_{j_1
\neq \cdots \neq j_d} \Big\langle x_{j_1}(\alpha) \big\langle
x_{j_2}(\alpha) \cdots x_{j_d}(\alpha) | \big| \Big\|_p \\ & + &
\Big\| \sum_{\alpha \in \Lambda} \sum_{j_1 \neq \cdots \neq j_d}
\big| x_{j_1}(\alpha) \big\langle x_{j_2}(\alpha) \cdots
x_{j_d}(\alpha) | \Big\rangle \Big\|_p \\ & + & \Big( \sum_{k=1}^n
\Big\| \sum_{\alpha \in \Lambda} \sum_{\begin{subarray}{c} j_1
\neq \cdots \neq j_d \\ [j_1=k] \end{subarray}} x_{k}(\alpha)
\big\langle x_{j_2}(\alpha) \cdots x_{j_d}(\alpha) | \Big\|_p^p
\Big)^{\frac1p}.
\end{eqnarray*}
According to Lemma \ref{Lemma-Side-Side} and the fact that $u: \Al
\to C_\infty(\Be)$ is a right $\Be$-module map, we easily obtain
\begin{eqnarray}
\label{Eq-A} \mathrm{A} & \sim_{c} & \Big\| \sum_{\alpha \in
\Lambda} \sum_{j_1 \neq \cdots \neq j_d} \big\langle
x_{j_1}(\alpha) x_{j_2}(\alpha) \cdots x_{j_d}(\alpha) | \Big\|_p
\\ \nonumber & + & \Big\| \sum_{\alpha \in \Lambda} \sum_{j_1 \neq \cdots
\neq j_d} | x_{j_1}(\alpha) \big\rangle \, \big\langle
x_{j_2}(\alpha) \cdots x_{j_d}(\alpha) | \Big\|_p
\\ \nonumber & + & \Big( \sum_{k=1}^n \Big\| \sum_{\alpha \in \Lambda}
\sum_{\begin{subarray}{c} j_1 \neq \cdots \neq j_d \\ [j_1=k]
\end{subarray}} x_{k}(\alpha)
\big\langle x_{j_2}(\alpha) \cdots x_{j_d}(\alpha) | \Big\|_p^p
\Big)^{\frac1p}.
\end{eqnarray}
On the other hand, the induction hypothesis gives $\mathrm{B}
\sim_{\mathcal{C}_p(d-1)} \mathrm{B}_1 + \mathrm{B}_2$ with
 $$\mathrm{B}_1 = \sum_{s=1}^{d} \Big\| \sum_{\alpha,
 j_1 \neq \cdots \neq j_d} \big| | x_{j_1}(\alpha) \big\rangle
 \cdots x_{j_s}(\alpha) \Big\rangle \big\langle x_{j_{s+1}}(\alpha)
 \cdots x_{j_d}(\alpha) | \Big\|_p,$$
and $\mathrm{B}_2$ given by
 $$\sum_{s=2}^{d} \Big( \sum_{k=1}^n \Big\| \sum_{\alpha \in \Lambda}
 \sum_{\begin{subarray}{c} j_1 \neq \cdots \neq j_d
 \\ [j_s=k] \end{subarray}} \big| | x_{j_1}(\alpha) \big\rangle
 \cdots x_{j_{s-1}}(\alpha) \Big\rangle x_{j_s}(\alpha) \big\langle
 x_{j_{s+1}}(\alpha) \cdots x_{j_d}(\alpha) | \Big\|_p^p
 \Big)^{\frac1p}_.$$
Moreover, the expressions above are simplified by means of Lemma
\ref{Lemma-Side-Side} as follows
\begin{eqnarray*}
\mathrm{B}_1 & = & \sum_{s=1}^{d} \Big\| \sum_{\alpha, j_1 \neq
\cdots \neq j_d} | x_{j_1}(\alpha) \cdots x_{j_s}(\alpha)
\big\rangle \big\langle x_{j_{s+1}}(\alpha) \cdots x_{j_d}(\alpha)
| \Big\|_p, \\ \mathrm{B}_2 & = & \sum_{s=2}^{d} \Big(
\sum_{k=1}^n \Big\| \sum_{\alpha \in \Lambda}
\sum_{\begin{subarray}{c} j_1 \neq \cdots \neq j_d \\ [j_s=k]
\end{subarray}} | x_{j_1}(\alpha) \cdots \big\rangle
x_{j_s}(\alpha) \big\langle \cdots x_{j_d}(\alpha) | \Big\|_p^p
\Big)^{\frac1p}_.
\end{eqnarray*}
Then we note that the first and third terms in (\ref{Eq-A}) are
the ones which are missing in $\mathrm{B}_1$ and $\mathrm{B}_2$
respectively to obtain $\Sigma_1 + \Sigma_2$, while the middle
term in (\ref{Eq-A}) already appears in $\mathrm{B}_1$. Thus, we
conclude that
 $$\|x\|_p \sim_{\mathcal{C}_p(d)} \Sigma_1 +\Sigma_2$$
where, after keeping track of the constants, we see that
$\mathcal{C}_p(d)$ is controlled by
 $$\mathcal{C}_p(d) \le c d^2 \mathcal{C}_p(d-1).$$
Therefore, the bound $\mathcal{C}_p(d) \le c^d d!^2$ follows from
the recurrence above. \fin

\begin{remark} \label{Remark-Duality-AC}
\emph{From a more functional analytic point of view, the right
hand side of Theorem C can be regarded as the norm of $x$ in an
operator space which is the result of intersecting $2d+1$ operator
spaces, see \cite{Bu1,PP,RX} for more explicit descriptions of
these constructions. We do not state this result in detail since
the notation becomes considerably more complicated. However,
equipped with the description given in \cite{RX} and with Theorem
C, it is not difficult to rephrase Theorem C as a complete
isomorphism between $\mathbf{P}_\Al(p,d)$ and certain $p$-direct
sum of Haagerup tensor products of (subspaces of) $L_p$-spaces.
Moreover, arguing as in \cite{PP} we could extend Theorem C to $1
\le p \le 2$ just replacing intersections by sums of operator
spaces. The same observation is valid for Theorem A.}
\end{remark}

\begin{remark} \label{Remark-Better-Constant}
\emph{The constant $c^d d!^2$ is far from being optimal.
Nevertheless, we can improve the constant in the lower estimate of
Theorem C. To that aim we use the projection
$\Gamma_{\Al}(p/2,2s): L_{p/2}(\Al) \to \mathbf{Q}_\Al(p/2,2s)$ so
that $\Gamma_{\Al}(p/2,2s)(xx^*)$ has the form
 $$\sum_{i_k,j_k,\alpha, \beta} x_{i_1}(\al) \cdots x_{i_s}(\al)
 \Es \big( \cdots x_{i_d}(\al) x_{j_d}(\beta)^* \cdots \big)
 x_{j_{s}}(\beta)^* \cdots x_{j_1}(\beta)^*.$$
A similar expression holds for
$\Gamma_\mathcal{A}(p/2,2(d-s))(x^*x)$
 $$\sum_{i_k,j_k,\alpha, \beta} x_{i_d}(\al)^* \cdots
 x_{i_{s+1}}(\al)^* \Es \big( \cdots x_{i_1}(\al)^* x_{j_1}(\beta)
 \cdots \big) x_{j_{s+1}}(\beta) \cdots x_{j_d}(\beta).$$
Therefore, since $\Gamma_{\Al}(p/2,2d)$ is bounded with constant
$4d+1$, we find
\begin{eqnarray*}
\Big\| \sum_{\alpha \in \Lambda} \sum_{j_1 \neq \cdots \neq j_d}
x_{j_1}(\al) \cdots x_{j_s}(\al) \big\langle x_{j_{s+1}}(\al)
\cdots x_{j_d}(\al) | \Big\|_p & \le & \sqrt{4s+1} \, \|x\|_p,
\\ \Big\| \sum_{\alpha \in \Lambda} \sum_{j_1
\neq \cdots \neq j_d} | x_{j_1}(\al) \cdots x_{j_s}(\al)
\big\rangle x_{j_{s+1}}(\al) \cdots x_{j_d}(\al) \Big\|_p & \le &
\sqrt{4(d-s)+1} \, \|x\|_p.
\end{eqnarray*}
In particular, since $\min(s,d-s) \le d/2$, we deduce
 $$\Big\| \sum_{\al \in \Lambda} \sum_{j_1 \neq \cdots \neq j_d}
 |x_{j_1}(\alpha) \cdots x_{j_{s}}(\alpha) \big\rangle \big\langle
 x_{j_{s+1}}(\alpha) \cdots x_{j_d}(\alpha)| \Big\|_p \le
 \sqrt{2d+1} \, \|x\|_p.$$
Therefore, we have proved the estimate
 $$\Sigma_1 \le (d+1) \sqrt{2d+1} \, \|x\|_p.$$
Similarly, using $\Gamma_\Al(p/2,2)$ as in the proof of Theorem B,
we obtain
\begin{eqnarray*}
\lefteqn{\Big\| \sum_{j_s=1}^n \sum_{\al \in \Lambda} \sum_{j_1
\neq \cdots \neq j_d} |x_{j_1}(\al) \cdots x_{j_{s-1}}(\al)
\big\rangle x_{j_s}(\al) \big\langle x_{j_{s+1}}(\al) \cdots
x_{j_d}(\al) | \Big\|_p} \\ & \le & \sqrt{5} \ \Big\|
\sum_{j_s=1}^n \sum_{\al \in \Lambda} \sum_{j_1 \neq \cdots \neq
j_d} x_{j_1}(\al) \cdots x_{j_s}(\al) \big\langle x_{j_{s+1}}(\al)
\cdots x_{j_d}(\al) | \Big\|_p
\\ & \le & \sqrt{10d+5} \, \|x\|_p.
\end{eqnarray*}
Hence, according to \eqref{Eq-QQpp} we deduce
 $$\Sigma_2 \le 12 d^2 \sqrt{10d+5} \, \|x\|_p.$$
Motivated by the results in \cite{RX}, we conjecture that the
growth of the constant in the upper estimate of Theorem C should
also be polynomial on $d$. However, at the time of this writing we
cannot prove this.}
\end{remark}

\begin{remark} \label{Remarkita}
\emph{Theorem C also generalizes the main results in
\cite{Bu1,PP}. Indeed, note that Theorem C uses $2d+1$ terms in
contrast with the $d+1$ terms in \cite{PP}. However, in the
particular case of free generators it is easily seen that the
terms associated to $\Sigma_1$ (exactly the $d+1$ terms appearing
in \cite{PP}) dominate the terms in $\Sigma_2$. We refer the
reader to the proofs of Lemma \ref{Lemma-Key-Triangle} and Theorem
F below for computations very similar to the ones we are omitting
here. Given $2 \le p \le \infty$ and as a consequence of Theorem C
and Remark \ref{Remark-Better-Constant}, we can rephrase the
Khintchine inequality in \cite{PP} as the following equivalence
for any operator valued $d$-homogeneous polynomial $x$ on the free
generators $\lambda(g_1), \lambda(g_2), \ldots, \lambda(g_n)$
 $$c d^{-3/2} \Sigma_1 \le \|x\|_p \le c^d d!^2 \Sigma_1.$$}
\end{remark}

\section{Square functions}
\label{Section4}

Now we apply our length-reduction formula to study the behavior of
the square function associated to free martingales. More
precisely, according to the Khintchine and Rosenthal inequalities
for free random variables, it is natural to ask whether or not the
noncommutative Burkholder-Gundy inequality \cite{PX1} holds in the
free setting for $p=\infty$, see also \cite{JX} for the
Burkholder-Gundy inequality over non semifinite von Neumann
algebras and \cite{PR,R3} for the weak type $(1,1)$ inequality
associated to it. In this section we find a counterexample to this
question. The following is the key step.

\begin{lemma} \label{Lemma-Key-Triangle}
Let $\mathsf{A}_k = L_\infty(-2,2)$ for $k = 0,1,2, \ldots$
equipped with the Wigner measure, and let $\Al = \mathsf{A}_0 *
\mathsf{A}_1 * \mathsf{A}_2 \cdots$ be the associated reduced free
product equipped with the n.f. tracial state $\phi$. Consider a
free family of semicircular elements $w_k \in \mathsf{A}_{2k-1}$
and $w_k' \in \mathsf{A}_{2k}$ for $k \ge 1$. Given an integer
$n$, fix a mean-zero element $f$ in $\mathsf{A}_0$ such that
 $$\|f\|_{L_2(\mathsf{A}_0)} = 1/\sqrt{n} \quad \mbox{and}
 \quad \|f\|_{L_\infty(\mathsf{A}_0)} = 1.$$
Let $a_{ij}\in\mathcal{B}(\ell_2)$ and
 $$x_{2n} = \sum_{1 \le i,j\le n} a_{ij} \ten w_i f \, w_j'
 \in \mathcal{B}(\ell_2) \ten\Al.$$
Then
 $$\|x_{2n}\|_{\mathcal{B}(\ell_2) \bar\ten \mathcal{A}} \sim_c
 \Big\| \sum_{1 \le i,j\le n} a_{ij} \ten e_{ij}
 \Big\|_{\mathcal{B}(\ell_2) \ten_{\min} \mathcal{B}(\ell_2)}.$$
 \end{lemma}

\dem By Remark \ref{Remark-Operator-Valued}, we have
 $$\mathcal{A} \, \bar\ten \, \mathcal{B}(\ell_2) = \big(
 \mathsf{A}_0 \bar\ten \mathcal{B}(\ell_2) \big)
 *_{\mathcal{B}(\ell_2)} \big( \mathsf{A}_1 \bar\ten
 \mathcal{B}(\ell_2) \big) *_{\mathcal{B}(\ell_2)} \big(
 \mathsf{A}_2 \bar\ten \mathcal{B}(\ell_2) \big)
 *_{\mathcal{B}(\ell_2)} \cdots$$
According to this isometry, we rewrite $x_{2n}$ as follows
\begin{eqnarray*}
x_{2n} = \summ_{i,j} a_{ij} \ten w_i f \, w_j' = \summ_{i,j}
(a_{ij} \ten w_i) (1 \ten f) (1 \ten w_j') = \summ_{i,j} x_{ij} y
z_j.
\end{eqnarray*}
In particular, Theorem C gives the following equivalence for
$\mathsf{E} = \phi \ten id_{\mathcal{B}(\ell_2)}$
\begin{eqnarray*}
\lefteqn{\|x_{2n}\|_{\mathcal{B}(\ell_2) \bar\ten \mathcal{A}}} \\
& \sim_{c} & \big\| \Es(x_{2n} x_{2n}^*) \big\|_\infty^{\frac12} +
\big\| \Es(x_{2n}^* x_{2n}) \big\|_\infty^{\frac12} \\
[5pt] & + & \Big\| \sum_{i,j=1}^n x_{ij} \big\langle y z_j |
\Big\|_{\infty} + \Big\| \sum_{i,j=1}^n | x_{ij} y \big\rangle z_j
\Big\|_{\infty} \\ & + & \Big\| \sum_{i,j=1}^n | x_{ij} y
\big\rangle \big\langle z_j | \Big\|_{\infty} + \Big\|
\sum_{i,j=1}^n | x_{ij} \big\rangle \big\langle y z_j |
\Big\|_{\infty} + \Big\| \sum_{i,j=1}^n | x_{ij} \big\rangle \, y
\, \big\langle z_j | \Big\|_{\infty}
\\ \\ & = & \mathrm{A} + \mathrm{B} + \mathrm{C} + \mathrm{D} +
\mathrm{E} + \mathrm{F} + \mathrm{G}.
\end{eqnarray*}
It is clear that
 $$\mathrm{A} = \Big\| \sum_{ijkl} a_{ij} a_{kl}^*
 \phi \big( w_i f w_{j}' w_{l}' f^* w_k \big)
 \Big\|_\infty^{\frac12} = \Big\| \sum_{ij} a_{ij} a_{ij}^* \phi
 \big( w_i f w_{j}^{_{'}2} f^* w_i \big) \Big\|_\infty^{\frac12}.$$
Since $\phi(w_k^2) = \phi(w_k^{_{'}2}) =1$, this gives
 $$\mathrm{A} = \|f\|_2 \Big\| \sum_{i,j=1}^n a_{ij} a_{ij}^*
 \Big\|_\infty^{\frac12} = \frac{1}{\sqrt{n}} \, \Big\|
 \sum_{i,j=1}^n a_{ij} \ten e_{1,ij} \Big\|_{\mathcal{B}(\ell_2)
 \ten_{\min} \mathcal{B}(\ell_2)}.$$
The same argument gives rise to the identity
 $$\mathrm{B} =
 \frac{1}{\sqrt{n}} \, \Big\| \sum_{i,j=1}^n a_{ij} \ten e_{ij,1}
 \Big\|_{\mathcal{B}(\ell_2) \ten_{\min} \mathcal{B}(\ell_2)}.$$
Let us estimate the term $\mathrm{C}$
\begin{eqnarray*}
\mathrm{C} & = & \Big\| \sum_{ijkl} x_{ij} \Es(y z_j z_l^* y^*)
x_{kl}^* \Big\|_\infty^{\frac12} \\ & = & \Big\| \sum_{ijkl}
a_{ij} a_{kl}^* \ten w_i \phi \big( f w_j' w_l' f^* \big) w_k
\Big\|_\infty^{\frac12} \\ & = & \|f\|_2 \Big\| \sum_{j=1}^n
e_{1j} \ten \Big(\sum_{i=1}^n a_{ij} \ten w_i \Big) \Big\|_\infty
\\ & = & \|f\|_2 \Big\| \sum_{i=1}^n \Big(\sum_{j=1}^n a_{ij} \ten
e_{1j} \Big) \ten w_i \Big\|_\infty.
\end{eqnarray*}
Now, applying the Khintchine inequality for free random variables
\cite{HP}
\begin{eqnarray*}
\mathrm{C} \!\! & \sim & \!\! \|f\|_2 \max \left\{ \Big\|
\sum_{i=1}^n \Big(\sum_{j=1}^n a_{ij} \ten e_{1j} \Big) \ten
e_{1i} \Big\|_\infty, \Big\| \sum_{i=1}^n \Big(\sum_{j=1}^n
a_{ij} \ten e_{1j} \Big) \ten e_{i1} \Big\|_\infty \right\} \\
\!\! & = & \!\! \|f\|_2 \max \left\{ \Big\| \sum_{i,j=1}^n a_{ij}
\ten e_{1,ij} \Big\|_{\mathcal{B}(\ell_2) \ten_{\min}
\mathcal{B}(\ell_2)}, \Big\| \sum_{i,j=1}^n a_{ij} \ten e_{ij}
\Big\|_{\mathcal{B}(\ell_2) \ten_{\min} \mathcal{B}(\ell_2)}
\right\}.
\end{eqnarray*}
Again the same argument gives
 $$\mathrm{D} \, \sim \ \|f\|_2 \max \left\{ \Big\|
 \sum_{i,j=1}^n a_{ij} \ten e_{ij,1} \Big\|_{\mathcal{B}(\ell_2)
 \ten_{\min} \mathcal{B}(\ell_2)}, \Big\| \sum_{i,j=1}^n a_{ij}
 \ten e_{ij} \Big\|_{\mathcal{B}(\ell_2) \ten_{\min}
 \mathcal{B}(\ell_2)} \right\}.$$ The term $\mathrm{E}$ is
calculated as follows
\begin{eqnarray*}
\mathrm{E} & = & \Big\| \sum_{ijkl} u(x_{ij} y) \Es(z_j z_l^*)
u(x_{kl} y)^* \Big\|_\infty^{\frac12}
\\ & = & \Big\| \sum_{j=1}^n \Big( \sum_{i=1}^n u \big( a_{ij}
\ten w_i f \big) \Big) \Big( \sum_{k=1}^n u \big( a_{kj} \ten w_k
f \big) \Big)^* \Big\|_\infty^{\frac12}
\\ & = & \Big\| \sum_{j=1}^n e_{1j} \ten \Big( \sum_{i=1}^n
u \big( a_{ij} \ten w_i f \big) \Big) \Big\|_\infty
\\ & = & \Big\| \sum_{i,j=1}^n e_{ij} \ten \Big(\sum_{r,s=1}^n
a_{ri}^* \phi \big( f^* w_r w_s f \big) a_{sj}  \Big)
\Big\|_\infty^{\frac12} \\ & = & \|f\|_2 \Big\| \sum_{i,j=1}^n
a_{ij} \ten e_{ij} \Big\|_{\mathcal{B}(\ell_2) \ten_{\min}
\mathcal{B}(\ell_2)} = \frac{1}{\sqrt{n}} \, \Big\| \sum_{i,j=1}^n
a_{ij} \ten e_{ij} \Big\|_{\mathcal{B}(\ell_2) \ten_{\min}
\mathcal{B}(\ell_2)}.
\end{eqnarray*}
The same identity holds for $\mathrm{F}$
 $$\mathrm{F} =\frac{1}{\sqrt{n}} \,
  \Big\| \sum_{i,j=1}^n a_{ij} \ten e_{ij}
 \Big\|_{\mathcal{B}(\ell_2) \ten_{\min} \mathcal{B}(\ell_2)}.$$
The calculation of $\mathrm{G}$ is very similar
\begin{eqnarray*}
\mathrm{G} & = & \Big\| \sum_{ijkl} u (x_{ij}) y \Es(z_j z_l^*) y^*
u (x_{kl})^* \Big\|_\infty^{\frac12} \\
& = & \Big\| \sum_{j=1}^n \Big( \sum_{i=1}^n u \big( a_{ij} \ten
w_i \big) y \Big) \Big( \sum_{k=1}^n u \big( a_{kj} \ten w_k \big)
y \Big)^* \Big\|_\infty^{\frac12} \\ & = & \Big\| \sum_{j=1}^n
e_{1j} \ten
\Big( \sum_{i=1}^n u \big( a_{ij} \ten w_i \big) y \Big) \Big\|_\infty \\
& = & \Big\| \sum_{i,j=1}^n e_{ij} \ten \Big( \sum_{r,s=1}^n
a_{ri}^* a_{sj} \ten f^* \phi (w_r w_s) f \Big)
\Big\|_\infty^{\frac12} \\ & = & \|f\|_\infty \Big\|
\sum_{i,j=1}^n a_{ij} \ten e_{ij} \Big\|_{\mathcal{B}(\ell_2)
\ten_{\min} \mathcal{B}(\ell_2)} = \Big\| \sum_{i,j=1}^n a_{ij}
\ten e_{ij} \Big\|_{\mathcal{B}(\ell_2) \ten_{\min}
\mathcal{B}(\ell_2)}.
\end{eqnarray*}
On the other hand, we observe that the maps on
$\mathcal{B}(\ell_2) \ten_{\min} \mathcal{B}(\ell_2)\ten_{\min}
\mathcal{B}(\ell_2)$
\begin{eqnarray*}
\sum_{i,j=1}^n a_{ij} \ten e_{i1} \ten e_{1j} & \mapsto &
\sum_{i,j=1}^n a_{ij} \ten e_{1i} \ten e_{1j}, \\
\sum_{i,j=1}^n a_{ij} \ten e_{i1} \ten e_{1j} & \mapsto &
\sum_{i,j=1}^n a_{ij} \ten e_{i1} \ten e_{j1},
\end{eqnarray*}
have norm $\sqrt{n}$. Indeed, this follows automatically from the
well-known fact that the natural mappings $R_n \to C_n$ and $C_n
\to R_n$ between $n$-dimensional row and column Hilbert spaces are
completely bounded with cb-norm $\sqrt{n}$, see e.g. \cite{ER} or
\cite{P3} for the proof. Thus we deduce
\begin{eqnarray*}
\Big\|\sum_{i,j=1}^n a_{ij} \ten e_{1,ij}
\Big\|_{\mathcal{B}(\ell_2) \ten_{\min} \mathcal{B}(\ell_2)} & \le
& \sqrt{n} \, \Big\|\sum_{i,j=1}^n a_{ij} \ten e_{ij}
\Big\|_{\mathcal{B}(\ell_2) \ten_{\min} \mathcal{B}(\ell_2)}, \\
\Big\|\sum_{i,j=1}^n a_{ij} \ten e_{ij,1}
\Big\|_{\mathcal{B}(\ell_2) \ten_{\min} \mathcal{B}(\ell_2)} & \le
& \sqrt{n} \, \Big\|\sum_{i,j=1}^n a_{ij} \ten e_{ij}
\Big\|_{\mathcal{B}(\ell_2) \ten_{\min} \mathcal{B}(\ell_2)}.
\end{eqnarray*}
The assertion follows easily from these inequalities and the
estimates above. \fin

The idea to find our counterexample follows an argument from
\cite{PX1}. We consider a suitable martingale for which the
Burkholder-Gundy inequality implies an upper estimate for the
triangular projection on $\mathcal{B}(\ell_2^n)$. This gives the
logarithmic growth stated in Theorem D. After the proof of our
counterexample or Theorem D, we shall study the reverse estimate
for free martingales whose martingale differences are polynomials
of a bounded degree.

\vskip5pt

\demD Let us define
 $$x_{2n} = \sum_{1 \le i,j \le n} a_{ij} w_i f
 w_j' \quad \mbox{with} \quad a_{ij} \in \C.$$
Here $w_i, f$ and $w_j'$ are defined as in Lemma
\ref{Lemma-Key-Triangle}. Moreover, the enumeration given in the
statement of Lemma \ref{Lemma-Key-Triangle} for the algebras
$\mathsf{A}_0, \mathsf{A}_1, \mathsf{A}_2, \ldots$ provides a
natural martingale structure for the $x_{2n}$'s, i.e. with respect
to the natural filtration $(\Al_k)$ defined by $\Al_k =
\mathsf{A}_0
* \mathsf{A}_1 * \mathsf{A}_2 * \cdots * \mathsf{A}_k$.
An easy inspection gives the following expressions valid for all
$k \ge 0$
\begin{equation} \label{Equation-Odd-Even}
dx_{2k} = \sum_{1 \le i \le k} a_{ik} w_i f w_k' \quad \mbox{and}
\quad d x_{2k-1} = \sum_{1 \le j < k} a_{kj} w_k f w_j'\,.
\end{equation}
We are interested in the best constant $\mathcal{K}_n$ for
 $$\max \left\{ \Big\| \Big( \sum_{k=1}^{2n} dx_k dx_k^*
 \Big)^{\frac12} \Big\|_\infty, \Big\| \Big( \sum_{k=1}^{2n} dx_k^*
 dx_k \Big)^{\frac12} \Big\|_\infty \right\} \le \mathcal{K}_n
 \Big\| \sum_{k=1}^{2n} dx_k \Big\|_\infty.$$
According to Lemma \ref{Lemma-Key-Triangle}, we have
 $$\Big\|\sum_{k=1}^{2n} dx_k \Big\|_{\infty}\sim_c
 \Big\| \sum_{i,j=1}^na_{ij} e_{ij} \Big\|_{\mathcal{B}(\ell_2)}.$$
On the other hand, we observe that
 $$\begin{array}{rcccl} \displaystyle \Big\| \Big( \sum_{k=1}^{n}
 dx_{2k} dx_{2k}^* \Big)^{\frac12} \Big\|_\infty \!\!\!\! & = &
 \!\!\!\! \displaystyle \Big\| \sum_{k=1}^{n} e_{1k} \ten dx_{2k}
 \Big\|_\infty \!\!\!\! & = & \!\!\!\! \displaystyle \Big\| \sum_{i
 \le k} a_{ik} \, e_{1k} \ten w_i f w_k' \Big\|_\infty, \\
 \displaystyle \Big\| \Big( \sum_{k=1}^{n} dx_{2k-1}^* dx_{2k-1}
 \Big)^{\frac12} \Big\|_\infty \!\!\!\! & = & \!\!\!\!
 \displaystyle \Big\| \sum_{k=1}^{n} e_{k1} \ten dx_{2k-1}
 \Big\|_\infty \!\!\!\! & = & \!\!\!\! \displaystyle \Big\| \sum_{k
 > j} a_{kj} e_{k1} \ten w_k f w_j' \Big\|_\infty.
 \end{array}$$
Thus, we may apply Lemma \ref{Lemma-Key-Triangle} one more time
and obtain
 $$\begin{array}{rcccl} \displaystyle \Big\| \Big(
 \sum_{k=1}^{n} dx_{2k} dx_{2k}^* \Big)^{\frac12} \Big\|_\infty
 \!\!\!\! & \sim_c & \!\!\!\! \displaystyle \Big\| \sum_{i \le j}
 a_{ij} e_{1j} \ten e_{ij} \Big\|_{\mathcal{B}(\ell_2 \ten \ell_2)}
 \!\!\!\! & = & \!\!\!\! \displaystyle \Big\|
 \sum_{\begin{subarray}{c} i,j =1 \\ i \le j
 \end{subarray}}^n a_{ij} e_{ij} \Big\|_{\mathcal{B}(\ell_2)}, \\
 \displaystyle \Big\| \Big( \sum_{k=1}^{n} dx_{2k-1}^* dx_{2k-1}
 \Big)^{\frac12} \Big\|_\infty \!\!\!\! & \sim_c & \!\!\!\!
 \displaystyle \Big\| \sum_{i > j} a_{ij} e_{i1} \ten e_{ij}
 \Big\|_{\mathcal{B}(\ell_2 \ten \ell_2)} \!\!\!\! & = & \!\!\!\!
 \displaystyle \Big\| \sum_{\begin{subarray}{c} i,j=1 \\ i > j
 \end{subarray}}^n a_{ij} e_{ij} \Big\|_{\mathcal{B}(\ell_2)}.
 \end{array}$$
That is, $\mathcal{K}_n$ is bounded from below by $c$ times the
norm of the triangular projection on $\Be(\ell_2^n)$. However, it
is well-known that the norm of the triangular projection grows
like $\log n$, see e.g. Kwapie\'n/Pelczynski \cite{KP}. This
completes the proof. \fin

After Theorem D, it remains open to see whether or not the reverse
estimate in the Burkholder-Gundy inequalities holds for free
martingales in $L_\infty(\mathcal{A})$. In the following result we
give a partial solution to this problem. We will work with free
martingales of the form $x_n = \sum_{k=1}^n d x_k$ with
\begin{equation} \label{Eq-diffs-dd}
d x_k = \sum_{\alpha \in \Lambda}^{\null} \sum_{j_1 \neq \cdots
\neq j_d} a_{j_1}^k(\alpha) \cdots a_{j_d}^k(\alpha) \quad
\mbox{and} \quad a_{j_s}^k(\alpha) \in \bubl_{j_s},
\end{equation}
where $1 \le j_1, j_2, \ldots, j_d \le k$. That is, we assume that
all the martingale differences are $d$-homogeneous free
polynomials. We shall refer to this kind of martingales as
\emph{$d$-homogeneous free martingales}. More generally, if $x$ is
a free martingale with $dx_k$ being a (not necessarily
homogeneous) free polynomial of degree $d$, we shall simply say
that $x$ is a \emph{$d$-polynomial free martingale}. We shall also
use the following notation
$$\mathcal{S}_\infty(x,n) = \max \left\{ \Big\| \Big(
\sum_{k=1}^{n} d x_k d x_k^* \Big)^{\frac12} \Big\|_\infty, \,
\Big\| \Big( \sum_{k=1}^{n} d x_k^* d x_k \Big)^{\frac12}
\Big\|_\infty \right\}.$$ Our main tools in the following result
are again Theorems A and B.

\begin{proposition} \label{Proposition-Sq1}
If $x$ is a $d$-polynomial free martingale, $$\Big\| \sum_{k=1}^n
d x_k \Big\|_\infty \le c^d d^2 \sqrt{d!} \,
\mathcal{S}_\infty(x,n).$$
\end{proposition}

\dem Let us consider the inequality
\begin{equation} \label{Eq-parahom}
\Big\| \sum_{k=1}^n d x_k \Big\|_\infty \le \mathcal{C}(d) \,
\mathcal{S}_\infty(x,n)
\end{equation}
valid for any $d$-homogeneous free martingale $x$ with $d \ge 0$.
To prove \eqref{Eq-parahom} and estimate $\mathcal{C}(d)$ we
proceed by induction on $d$. Namely, for $d=0$ we have $dx_1 =
\Es(x)$ and $dx_k = 0$ for $k=2,3,\ldots$ In particular,
 $$\Big\| \sum_{k=1}^n dx_k \Big\|_\infty = \|dx_1\|_\infty \le
 \mathcal{S}_\infty(x,n).$$
Therefore, \eqref{Eq-parahom} holds for $d=0$ with $\mathcal{C}(0)
= 1$. If $d=1$ we observe that
 $$\sum_{k=1}^n dx_k = \sum_{k=1}^n \LL_k(dx_k).$$
Thus, Proposition \ref{Proposition-Sum-L&R} gives
 $$\Big\|\sum_{k=1}^n dx_k \Big\|_\infty
 \le 3 \max \left\{ \Big\|\sum_{k=1}^n \LL_k(dx_k)
  \LL_k(dx_k)^* \Big\|_\infty^{\frac12},
 \Big\| \sum_{k=1}^n \LL_k(dx_k)^* \LL_k(dx_k)
 \Big\|_\infty^{\frac12} \right\}.$$
This, combined with the proof of Lemma
\ref{Lemma-Complementation-Boundedness}, gives rise to
 $$\Big\| \sum_{k=1}^n dx_k \Big\|_\infty \le 9 \,
 \mathcal{S}_\infty(x,n).$$
In particular, \eqref{Eq-parahom} holds for $d=1$ with
$\mathcal{C}(1) \le 9$. Now we assume that \eqref{Eq-parahom}
holds for $(d-1)$-homogeneous free martingales with some constant
$\mathcal{C}(d-1)$. To prove \eqref{Eq-parahom} for a
$d$-homogeneous free martingale $x$, we decompose the martingale
differences by means of the mappings $\LL_k$ as follows
 $$\Big\| \sum_{k=1}^n dx_k \Big\|_\infty \le \Big\| \sum_{k=1}^n
 \LL_k(dx_k) \Big\|_\infty + \Big\| \sum_{k=1}^n (id_\mathcal{A} -
 \LL_k)(dx_k) \Big\|_\infty = \mathrm{A} + \mathrm{B}.$$
The estimate
 \begin{equation} \label{Equation-AAA} \mathrm{A}
 \le 9 \, \mathcal{S}_\infty(x,n),
 \end{equation}
follows as the inequality $\mathcal{C}(1) \le 9$ above. On the
other hand, we have
 $$(id_\mathcal{A} - \LL_k)(dx_k) = \sum_{\alpha \in \Lambda}
 \sum_{j=1}^{k-1} x_j^k(\alpha) w_j^k(\alpha)$$
with $x_j^k(\alpha) \in \bubl_j$ and $w_j^k(\alpha) \in
\mathbf{P}_\Al(d-1)$ satisfying $\LL_j(w_j^k(\alpha)) = 0$.
Indeed, this follows from the fact that no word in
$(id_\mathcal{A} - \LL_k)(dx_k)$ starts with a mean-zero letter in
$\mathsf{A}_k$ and that $dx_k\in\mathcal{A}_k$. Thus, we may write
$\mathrm{B}$ in the form
 $$\mathrm{B} = \Big\| \sum_{(\alpha,k) \in \Delta} \,
 \sum_{j=1}^{n-1}x_j(\alpha,k) w_j(\alpha,k) \Big\|_\infty$$
with $\Delta = \Lambda \times \big\{ 1,2, \ldots, n \big\}$ and
 $$x_j(\alpha,k) w_j(\alpha,k) =
 \begin{cases} 0 & \mbox{if} \ j \ge k,
 \\ x_j^k(\alpha) w_j^k(\alpha) & \mbox{if} \ j < k. \end{cases}$$
According to Theorem B we obtain
\begin{eqnarray*}
\Big\| \sum_{(\alpha,k),j} x_j(\alpha,k) w_j(\alpha,k)
\Big\|_\infty & \le & \Big\| \sum_{(\alpha,k),j} x_j(\alpha,k)
\big\langle w_j(\alpha,k) | \Big\|_\infty \\ & + & \Big\|
\sum_{(\alpha,k),j} | x_j(\alpha,k) \big\rangle w_j(\alpha,k)
\Big\|_\infty = \mathrm{B}_1 + \mathrm{B}_2.
\end{eqnarray*}
Note that the constant $1$ in the inequality above holds since we
are only considering the case $(p,q) = (\infty,2)$ in the proof of
Theorem B. Let us start by estimating the first term
$\mathrm{B}_1$. We claim that
 $$\mathrm{B}_1^2 = \Big\|\sum_{k=1}^n
 \sum_{\alpha,\beta \in \Lambda} \sum_{j_1, j_2 =1}^{n-1}
  x_{j_1}(\alpha,k) \Es \big( w_{j_1}(\alpha,k)
  w_{j_2}(\beta,k)^* \big) x_{j_2}(\beta,k)^* \Big\|_{\infty}.$$
To see this, it suffices to show that
 $$\Es \big(w_{j_1}(\alpha,k_1) w_{j_2}(\beta,k_2)^* \big) = 0$$
for $k_1 \neq k_2$. Indeed, let us assume without lost of
generality that $k_1 < k_2$. Then we know by construction that
$w_{j_1}(\alpha,k_1) \in \Al_{k_1}$ and that $w_{j_2}(\beta,k_2)$
contains a mean-zero letter in $\mathsf{A}_{k_2}$ with $k_2
> k_1$. Thus, our claim follows easily by freeness.
Hence, we may write the identity above as follows
 $$\mathrm{B}_1 = \Big\|
 \sum_{k=1}^n e_{1k} \ten \sum_{\alpha \in \Lambda}
 \sum_{j=1}^{k-1} x_j^k(\alpha) \big\langle w_j^k(\alpha) |
 \Big\|_\infty.$$
Arguing as in the proof of Theorem B, we obtain
\begin{eqnarray}
\label{Equation-B1B1B1} \mathrm{B}_1 \!\!\! & \le & \!\!\!
\sqrt{5} \ \Big\| \sum_{k=1}^n e_{1k} \ten \sum_{\alpha \in
\Lambda} \sum_{j=1}^{k-1} x_j^k(\alpha) w_j^k(\alpha) \Big\|_\infty \\
\nonumber \!\!\! & = & \!\!\! \sqrt{5} \ \Big\| \sum_{k=1}^n
e_{1k} \ten (id_\mathcal{A} - \LL_k) (d x_k) \Big\|_\infty \\
\nonumber \!\!\! & \le & \!\!\! \sqrt{5} \, \Big[ \Big\|
\sum_{k=1}^n e_{1k} \ten d x_k \Big\|_\infty + \Big\| \sum_{k=1}^n
e_{1k} \ten \LL_k (d x_k) \Big\|_\infty \Big] \le 4 \sqrt{5} \,
\mathcal{S}_\infty(x,n),
\end{eqnarray}
where the last inequality follows from Lemma
\ref{Lemma-Complementation-Boundedness} one more time.

\vskip5pt

To estimate $\mathrm{B}_2$ we observe that
\begin{equation} \label{Eq-martingaled-1}
\sum_{k=1}^n \sum_{\alpha \in \Lambda} \sum_{j=1}^{k-1} |
x_j^k(\alpha) \big\rangle\, w_j^k(\alpha)
\end{equation}
can be regarded as a sum of martingale differences on the von
Neumann algebra $\Al \bar\ten \mathcal{B}(\ell_2)$ with respect to
the index $k$ and the filtration $\Al_1 \bar\ten \Be(\ell_2),
\Al_2 \bar\ten \Be(\ell_2), \ldots$ Indeed, we have
 $$\Es_{k-1} \ten id_{\Be(\ell_2)} \Big( \sum_{\alpha \in
 \Lambda} \sum_{j=1}^{k-1} | x_j^k(\alpha) \big\rangle
 w_j^k(\alpha) \Big) = \sum_{\alpha \in \Lambda} \sum_{j=1}^{k-1} |
 x_j^k(\alpha) \big\rangle \Es_{k-1}(w_j^k(\alpha)) = 0.$$
Then, since \eqref{Eq-martingaled-1} forms a $(d-1)$-homogeneous
free martingale, we may apply the induction hypothesis and obtain
in this way the following upper bound for $\mathrm{B}_2$
 $$\mathcal{C}(d-1) \max \left\{ \Big\| \sum_{k=1}^n e_{1k} \ten
 \sum_{\alpha,j} | x_j^k(\alpha) \big\rangle w_j^k(\alpha)
 \Big\|_\infty, \Big\| \sum_{k=1}^n e_{k1} \ten \sum_{\alpha,j} |
 x_j^k(\alpha) \big\rangle w_j^k(\alpha) \Big\|_\infty \right\}.$$
Then, arguing as in the proof of Theorem B $(2(2(d-1))+1 = 4d-3)$,
we deduce
\begin{eqnarray*}
\mathrm{B}_2 & \le & \sqrt{4d-3} \, \mathcal{C}(d-1) \\ & \times &
\max \left\{ \Big\| \sum_{k=1}^n e_{1k} \ten (id_\mathcal{A} -
\LL_k)(dx_k) \Big\|_\infty, \Big\| \sum_{k=1}^n e_{k1} \ten
(id_\mathcal{A} - \LL_k)(dx_k) \Big\|_\infty \right\}.
\end{eqnarray*}
The triangle inequality and Lemma
\ref{Lemma-Complementation-Boundedness} produce
\begin{equation} \label{Equation-B2B2B2}
\mathrm{B}_2 \le 4 \sqrt{4d-3} \, \mathcal{C}(d-1) \,
\mathcal{S}_\infty(x,n).
\end{equation}
Now (\ref{Equation-AAA}, \ref{Equation-B1B1B1},
\ref{Equation-B2B2B2}) give
 $$\mathcal{C}(d) \le (9 + 4 \sqrt{5})
 + 4 \sqrt{4d-3} \, \mathcal{C}(d-1) \le c \sqrt{d} \,
 \mathcal{C}(d-1).$$
Iterating the recurrence and using $\mathcal{C}_0 = 1$ we find
$\mathcal{C}(d) \le c^d \sqrt{d!}$. Therefore,
\begin{equation} \label{Eq-parahom2}
\Big\| \sum_{k=1}^n d x_k \Big\|_\infty \le c^d \sqrt{d!}
\mathcal{S}_\infty(x,n)
\end{equation}
for $d$-homogeneous free martingales.

\vskip5pt

Now let $x$ be any $d$-polynomial free martingale $x$. We may
decompose $x$ into its homogeneous parts $dx_k = \sum_s dx_k^s$
with $0 \le s \le d$. It is clear that $dx_1^s, dx_2^s, dx_3^s,
\ldots$ are the martingale differences of an $s$-homogeneous free
martingale $x^s$. Therefore, applying \eqref{Eq-parahom2} we
deduce
 $$\Big\| \sum_{k=1}^n dx_k \Big\|_\infty \le \sum_{s=0}^d
 \Big\| \sum_{k=1}^n dx_k^s \Big\|_\infty \le \|\Es(x)\|_\infty +
 \sum_{s=1}^d c^s \sqrt{s!} \, \mathcal{S}_\infty(x^s,n).$$
For the first term we have
 $$\|\Es(x)\|_\infty = \|\Es(\Es_1(x))\|_\infty \le
 \|\Es_1(x)\|_\infty = \|dx_1\|_\infty \le
 \mathcal{S}_\infty(x,n).$$
The rest of the terms are estimated by Theorem
\ref{Theorem-Complementation}
\begin{eqnarray*}
\mathcal{S}_\infty(x^s,n) & \sim & \Big\| \sum_{k=1}^n e_{1k} \ten
dx_k^s \Big\|_\infty + \Big\| \sum_{k=1}^n e_{k1} \ten dx_k^s
\Big\|_\infty \\ & = & \Big\| \big( id_{\mathcal{B}(\ell_2)} \ten
\Pi_\mathcal{A}(\infty,s) \big) \Big( \sum_{k=1}^n e_{1k} \ten
dx_k \Big) \Big\|_\infty \\ & + & \Big\| \big(
id_{\mathcal{B}(\ell_2)} \ten \Pi_\mathcal{A}(\infty,s) \big)
\Big( \sum_{k=1}^n e_{k1} \ten dx_k \Big) \Big\|_\infty \le  4s \,
\mathcal{S}_\infty(x,n).
\end{eqnarray*}
Our estimates give rise to
 $$\Big\| \sum_{k=1}^n dx_k
 \Big\|_\infty \le \Big( 1 + 4 \sum_{s=1}^d c^s s \sqrt{s!} \Big)
 \mathcal{S}_\infty(x,n) \le c^d d^2 \sqrt{d!} \,
 \mathcal{S}_\infty(x,n).$$
This is the desired estimate. The proof is complete. \fin

\begin{remark}
\emph{Proposition \ref{Proposition-Sq1} extends to the case $2 \le
p \le \infty$. Indeed, we just need to replace Proposition
\ref{Proposition-Sum-L&R} by Corollary \ref{Corollary-Voiculescu1}
and apply Theorem B in full generality. Of course, this would
provide a worse constant. The relevance of Proposition
\ref{Proposition-Sq1} lies however in the fact that the resulting
constants are uniformly bounded as $p \to \infty$, in contrast
with the non-free setting \cite{PX1}.}
\end{remark}

\section{Generalized circular systems}
\label{Section-Generalized-Circular}

In this last section we illustrate our results by investigating
Khintchine type inequalities for Shlyakhtenko's generalized
circular systems and Hiai's generalized $q$-gaussians. Given an
infinite dimensional and separable Hilbert space $\mathcal{H}$
equipped with a distinguished unit vector or vacuum $\Omega$, we
denote by $\mathcal{F}(\mathcal{H})$ the associated Fock space
 $$\mathcal{F}(\mathcal{H}) = \C \Omega \oplus \bigoplus_{n \ge 1}
 \mathcal{H}^{\otimes n}.$$
Given any vector $e \in \mathcal{H}$, we denote by $\ell(e)$ the
left creation operator on $\mathcal{F}(\mathcal{H})$ associated
with $e$, which acts by tensoring from the left. The adjoint map
$\ell^*(e)$ is called the annihilation operator on
$\mathcal{F}(\mathcal{H})$, see \cite{VDN} for more details. Let
us fix an orthonormal basis $(e_{\pm k})_{k \ge 1}$ in
$\mathcal{H}$ and two sequences $(\lambda_k)_{k \ge 1}$ and
$(\mu_k)_{k \ge 1}$ of positive numbers. Set
 $$g_k = \lambda_k \ell(e_k) + \mu_k \ell^*(e_{-k}).$$
The $g_k$'s are \emph{generalized circular random variables}
studied by Shlyakhtenko \cite{S}. Let $\Gamma$ denote the von
Neumann generated by the generalized circular system $(g_k)_{k \ge
1}$.  $\Gamma$ is equipped with the vacuum state $\phi$ given by
$\phi(x) = \langle \Omega,\; x \Omega \rangle$. According to
\cite{S}, $\phi$ is faithful and the $g_k$'s are free with respect
to $\phi$. In fact, if $\Gamma_k$ is the von Neumann subalgebra of
$\Gamma$ generated by $g_k$, then $(\Gamma, \phi) = *_{k \ge 1}
(\Gamma_k, \phi_{\Gamma_k})$. Shlyakhtenko also calculated in
\cite{S} the modular group and showed that $\sigma_t(g_k) =
(\lambda_k^{-1} \mu_k)^{2it} g_k$. In particular, the $g_k$'s are
analytic elements of $\Gamma$ and eigenvectors of the modular
automorphism group $\sigma$. Let us write $d_\phi$ for the density
associated to the state $\phi$ on $\Gamma$. We shall also need the
elements
\begin{equation} \label{Eq-Generalized-p}
g_{k,p} = d_\phi^{\frac{1}{2p}} g_k d_\phi^{\frac{1}{2p}} =
(\lambda_k^{-1} \mu_k)^{\frac1p}\, g_k\, d_\phi^{^\frac{1}{p}} =
(\lambda_k \mu_k^{-1})^{\frac1p}\, d_\phi^{^\frac{1}{p}}\, g_k.
\end{equation}

\vskip5pt

The following is the Khintchine type inequality for
$1$-homogeneous polynomials on generalized circular random
variables. Its proof can be found in \cite{X2}, where the
third-named author used Theorem A to obtain constants independent
of $p$. When $\lambda_k = \mu_k$ for $k \ge 1$, the $g_k$'s become
a usual circular system in Voiculescu's sense and the result below
reduces to Theorem 8.6.5 in \cite{P2}. On the other hand, the case
$p=\infty$ was already proved by Pisier and Shlyakhtenko in
\cite{PS}.

\begin{theorem} \label{Th-OldE}
Let $\mathcal{N}$ be a von Neumann algebra and $1 \le p \le
\infty$. Let us consider a finite sequence $x_1, x_2, \ldots, x_n$
in $L_p(\mathcal{N})$. Then, the following equivalences hold up to
an absolute constant $c$ independent of $n$
\begin{itemize}
\item[i)] If $1 \le p \le 2$, then
\begin{eqnarray*}
\lefteqn{\Big\| \sum_{k=1}^n x_k \ten g_{k,p} \Big\|_p} \\ &
\sim_c & \inf_{x_k = a_k + b_k} \Big\| \Big( \sum_{k=1}^n
\lambda_k^{\frac{2}{p}} \mu_k^{\frac{2}{p'}} a_ka_k^*
\Big)^{\frac12} \Big\|_p + \Big\| \Big( \sum_{k=1}^n
\lambda_k^{\frac{2}{p'}} \mu_k^{\frac{2}{p}} b_k^*b_k
\Big)^{\frac12} \Big\|_p
\end{eqnarray*}

\item[ii)] If $2\le p\le \infty$, then
\begin{eqnarray*}
\lefteqn{\Big\| \sum_{k=1}^n x_k \ten g_{k,p} \Big\|_p} \\ &
\sim_c & \max \left\{ \Big\| \Big( \sum_{k=1}^n
\lambda_k^{\frac{2}{p}} \mu_k^{\frac{2}{p'}} x_k x_k^*
\Big)^{\frac12} \Big\|_p \, , \, \Big\| \Big( \sum_{k=1}^n
\lambda_k^{\frac{2}{p'}} \mu_k^{\frac{2}{p}} x_k^* x_k
\Big)^{\frac12} \Big\|_p \right\}.
\end{eqnarray*}
\end{itemize}
Moreover, let us write $\mathcal{G}_p$ for the closed subspace of
$L_p(\Gamma)$ generated by the system of generalized circular
variables $(g_{k,p})_{k \ge 1}$. Then, there exists a completely
bounded projection $\gamma_p: L_p(\Gamma) \to \mathcal{G}_p$
satisfying
 $$\|\gamma_p\|_{cb} \le 2^{|1-\frac2p|}.$$
\end{theorem}

\begin{remark}
\emph{It is worthy of mention that Theorem \ref{Th-OldE} improves
Theorem C in the case of generalized circular systems. Indeed, we
have only used two terms while Theorem C needs three terms in the
general case of $1$-homogeneous polynomials. This phenomenon will
also occur in the case of degree $2$, see below.}
\end{remark}

As application, we collect some interpolation identities that
arise from Theorem \ref{Th-OldE}. Indeed, we consider the spaces
$\mathcal{J}_p$ and $\mathcal{K}_p$, respectively defined as the
closure of finite sequences in $L_p(\mathcal{N})$ with respect to
the following norms
\begin{eqnarray*}
\|(x_k)\|_{\mathcal{K}_p} \!\!\!\! & = & \!\!\!\! \inf_{x_k =
a_k+b_k} \Big\| \Big( \sum_{k=1}^n \lambda_k^{\frac2p}
\mu_k^{\frac{2}{p'}} a_k a_k^* \Big)^{\frac12} \Big\|_p + \Big\|
\Big( \sum_{k=1}^n \lambda_k^{\frac{2}{p'}} \mu_k^{\frac2p} b_k^*
b_k \Big)^{\frac12} \Big\|_p, \\ \|(z_k)\|_{\mathcal{J}_p \,
\null} \!\!\!\! & = & \!\!\!\! \max \left\{ \Big\| \Big( \summ_k
\lambda_k^{\frac2p} \mu_k^{\frac{2}{p'}} z_k z_k^* \Big)^{\frac12}
\Big\|_p \, , \, \Big\| \Big( \summ_k \lambda_k^{\frac{2}{p'}}
\mu_k^{\frac2p} z_k^* z_k \Big)^{\frac12} \Big\|_p \right\}.
\end{eqnarray*}
 Given $1 \le p \le \infty$, we define the spaces
 $$L_p \big( \mathcal{N}; RC_p(\lambda,\mu) \big) = \begin{cases}
 \mathcal{K}_p & \mbox{for} \ 1 \le p \le 2, \\ \mathcal{J}_p &
 \mbox{for} \ 2 \le p \le \infty, \end{cases}$$
and the maps
 $$u_p : (x_k) \in L_p \big( \mathcal{N}; RC_p(\lambda,\mu) \big)
 \mapsto \summ_k x_k\ten g_{k,p} \in L_p(\mathcal{N} \bar\ten
 \Gamma).$$

\begin{corollary} \label{Corollary-Interpolation}
If $1 \le p_0,p_1 \le \infty$, $0<\theta<1$ and $1/p =
(1-\theta)/p_0 + \theta/p_1$, then
 $$\big[ L_{p_0} \big( \mathcal{N};
 RC_{p_0}(\lambda,\mu) \big), L_{p_1} \big( \mathcal{N};
 RC_{p_1}(\lambda,\mu) \big) \big]_\theta \simeq L_p \big(
 \mathcal{N}; RC_p(\lambda,\mu) \big).$$
Moreover, the relevant constants are majorized by a universal
constant.
\end{corollary}

\dem Let us recall Kosaki's theorem \cite{Ko}
 $$[L_{p_0}(\mathcal{N} \bar\ten \Gamma),\;
 L_{p_1}(\mathcal{N} \bar\ten\Gamma)]_{\theta}
 = L_p(\mathcal{N} \bar\ten \Gamma).$$
More precisely, if the von Neumann algebra $\mathcal{N}$ is
equipped with the \emph{n.f.} state $\psi$ and $d_{\psi \ten
\phi}$ denotes the density associated to $\psi \ten \phi$, we use
in the interpolation isometry above the symmetric inclusions
\begin{eqnarray*}
d_{\psi \ten \phi}^{1/2p_0'} L_{p_0}(\mathcal{N} \bar\ten \Gamma)
d_{\psi \ten \phi}^{1/2p_0'} & \subset & L_1(\mathcal{N} \bar\ten
\Gamma), \\ d_{\psi \ten \phi}^{1/2p_1'} L_{p_1}(\mathcal{N}
\bar\ten \Gamma) d_{\psi \ten \phi}^{1/2p_1'} & \subset &
L_1(\mathcal{N} \bar\ten \Gamma).
\end{eqnarray*}
Then we recall from Theorem \ref{Th-OldE} that the maps $u_p$
defined above are isomorphic embeddings. Using in addition the
projection $\gamma_p$ introduced in Theorem \ref{Th-OldE}, we
deduce the assertion. The proof is complete. \fin

Corollary \ref{Corollary-Interpolation} provides interesting
applications in the theory of operator spaces. Given two sequences
$(\xi_k)_{k \ge 1}$ and $(\rho_k)_{k \ge 1}$ of positive numbers
we introduce the operator space $R_p(\xi) \cap C_p(\rho)$ as the
span of the sequence $f_k = \xi_k e_{1k} + \rho_k e_{k1}$ in the
Schatten class $S_p$. Note that
\begin{eqnarray*}
\lefteqn{\Big\| \summ_k x_k \ten f_k \Big\|_{L_p(\mathcal{N};
R_p(\xi) \cap C_p(\rho))}} \\ & \sim & \max \left\{ \Big\| \Big(
\summ_k \xi_k^2 x_kx_k^* \Big)^{1/2} \Big\|_p \, , \, \Big\| \Big(
\summ_k \rho_k^2 x_k^*x_k \Big)^{1/2} \Big\|_p \right\}.
\end{eqnarray*}
By duality we understand the sum $R_p(\xi) + C_p(\rho)$ as a
quotient space. Indeed, we consider the subspace $R_p \oplus C_p$
of $S_p$ as the span of $(e_{1k}; e_{k1})_{k \ge 1}$ in $S_p$.
Then we have
$$R_p(\xi) + C_p(\rho) = R_p \oplus C_p /
\Delta(\xi,\rho),$$ where $\Delta$ is the weighted diagonal
$\Delta(\xi,\rho) = {\rm span} \big\{ \xi_k e_{1k} - \rho_k e_{k1}
\, | \ k \ge 1 \big\}$. Let $\pi$ be the natural quotient map and
let us consider the sequence $f_k = \pi(\xi_k e_{1k}) = \pi(\rho_k
e_{k1})$ in $R_p(\xi) + C_p(\rho)$. Then we find
\begin{eqnarray*}
\lefteqn{\Big\| \summ_k x_k \ten f_k \Big\|_{L_p(\mathcal{N};
R_p(\xi) + C_p(\rho))}} \\ & \sim & \inf_{x_k = a_k + b_k} \Big\|
\Big( \summ_k \xi_k^2 a_k a_k^* \Big)^{1/2} \Big\|_p + \Big\|
\Big( \summ_k \rho_k^2 b_k^* b_k \Big)^{1/2} \Big\|_p.
\end{eqnarray*}

\begin{corollary}
Let $(\lambda_k)_{k \ge 1}$ and $(\mu_k)_{k \ge 1}$ be two
sequences in $\R_+$ and $1 < p < \infty$. Then, the following
cb-isomorphisms hold according to the value of $\theta = 1/p$
$$\big[ R(\lambda) \cap
C(\mu), R(\lambda) + C(\mu) \big]_{\theta} \simeq_{cb}
\begin{cases} R_p(\lambda^{\theta} \mu^{1-\theta}) +
C_p(\lambda^{1-\theta} \mu^{\theta}), & \mbox{if} \ 1 < p \le 2,
\\ R_p(\lambda^{\theta} \mu^{1-\theta}) \cap
C_p(\lambda^{1-\theta} \mu^{\theta}), & \mbox{if} \ 2 \le p <
\infty. \end{cases}$$ The relevant constants are majorized by an
absolute constant.
\end{corollary}

\dem This is a reformulation of Corollary
\ref{Corollary-Interpolation} in operator space terms. \fin

We now discuss the analogue of Theorem \ref{Th-OldE} for
$q$-gaussians. We refer to \cite{BKS} for the basic definitions on
$q$-deformation and to Hiai's paper \cite{Hi} for the quasi-free
$q$-deformation. Given an infinite dimensional separable Hilbert
space $\mathcal{H}$ equipped with an orthonormal basis $(e_{\pm
k})_{k \ge 1}$ and given $-1<q<1$, we denote by
$\mathcal{F}_q(\mathcal{H})$ the associated $q$-Fock space
 $$\mathcal{F}_q(\mathcal{H})
 = \C \Omega \oplus \bigoplus_{n \ge 1} \mathcal{H}^{\ten n}$$
equipped with the $q$-scalar product induced by
 $$\big\langle f_1\ten \cdots \ten f_n,\; g_1 \ten \cdots \ten g_m
 \big\rangle_q =\delta_{nm} \sum_{\pi \in \mathcal{S}_n}^{\null}
 q^{i(\pi)}\langle f_1,\; g_{\pi(1)} \rangle \cdots \langle f_n,\;
 g_{\pi(n)}\rangle,$$
where $\mathcal{S}_n$ denotes the symmetric group of permutations
of $n$ elements and $i(\pi)$ stands for the number of inversions
of $\pi$. Given a vector $e \in \mathcal{H}$, we write $\ell_q(e)$
for the left creation operator and $\ell_q^*(e)$ for the left
annihilation, see \cite{BKS} for the precise definitions. As in
the free case we define
 $$gq_{k} = \lambda_k \ell_q(e_k) + \mu_k\ell_q^*(e_{-k})$$
after having fixed two sequences $(\lambda_k)_{k \ge 1}$ and
$(\mu_k)_{k \ge 1}$ of positive numbers. The $gq_k$'s are
\emph{$q$-generalized circular variables}. The von Neumann algebra
generated by these variables in the GNS-construction with respect
to the vacuum state $\phi_q(\cdot) = \langle \Omega, \cdot \,
\Omega \rangle_q$ will be denoted by $\Gamma_q$. A discussion of
the modular group of $\phi_q$ and important properties of these
von Neumann algebras can be found in Hiai's paper. Indeed, we
still have
\[ \sigma_t(gq_k)= (\lambda_k^{-1}\mu_k)^{2it} gq_k .\]
Therefore, $gq_k$ is an analytic element and we find as above
\begin{equation} \label{Eqgqkp}
gq_{k,p} = d_{\phi_q}^{\frac{1}{2p}}\, gq_k\,
d_{\phi_q}^{\frac{1}{2p}} = (\lambda_k^{-1} \mu_k)^{\frac1p}\,
gq_k\, d_{\phi_q}^{\frac{1}{p}} =  (\lambda_k
\mu_k^{-1})^{\frac1p}\, d_{\phi_q}^{\frac{1}{p}}\, gq_k.
\end{equation}

\vskip5pt

\demE Let us first see that the map
\begin{equation} \label{Eq-Boundup1}
u_p: (x_k) \in \mathcal{K}_p \mapsto \summ_k x_k \ten gq_{k,p} \in
L_p(\mathcal{N} \bar\ten \Gamma_q)
\end{equation}
is a contraction for $1 \le p \le 2$. According to \cite{JX}, we
have
 $$\|x\|_p \le \min \Big\{ \big\| \mathsf{E}(xx^*)^{\frac12}
 \big\|_p, \big\| \mathsf{E}(x^*x)^{\frac12} \big\|_p \Big\} \quad
 \mbox{for} \quad 1 \le p \le 2.$$
Taking $x = \sum_k x_k \ten gq_{k,p}$, the $L_p$-norm of $x$ is
bounded above by
 $$\min \left\{ \Big\| \Big( \summ_{i,j} \mathsf{E} \big( x_i x_j^*
 \ten gq_{i,p} gq_{j,p}^* \big) \Big)^{\frac12} \Big\|_p \, , \,
 \Big\| \Big( \summ_{i,j} \mathsf{E} \big( x_i^* x_j \ten
 gq_{i,p}^* gq_{j,p} \big) \Big)^{\frac12} \Big\|_p \right\},$$
where $\mathsf{E} = id_{\mathcal{N}} \ten \phi_q$ in our case.
Therefore, recalling from \eqref{Eqgqkp} that
 $$gq_{i,p}\,
 gq_{j,p}^* = d_{\phi_q}^{\frac{1}{2p}}\, gq_i\,
 d_{\phi_q}^{\frac{1}{p}}\, gq_j^* \, d_{\phi_q}^{\frac{1}{2p}} =
 (\lambda_i \mu_i^{-1} \lambda_j \mu_j^{-1})^{\frac1p}\,
 d_{\phi_q}^{\frac{1}{p}}\, gq_i\, gq_j^*\, d_{\phi_q}^{\frac{1}{p}},$$
 $$gq_{i,p}^*\, gq_{j,p} = d_{\phi_q}^{\frac{1}{2p}}\, gq_i^*\,
 d_{\phi_q}^{\frac{1}{p}}\, gq_j\, d_{\phi_q}^{\frac{1}{2p}} =
 (\lambda_i^{-1} \mu_i \lambda_j^{-1} \mu_j)^{\frac1p}\,
 d_{\phi_q}^{\frac{1}{p}}\, gq_i^*\, gq_j\, d_{\phi_q}^{\frac{1}{p}},$$
and using the identities $\phi_q(gq_i\, gq_j^*) = \delta_{ij}
\mu_i^2$ and $\phi_q(gq_i^*\,gq_j) = \delta_{ij} \lambda_i^2$, we
deduce
\begin{eqnarray*}
\mathsf{E} \big( x_i x_j^* \ten gq_{i,p}\, gq_{j,p}^* \big) & = &
\delta_{ij} \lambda_i^{\frac2p} \mu_i^{\frac{2}{p'}} x_ix_i^*, \\
\mathsf{E} \big( x_i^* x_j \ten gq_{i,p}^*\, gq_{j,p} \big) & = &
\delta_{ij} \lambda_i^{\frac{2}{p'}} \mu_i^{\frac{2}{p}} x_i^*x_i.
\end{eqnarray*}
Therefore, the triangle inequality yields
 $$\Big\| \sum_{k=1}^n x_k \otimes gq_{k,p} \Big\|_p \le \min
 \left\{ \Big\| \Big( \sum_{k=1}^n \lambda_k^{\frac2p}
 \mu_k^{\frac{2}{p'}} x_k x_k^* \Big)^{\frac12} \Big\|_p , \Big\|
 \Big( \sum_{k=1}^n \lambda_k^{\frac{2}{p'}} \mu_k^{\frac{2}{p}}
 x_k^* x_k \Big)^{\frac12} \Big\|_p \right\}.$$
This proves the contractivity of \eqref{Eq-Boundup1} for $1 \le p
\le 2$. Now we show that
\begin{equation} \label{Eq-Boundup2}
u_p: (x_k) \in \mathcal{J}_p \mapsto \summ_k
x_k \ten gq_{k,p} \in L_p(\mathcal{N} \bar\ten \Gamma_q)
\end{equation}
is bounded for $2 \le p \le \infty$ with a constant $c_q$
depending only on $q$. If $p=2$, the result follows by the
orthogonality of the $gq_{k,2}$'s in $L_2(\Gamma_q)$. Therefore,
according to Corollary \ref{Corollary-Interpolation}, it suffices
to estimate the norm of $u_{\infty}: \mathcal{J}_{\infty}\to
\mathcal{N} \bar\ten \Gamma_q$ and apply complex interpolation. By
the definition of $gq_k$ we have
 $$\summ_k x_k \ten gq_k = \summ_k
 \lambda_k x_k \ten \ell_{q}(e_k) + \summ_k \mu_k x_k \ten
 \ell_q^*(e_{-k}).$$
By Cauchy-Schwartz,
\begin{eqnarray*}
\Big\| \summ_k \lambda_k  x_k \ten \ell_{q}(e_k) \Big\|_\infty &
\le & \Big\| \Big( \summ_k \lambda_k^2 x_k^*x_k \Big)^{\frac12}
\Big\|_\infty \, \Big\| \Big( \summ_k \ell_{q}(e_k)
\ell_{q}(e_k)^* \Big)^{\frac12} \Big\|_\infty \\ & \le &
\frac{1}{\sqrt{1-|q|}} \, \Big\| \Big( \summ_k \lambda_k^2
x_k^*x_k \Big)^{\frac12} \Big\|_\infty,
\end{eqnarray*}
where the last inequality follows from \cite{BS}. Similarly, we
have
 $$\Big\| \summ_k \mu_k x_k \ten \ell_q^*(e_{-k})
 \Big\|_\infty \le \frac{1}{\sqrt{1-|q|}} \, \Big\| \Big( \summ_k
 \mu_k^2 x_k x_k^* \Big)^{\frac12} \Big\|_\infty.$$
Thus we obtain $\|u_{\infty}\| \le 2/\sqrt{1-|q|}$ and
 $$\big\|u_p: \mathcal{J}_p \to L_p(\mathcal{N} \bar\otimes \Gamma_q)
 \big\| \le \Big( \frac{2}{\sqrt{1-|q|}} \Big)^{1-\frac2p} \quad
 \mbox{for} \quad 2 \le p \le \infty.$$
The crucial observation here is that
\begin{eqnarray}
\label{EqDualforProj}
\lefteqn{\big\langle u_p((x_k)), u_{p'}((z_k)) \big\rangle} \\
\nonumber & = & \summ_{i,j} \mbox{tr}_{\mathcal{N}}(x_i^* z_j) \,
\mbox{tr}_{\Gamma_q} (gq_{i,p}^* gq_{j,p'}) \\ \nonumber & = &
\summ_{i,j} \mbox{tr}_{\mathcal{N}} (x_i^* z_j) \, (\lambda_i^{-1}
\mu_i)^{1/p} \, (\lambda_j^{-1} \mu_j)^{1/p'} \,
\mbox{tr}_{\Gamma_q} (d_\phi^{\frac1p} gq_i^* gq_j d_\phi^{\frac{1}{p'}}) \\
\nonumber & = & \summ_{i,j} \mbox{tr}_{\mathcal{N}} (x_i^* z_j) \,
(\lambda_i^{-1} \mu_i)^{1/p} \, (\lambda_j^{-1} \mu_j)^{1/p'} \,
\phi_q(gq_i^* gq_j) \\ \nonumber & = & \summ_k \ \lambda_k \mu_k
\, \mbox{tr}_{\mathcal{N}} (x_k^*y_k) = \big\langle (x_k), (z_k)
\big\rangle.
\end{eqnarray}
This relation and the boundedness of the maps \eqref{Eq-Boundup1}
and \eqref{Eq-Boundup2} immediately imply the inequalities stated
in i) and ii).

\vskip5pt

On the other hand, according to \eqref{EqDualforProj} we know that
$u_{p'}^* u_p$ is the identity map and we may construct the
following projection for every index $1 \le p \le \infty$
 $$u_p u_{p'}^* = id_{L_p(\mathcal{N})}
 \ten \gamma q_p: L_p(\mathcal{N} \bar\ten \Gamma_q) \to
 L_p(\mathcal{N}; \mathcal{G}q_p).$$
By elementary properties from \cite{P2} of vector-valued
noncommutative $L_p$ spaces, it suffices to prove that the maps
above are bounded with the following constants for $1 \le p \le 2
\le p' \le \infty$
 $$\max \Big\{ \|\gamma q_p\|_{cb}, \|\gamma
 q_{p'}\|_{cb} \Big\} = \max \Big\{ \big\| u_p u_{p'}^* \big\|,
 \big\| u_{p'}u_{p}^* \big\| \Big\} \le \Big(
 \frac{2}{\sqrt{1-|q|}} \Big)^{\frac{2}{p}-1}.$$
Recalling that the second estimate follows from the first by
taking adjoints and that the estimate for $p=2$ is trivial, it
suffices to prove the estimate for $u_1u_\infty^*$ and apply
complex interpolation. However, according to our previous
estimates we find $\|u_1 u_\infty^*\| \le 2/\sqrt{1-|q|}$, as
desired. This completes the proof. \fin

After this intermezzo on $q$-gaussians, we conclude by
illustrating our inequalities for $2$-homogeneous polynomials on
generalized circular variables. Again, our result in this
particular case improves Theorem C since we obtain only three
terms out of the five given there.

\vskip5pt

\demF Following the arguments in Theorem E, it suffices to prove
the assertion for $2 \le p \le \infty$ since the case $1 \le p \le
2$ and the complementation result follow from the same duality
arguments. In order to prove the assertion for $2 \le p \le
\infty$, we first consider a finite index set $\Lambda$ to
factorize
 $$\sum_{i \neq j} x_{ij} \ten d_\phi^{\frac{1}{2p}} g_i
 g_j d_\phi^{\frac{1}{2p}} = \sum_{i \neq j} \big( x_{ij} \ten
 d_\phi^{\frac{1}{2p}} g_i \big) \big( 1 \ten g_j d^{\frac{1}{2p}}
 \big) = \sum_{i \neq j} \alpha_{ij} \beta_j.$$
According to Theorem B we have for $2 \le p \le \infty$
 $$\Big\|\sum_{i \neq j} \alpha_{ij} \beta_j \Big\|_p
 \sim_c \Big\| \sum_{i\neq j} |\alpha_{ij}\big\rangle\beta_j\Big\|_p
 + \Big\|\sum_{i \neq j} \alpha_{ij} \big\langle \beta_j | \Big\|_p.$$
Let us denote the terms on the right by $\mathrm{A}$ and
$\mathrm{B}$ respectively. To simplify the expressions for
$\mathrm{A}$ and $\mathrm{B}$ we need to calculate
$\mathsf{E}(\alpha_{ij}^* \alpha_{kl})$ and $\mathsf{E}(\beta_j
\beta_l^*)$. According to \eqref{Eq-Generalized-p}, we easily find
\begin{eqnarray*}
\mathsf{E}(\alpha_{ij}^* \alpha_{kl}) & = & \delta_{ik} \, \Big(
x_{ij}^* x_{il} \ten \lambda_{i}^{\frac{2}{p'}} \mu_{i}^{\frac2p}
d_\phi^{\frac{1}{p}} \Big), \\ \mathsf{E} (b_j b_l^*) & = &
\delta_{jl} \ \Big( 1 \ten \lambda_{j}^{\frac{2}{p}}
\mu_{j}^{\frac{2}{p'}} d_\phi^{\frac{1}{p}} \Big).
\end{eqnarray*}
Using these relations and recalling the factorization above of
$x_{ij}$, we obtain
\begin{eqnarray*}
\mathrm{A} & = & \Big\| \Big( \summ_i \lambda_i^{\frac{2}{p'}}
\mu_i^{\frac2p} \summ_{j_1,j_2} x_{ij_1}^* x_{ij_2} \ten
g_{j_1,p}^* g_{j_2,p} \Big)^{\frac12} \Big\|_p, \\
\mathrm{B} & = & \Big\| \Big( \summ_j \lambda_j^{\frac{2}{p}}
\mu_j^{\frac{2}{p'}} \summ_{i_1,i_2} x_{i_1j} x_{i_2j}^* \ten
g_{i_1,p} g_{i_2,p}^* \Big)^{\frac12} \Big\|_p.
\end{eqnarray*}
Equivalently, we have
\begin{eqnarray*}
\mathrm{A} & = & \Big\| \summ_i \lambda_i^{\frac{1}{p'}}
\mu_i^{\frac1p} \Big( \summ_j x_{ij} \ten g_{j,p} \Big) \ten
e_{i1} \Big\|_p = \Big\| \summ_k a_k \ten g_{k,p}
\Big\|_p, \\
\mathrm{B} & = & \Big\| \summ_j \lambda_j^{\frac{1}{p}}
\mu_j^{\frac{1}{p'}} \Big( \summ_i x_{ij} \ten g_{i,p} \Big) \ten
e_{1j} \Big\|_p = \Big\| \summ_k b_k \ten g_{k,p} \Big\|_p,
\end{eqnarray*}
where $a_k$ and $b_k$ are respectively given by $$a_k = \summ_i
\lambda_i^{\frac{1}{p'}} \mu_i^{\frac1p} x_{ik} \ten e_{i1} \quad
\mbox{and} \quad b_k = \summ_j \lambda_j^{\frac{1}{p}}
\mu_j^{\frac{1}{p'}} x_{kj} \ten e_{1j}.$$ According to Theorem
\ref{Th-OldE} we obtain
\begin{eqnarray*}
\mathrm{A} & \sim_c & \Big\| \Big( \sum_{k=1}^n
\lambda_k^{\frac{2}{p}} \mu_k^{\frac{2}{p'}} a_k a_k^*
\Big)^{\frac12} \Big\|_p + \Big\| \Big( \sum_{k=1}^n
\lambda_k^{\frac{2}{p'}} \mu_k^{\frac{2}{p}} a_k^* a_k
\Big)^{\frac12}
\Big\|_p = \mathrm{A}_1 + \mathrm{A}_2, \\
\mathrm{B} & \sim_c & \Big\| \Big( \sum_{k=1}^n
\lambda_k^{\frac{2}{p}} \mu_k^{\frac{2}{p'}} b_k b_k^* \,
\Big)^{\frac12} \Big\|_p + \Big\| \Big( \sum_{k=1}^n
\lambda_k^{\frac{2}{p'}} \mu_k^{\frac{2}{p}} b_k^* b_k \,
\Big)^{\frac12} \Big\|_p\,  = \mathrm{B}_1 + \mathrm{B}_2.
\end{eqnarray*}
Finally, using the terminology introduced in the statement of
Theorem F, we have
$$\mathrm{B}_1 = \mathcal{R}_p(x), \quad \mathrm{A}_1 =
\mathcal{M}_p(x) = \mathrm{B}_2, \quad \mathrm{A}_2 =
\mathcal{C}_p(x).$$ Details of the identities above are left to
the reader. This completes the proof. \fin

\begin{remark}
\emph{Although it is out of the scope of this paper, the methods
used in the proof of Theorem F are also valid for any degree $d
\ge 1$. In this way, the $L_p$ norm of an operator-valued
$d$-homogeneous polynomial on generalized circular random
variables behaves as the \emph{asymmetric} version of the main
result in \cite{PP}. Thus, we obtain $d+1$ terms instead of the
$2d+1$ provided by Theorem C.}
\end{remark}

\bibliographystyle{amsplain}

\vskip-2pt

\end{document}